\documentclass[journal]{IEEEtran}
\IEEEoverridecommandlockouts
\usepackage{cite}
\usepackage{amsmath,amssymb,amsfonts}
\AtBeginDocument{\pagenumbering{arabic}}
\usepackage{nccmath}
\usepackage{romannum}
\usepackage{graphicx}
\usepackage{mathtools}
\usepackage[percent]{overpic}  
\usepackage{amsmath}    
\usepackage{amssymb}    
\usepackage{bm}         
\usepackage{mathtools}  
\usepackage{textcomp}
\usepackage{svg}
\usepackage{xcolor}
\usepackage{mdwlist}
\usepackage{bm}
\usepackage{bbm}
\usepackage{eurosym}
\usepackage{tikz}
\usepackage{booktabs}  
\usepackage{array}     
\usepackage{dblfloatfix}
\usepackage{makecell}
\usepackage{comment}
\usepackage{amsmath}
\usepackage{tabularx}
\usepackage{multirow}
\usepackage{enumitem}
\usepackage{graphicx}
\usepackage{bbm}
\usepackage{booktabs}
\usepackage{tabularx}
\usepackage{threeparttable}
\newcolumntype{C}{>{\centering\arraybackslash}X} 
\newtheorem{theorem}{Theorem}

\newtheorem{proposition}{Proposition}
\newtheorem{corollary}{Corollary}

\newtheorem{lemma}{Lemma}
\usepackage{algorithm}
\usepackage{array}
\usepackage{algpseudocode}
\usepackage{tabularx,booktabs}
\usepackage{textcomp}
\usepackage{microtype}
\usepackage{nomencl}
\usepackage{etoolbox}

\renewcommand{\nomgroup}[1]{%
\ifthenelse{\equal{#1}{S}}{\item[\textbf{Sets}]}{%
\ifthenelse{\equal{#1}{P}}{\item[\textbf{Parameters}]}{%
\ifthenelse{\equal{#1}{V}}{\item[\textbf{Variables}]}{%
\ifthenelse{\equal{#1}{R}}{\item[\textbf{Random Variables}]}{}}}}
}

\makenomenclature
\usepackage{booktabs,makecell,tabularx}

\usepackage{tikz}
\tikzset{
    block/.style = {rectangle, draw, text width=4cm, text centered, minimum height=1.5cm, node distance=1.5cm},
    smallblock/.style = {rectangle, draw, text width=3.2cm, text centered, minimum height=1cm, node distance=1.5cm},
    line/.style = {draw, -latex'},
}

\usepackage{pgfplots}
\usetikzlibrary{intersections}
\usepgfplotslibrary{fillbetween}
\usetikzlibrary{arrows, arrows.meta, positioning, backgrounds, automata, decorations, shapes, fit, calc}

\mathchardef\mhyphen="2D 

\usepackage{siunitx}
\usepackage{etoolbox}
\newrobustcmd{\B}{\bfseries}
\allowdisplaybreaks
\usetikzlibrary{positioning}
\usetikzlibrary{decorations.pathreplacing}
\usepackage{amssymb}
\usepackage{pifont}
\newcommand{\xmark}{\ding{55}}%
\usetikzlibrary{calc}
\usepackage[americanresistors,americaninductors,americancurrent,RPvoltages]{circuitikz}
\def\BibTeX{{\rm B\kern-.05em{\sc i\kern-.025em b}\kern-.08em
    T\kern-.1667em\lower.7ex\hbox{E}\kern-.125emX}}
\ifCLASSOPTIONcompsoc
    \usepackage[caption=false, font=normalsize, labelfont=sf, textfont=sf]{subfig}
\else
\usepackage[caption=false, font=footnotesize]{subfig}
\fi
\algnewcommand\algorithmicinput{\textbf{Input:}}
\algnewcommand\INPUT{\item[\algorithmicinput]}
\algdef{SE}[IF]{If}{EndIf}[1]
  {\algorithmicif\ #1\ \algorithmicthen}{\algorithmicend\ \algorithmicif}

\usepackage{hyperref}

\PassOptionsToPackage{hyphens}{url}\usepackage{hyperref}
\definecolor{lime}{HTML}{A6CE39}
\DeclareRobustCommand{\orcidicon}{%
	\begin{tikzpicture}
	\draw[lime, fill=lime] (0,0) 
	circle [radius=0.16] 
	node[white] {{\fontfamily{qag}\selectfont \tiny ID}};
	\draw[white, fill=white] (-0.0625,0.095) 
	circle [radius=0.007];
	\end{tikzpicture}
	\hspace{-2mm}
}

\foreach \x in {A, ..., Z}{%
	\expandafter\xdef\csname orcid\x\endcsname{\noexpand\href{https://orcid.org/\csname orcidauthor\x\endcsname}{\noexpand\orcidicon}}
}

\ifCLASSINFOpdf

\else

\fi

\begin{document}

\title{Bilevel Transmission Expansion Planning with Joint Chance-Constrained Dispatch}

\author{
        Yuxin~Xia\orcidA{},~\IEEEmembership{Student,~IEEE,}
        Yihong~Zhou\orcidD{},~\IEEEmembership{Student,~IEEE,}
        Iacopo~Savelli\orcidB{},
        and Thomas~Morstyn\orcidC{},~\IEEEmembership{Senior,~IEEE} 
\thanks{Y. Xia and Y. Zhou are with the School
of Engineering, the University of Edinburgh, Edinburgh, EH9 3JL, U.K. (e-mail:
yuxin.xia@ed.ac.uk, yihong.zhou@ed.ac.uk).
I. Savelli is with the Centre for Research on Geography, Resources, Environment, Energy $\&$ Networks, Bocconi University, Milano, Italy (e-mail: iacopo.savelli@unibocconi.it). T. Morstyn is with Department of Engineering Science, University of Oxford, U.K. (e-mail:
thomas.morstyn@eng.ox.ac.uk).}
\thanks{Manuscript received XXXX; revised XXXX.}}


\maketitle

\begin{abstract}
In transmission expansion planning (TEP), network planners make long-term investment decisions while anticipating market clearing outcomes that are increasingly affected by renewable generation uncertainty. Additionally, market participants' sensitivity to network charges and the requirement for cost recovery by the network planner introduce further complexity. Since the day-ahead market clears before uncertainty realizes, explicitly modelling these uncertainties at the lower-level market clearing becomes important in bilevel TEP problems. In this paper, we introduce a novel bilevel TEP framework with lower-level joint chance-constrained market clearing that manages line flow constraints under wind uncertainty and accounts for the effect of network tariffs on participants' actual marginal costs and utility. To solve this complex problem, we propose a Strengthened Linear Approximation (SLA) technique for handling Wasserstein distributionally robust joint chance constraints with right-hand-side uncertainties (RHS-WDRJCC). The proposed method offers more efficient approximations without additional conservativeness and avoids the numerical issues encountered in existing approaches by introducing valid inequalities. The case study demonstrates that the proposed model achieves the desired out-of-sample constraint satisfaction probability. Moreover, the numerical results highlight the significant computational advantage of SLA, achieving up to a 26× speedup compared to existing methods such as worst-case conditional value-at-risk, while maintaining high solution quality.

\looseness = -1
\end{abstract}

\begin{IEEEkeywords}
Transmission expansion planning (TEP), bilevel optimization, Wasserstein distributionally robust joint chance constraint (WDRJCC), Strengthened Linear Approximation (SLA), out-of-sample analysis.
\end{IEEEkeywords}

\IEEEpeerreviewmaketitle

\section{Introduction}
Market-driven transmission expansion planning (TEP) is attractive for reducing the excess costs associated with congestion \cite{6345475}, yet it faces increasing challenges due to several evolving factors. Firstly, renewable energy sources introduces substantial uncertainty into electricity market operations \cite{9242289}, requiring TEP methodologies to explicitly account for these uncertainties and perform robust grid investments. Secondly, the impact of network tariffs on the market behaviour of price-sensitive producers and consumers is often neglected, which can lead to inefficiencies \cite{SAVELLI2020113979}. Lastly, given the time-consuming nature and practical challenges associated with constructing new transmission lines, a comprehensive approach to grid expansion should incorporate alternatives such as reconductoring \cite{2411207121}. Consequently, an effective TEP framework should address renewable generation uncertainties and account for the impact of price responsiveness, enabling timely, cost-effective, and robust grid expansions to meet evolving power system needs.

Power system operation involves critical constraints that must be satisfied under uncertainty to ensure system safety and reliability. To tackle uncertainty in constraints, formulations use non-deterministic approaches such as robust optimization (RO) \cite{8275041} and chance-constrained programming (CCP) \cite{10607938}. CCP is less restrictive than RO by allowing a controlled probability of constraint violation, which can yield less conservative solutions. However, classic CCP requires knowledge of the exact distribution of uncertainty \cite{zhou2024strengthenedfasterlinearapproximation}, which is typically unavailable, as decision-makers usually only have access to historical data. To address this limitation, distributionally robust chance-constrained programming (DRCCP) has been proposed. DRCCP is a specific form of distributionally robust optimization (DRO) applied to chance constraints. It hedges against distributional uncertainty by considering an ``ambiguity set'', which represents a family of potential distributions of the uncertainty. These sets are typically classified as either moment-based, which can be defined by mean and variance \cite{Zymler2013}, or distance-based, such as those based on $\phi$-divergence or Wasserstein distance \cite{chen2023approximations, mohajerin2018data}. 

DRO, including DRCCP formulations, has been applied in many power system planning problems. A comparison of these works is presented in Table~\ref{tab:comparison_works}, where most studies utilize moment-based ambiguity set due to the existence of tractable solution techniques \cite{ZHAN2022107417, 8933104, 10607938, 8447238, 9028125, 8442871}. However, moment-based ambiguity sets provide a limited representation of the true probability distribution \cite{8294298}, often leading to poor out-of-sample performance or excessive conservativeness. In contrast, Wasserstein ambiguity sets have gained popularity due to their superior out-of-sample performance \cite{mohajerin2018data}. While Wasserstein distributionally robust chance-constrained (WDRCC) programming has been extensively used in power system operations \cite{9026959, zhou2024strengthenedfasterlinearapproximation}, only one paper \cite{9961917} constructs Wasserstein ambiguity sets for planning problems, solving them using conditional value-at-risk (CVaR) approximation schemes proposed in \cite{Xie2021}. Moreover, existing TEP literature employs distributionally robust individual chance-constrained (DRICC) programming which fails to provide the stronger guarantees offered by joint chance constraints (JCC) for simultaneous satisfaction of multiple system constraints \cite{10607938}.

\begin{table}[t]
\centering
\caption{Comparison of Relevant Literature and Proposed Model}
\label{tab:comparison_works}
\footnotesize 
\begin{tabularx}{\linewidth}{l *{5}{X}} 
\toprule\midrule
\shortstack[c]{\\Ref.\\~} & \shortstack[c]{\\Problem\\~} & \shortstack[c]{\\Uncertainty\\~} & \shortstack{Ambiguity\\Set} & \shortstack{Tariff\\Response} & \shortstack{Cost\\Recovery} \\
\midrule
\cite{ZHAN2022107417} & IESP & DRICC & Moment & N/A & N/A \\
\cite{8933104} & GEP & DRICC & Moment & N/A & N/A \\
\cite{9961917} & GEP & DRICC & Wasserstein & N/A & N/A \\
\cite{8447238,9028125,8442871} & TEP & DRO & Moment & N & N \\
\cite{10607938} & TEP & DRICC & Moment & N & N \\
This paper & TEP & DRJCC & Wasserstein & Y & Y \\
\midrule\bottomrule
\end{tabularx}
\vspace{0.5mm}
\begin{minipage}{\linewidth}
\scriptsize 
IESP: Integrated energy system planning; GEP: Generation expansion planning; TEP: Transmission expansion planning; DRO: Distributionally robust optimization; DRICC: Distributionally robust individual chance constraint; DRJCC: Distributionally robust joint chance constraint.
\end{minipage}
\end{table}

In addition, to the best of the authors' knowledge, existing distributionally robust TEP studies have not addressed the revenue adequacy or cost recovery of the network planner. Addressing these aspects typically requires hierarchical decision-making structures (e.g., bilevel or multilevel optimization \cite{XIA2025124721}) that allow the planner to anticipate \emph{lower-level} market responses and access the dual variables of the power balance constraints, which represent market prices under marginal pricing \cite{8031353}. These dual variables are important for calculating potential revenues from congestion rents/merchandising surplus. However, the combination of hierarchical bilevel or multilevel optimization and DRO presents several open research challenges \cite{BECK2023401}. A key challenge arises because the follower often makes decisions before the uncertainty is realized. For example, the day-ahead market clears before uncertainties such as wind generation or load are realized. The method in \cite{10038580} assumes uncertainty at the leader stage, but in many power system applications leaders first set investment or bidding strategies and then the day-ahead market clears before uncertainty realizes. CCP has been intensively shown to internalize the stochasticity of renewable resources in market clearing, leading to formulations that enforce network and generation limits under uncertainty \cite{8937824, zhou2024strengthenedfasterlinearapproximation}. This makes CCP-based uncertainty modelling at the \emph{lower-level} in bilevel optimization necessary for effective risk management; however, research on this subject remains under-explored \cite{Heitsch2022,BECK2023401}. Moreover, standard solution methods for bilevel problems require the lower‐level problem to be convex to derive optimality conditions, which becomes problematic when CCP or DRCCP appears at the lower level.

Our prior work \cite{zhou2024strengthenedfasterlinearapproximation} presents an initial contribution to integrating bilevel optimization with DRCCPs at the lower level by introducing the Strengthened and Faster Linear Approximation (SFLA) for Wasserstein distributionally robust joint chance-constrained problems with right-hand-side uncertainty (RHS-WDRJCC). We demonstrate the effectiveness of SFLA by applying it to a bilevel storage bidding problem. In this setting, the lower-level involves RHS-WDRJCC market clearing under wind power uncertainty, where the coefficients of the random variables are fixed parameters, independent of upper-level decisions. However, when the coefficients of random variables depend on upper-level decision, existing exact and approximation schemes face significant limitations. Specifically, the Bonferroni approximation~\cite{chen2023approximations} and exact formulations \cite{2022_OR_exact_DRO,exact_milp_strengthened} are not applicable; the Linear Approximation (LA) \cite{chen2023approximations} and worst-case Conditional Value-at-Risk (W-CVaR) \cite{mohajerin2018data, chen2023approximations} are computationally intensive; and SFLA~\cite{zhou2024strengthenedfasterlinearapproximation} may introduce numerical issues, see Section~\ref{section:discussion_appro} for discussion. This highlights the need for computationally efficient approximation schemes for RHS-WDRJCC that avoid additional numerical difficulties and conservativeness for this case.

The central research question addressed in this paper is: \emph{How can TEP integrate distributionally robust joint chance‐constrained dispatch while explicitly modelling market participants’ price‐sensitive behaviour and ensuring revenue adequacy for the network planner?} To address this, we propose a novel bilevel TEP model that integrates dynamic network charges with a lower-level WDRJCC dispatch, solved efficiently using the proposed Strengthened Linear Approximation (SLA) scheme. To the best of the authors' knowledge, this study represents the first effort to develop a bilevel TEP problem within the context of chance-constrained dispatch, explicitly accounting for the dependency of random variables' coefficients on upper-level decision variables in RHS-WDRJCC problems. In summary, the contributions of this paper are as follows:
\begin{itemize}
    \item We propose a bilevel TEP framework with joint chance-constrained market clearing at the lower level that captures wind power uncertainty in the dispatch problem.
    \item We propose a SLA scheme for RHS-WDRJCC  that overcomes limitations of existing approaches when random variable coefficients depend on upper-level decisions. This novel approximation achieves superior computational efficiency without sacrificing conservativeness through valid inequalities.
    \item Numerical results demonstrate the robustness of our approach by verifying desired out-of-sample constraint satisfaction probabilities. Comparisons with benchmarks including existing linear approximations, LA and W-CVaR, and comprehensive analysis demonstrate the superior performance of the proposed SLA method.
\end{itemize}

The rest of the paper is organized as follows: Section~\ref{section2} presents the formulation of the bilevel TEP problem with lower-level WDRJCC market clearing problems. Section~\ref{section3} presents existing solution methodologies for RHS-WDRJCC problems, introduces the proposed SLA method, and discusses the reformulation of the bilevel problem under different approximation schemes. Section~\ref{section4} presents the case study results. Finally, Section~\ref{section5} summarizes the key findings and concludes the paper. The full nomenclature and appendices are provided in the Supplementary Material.

\section{Problem Formulation}\label{section2}
\begin{figure}[t] \label{Fig:bilevel_frame}
  \centering
  \includegraphics[width=0.7\linewidth]{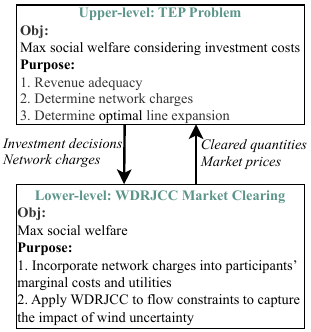}
  \caption{ Bilevel structure of the proposed formulation for optimal network planning. The upper-level problem represents the network planner, who determines investment decisions on reconductoring lines ($b^R_{t,l,j}$, $z^R_{t,l}$), parallel lines ($z^P_{t,l,m}$), and sets network charges—both volumetric-based ($\tau^V_l$) and capacity-based ($\tau^C$). Investment decisions affect the transmission capacity and the system susceptance matrix as the network topology changes. The volumetric-based charges $\tau^V_l$ also influence market participants’ marginal costs and willingness to pay, thereby affecting the outcomes of the lower-level market clearing problem. The lower-level problem is formulated as a RHS-WDRJCC problem to account for wind uncertainty. It determines cleared quantities ($g_{t,s,k,b}$, $d_{t,s,k,b}$, $p^{w,sch}_{t,s,k,b}$) and market prices ($\pi_{t,s,b}$). These outcomes are then used by the network planner to compute revenues, which consist of merchandising surplus and network charges, and are compared against investment costs to ensure revenue adequacy for the network planner.}
\end{figure}

This section presents the proposed bilevel RHS-WDRJCC TEP problem. The proposed framework is structured as a bilevel programming problem and is sketched in Fig.~\ref{Fig:bilevel_frame}. The upper level addresses the long-term network planning, whereas the lower level focuses on the joint chance constrained day-ahead market clearing. 

Specifically, the upper-level investment problem aims to maximize social welfare (see the definition of social welfare in \cite{SAVELLI2025101937}) by determining optimal investment strategies for transmission lines, alongside setting network charges, while ensuring revenue adequacy for the network planner. The planner can invest in expanding the transmission network by adding new candidate lines to existing corridors and by enhancing the capacity of existing lines through reconductoring \cite{2411207121}. Furthermore, the planner can levy both volumetric and capacity-based network charges. Notably, the imposition of volumetric network charges directly influences the marginal costs and the willingness-to-pay of price-sensitive market participants, consequently affecting the outcomes of the market clearing process \cite{SAVELLI2020113979}. The network planner recovers investment costs through a combination of merchandising surplus based on market clearing results as well as the volumetric and capacity-based network charges collected from market participants. The lower-level problem represents the day-ahead market clearing and explicitly captures how market participants respond to volumetric network charges set at the upper level. Moreover, to effectively manage the uncertainty associated with non-dispatchable wind power, joint chance constraints are employed to secure the thermal limits of transmission lines against uncertain power flows due to wind forecast errors. Further assumptions regarding the model can be found in Appendix~\ref{section:assumption}.

\subsection{Upper-level Problem}\label{section2c}
Problem~\eqref{eq:upper} specifies the long term transmission investment planning problem. The primal variables of the upper-level problem are \vspace{-5pt}$$\Xi_{UL} = \{b_{t,l, j}^{R},z^R_{t,l},o_{t,c},z^P_{t,l},\tau^C,\tau^V_l,\hat{c}_{t,s, k,b}^{d},\hat{c}_{t,s, k,b}^{g},\hat{c}_{t,s, k,b}^{w}\},$$
as well as all lower-level primal and dual variables, which will be defined later.

\begin{fleqn}
\begin{subequations}\label{eq:upper}
\vspace{-0.7em}
\subsubsection{Objective function}
\begin{align}
    &\max_{\Xi_{UL}} \sum_{t \in \mathcal{T}} \frac{1}{(1+r)^{t-1}}\bigg[\sum_{s \in \mathcal{S}}\sum_{b \in \mathcal{B}}\Psi\bigg(\sum_{k \in \Omega_{t, s,b}^{D}} c_{t,s, k,b}^{d} d_{t,s, k,b}\label{eq:upper_obj}\\[-5pt]
    &-\sum_{k \in \Omega_{t, s,b}^{G}} c_{t,s,k,b}^{g} g_{t,s,k,b} -\sum_{k \in \Omega_{t, s,b}^{W}} c^{cur}_{t,s,k,b}p^{w,cur}_{t,s,k,b}\bigg)-C_t \bigg] \notag \end{align}
The upper-level objective function~\eqref{eq:upper_obj} aims to maximize the present value of long-term social welfare considering investment costs $C_t$ and wind curtailment costs $c^{cur}_{t,s,k,b}p^{w,cur}_{t,s,k,b}$, using discount rate $r$. Volumetric ($\tau^V_l$) and capacity-based ($\tau^C$) network charges are not included in the objective function. This is because these charges represent monetary transfers from consumers and producers to the network operator for the purpose of recovering investment costs, thus constituting offsetting financial flows within the overall system. The parameter \(\Psi\) serves as a scaling factor to align costs and benefits across investment and operational scales\footnote{For example, if for each one-year investment period, we consider a single representative one-hour market‐clearing period (\(\mathcal{S} = \{1\}\)), then \(\Psi = 8760\).}.
\vspace{-0.7em}
\subsubsection{Actual marginal cost and benefit under volumetric tariff constraints}
    \begin{align}
&\hat{c}_{t,s,k,b}^{d} = c_{t,s,k,b}^{d} - (\sum_{l\in\mathcal{L}^R}z^R_{t,l}\tau^V_{l}\delta_{l,b} + \sum_{l\in\mathcal{L}^P}\sum_{m \in \mathcal{M}}z^P_{t,l,m}\tau^V_{l}\delta_{l,b} ), \notag\\[-5pt]
& \quad \forall t\in \mathcal{T} ,s \in \mathcal{S},k \in \Omega_{t, s,b}^{D},b \in \mathcal{B}, \label{eq:upperb}
\end{align}
\vspace{-22pt}
\begin{align}
&\hat{c}_{t,s,k,b}^{g} = c_{t,s,k,b}^{g} + (\sum_{l\in\mathcal{L}^R}z^R_{t,l}\tau^V_{l}\delta_{l,b} + \sum_{l\in\mathcal{L}^P}\sum_{m \in \mathcal{M}} z^P_{t,l,m}\tau^V_{l}\delta_{l,b}), \notag\\[-5pt]
&\quad \forall t\in \mathcal{T} ,s \in \mathcal{S},k \in \Omega_{t, s,b}^{G},b \in \mathcal{B},  \label{eq:upperc}
\end{align}
\vspace{-22pt}
\begin{align}
&\hat{c}_{t,s,k,b}^{w} = (\sum_{l\in\mathcal{L}^R}z^R_{t,l}\tau^V_{l}\delta_{l,b} + \sum_{l\in\mathcal{L}^P}\sum_{m \in \mathcal{M}}z^P_{t,l,m}\tau^V_{l}\delta_{l,b}),\notag\\[-5pt]
& \quad \forall t\in \mathcal{T} ,s \in \mathcal{S},k \in \Omega_{t, s,b}^{W},b \in \mathcal{B}, \label{eq:upperd}
\end{align}
Constraints~\eqref{eq:upperb}--\eqref{eq:upperd} introduce the effect of volumetric network charges into the actual marginal costs of generators and actual marginal benefits of consumers. For example, the actual marginal benefit of the consumer $k$ at time $t$ and scenario $s$ (i.e., $\hat{c}_{t,s,k,b}^{d}$) is reduced by a quantity (i.e. the imposed volumetric network charges $\tau^V_l$) with respect to its initial gross benefit (i.e., $c_{t,s,k,b}^{d}$), as shown in \eqref{eq:upperb}. Thus, the market demand curve is reshaped to account for the additional costs. A similar formulation is used for conventional generators in \eqref{eq:upperc} and wind generators in \eqref{eq:upperd}, where appears with a plus sign since costs (rather than benefits) are considered. The actual market curves determined after volumetric network costs will be used by the lower level problem~\eqref{eq:lower} to clear the market.\vspace{-0.7em}
\subsubsection{Merchandising surplus constraints}
    \begin{align}
    &MS_t=\Psi\sum_{s \in \mathcal{S}} \sum_{b \in \mathcal{B}}(\sum_{k \in \Omega_{t, s,b}^{D}}\pi_{t,s,b}^{*}d_{t, s,k,b}^{*} - \sum_{k \in \Omega_{t, s,b}^{G}}
 \pi_{t,s,b}^{*}g_{t, s,k,b}^{*} \notag\\[-5pt]
&- \sum_{k \in \Omega_{t, s,b}^{W}} \pi_{t,s,b}^{*}p_{t, s,k,b}^{w,sch*}),\quad\forall t\in\mathcal{T},\label{eq:uppere} \end{align}
Here, $MS_t$ in~\eqref{eq:uppere} represents the merchandising surplus at time $t$, calculated as the difference between total payments by demand and total revenues received by generators and wind producers, which is calculated using optimal lower-level market dispatch quantities ($g_{t,s,k,b}$, $d_{t,s,k,b}$, $p^{w,sch}_{t,s,k,b}$) and market prices ($\pi_{t,s,b}$). \vspace{-0.7em}
\subsubsection{Investment costs constraints}
\begin{align}
&C_t = \sum_{l \in \mathcal{L}^{R}}
      \sum_{j \in \mathcal{J}\setminus\{0\}}
      \bigl(K_{l}^{fix}\,b_{t,l,j}^{R}
           +K_{l}^{var}\,j\,b_{t,l,j}^{R}\,\mathcal{F}^{0}_{l}\bigr)\notag\\
    &+ \sum_{l \in \mathcal{L}^P}
        K_l^{fix}\Bigl(
          \sum_{m \in \mathcal{M}} m\,z^P_{t,l,m}
         - \mathbf{1}_{t>1}\sum_{m \in \mathcal{M}} m\,z^P_{t-1,l,m}
        \Bigr),\label{eq:upperf}\\
        &
    \quad \forall\,t \in \mathcal{T},\notag
\end{align}
Total investment costs are modelled in constraints~\eqref{eq:upperf}. For reconductoring existing lines $\mathcal{L}^R$, these costs include a fixed component $K_{l}^{f i x}$ and a variable component $K_{l}^{var}$, where the capacity expansion factor $j \in \mathcal{J}$ represents discrete increments relative to the original line rating $\mathcal{F}^{0}_{l}$. A binary variable \(b_{t, l, j}^{R}\) is equal to one if a lumpy reconductoring investment corresponding to an additional capacity of \(j \mathcal{F}_{l}^0\) is made for line \(l\) in investment period \(t\), and zero otherwise. Conversely, a fixed investment cost $K_{l}^{f i x}$ is assumed for candidate parallel lines $\mathcal{L}^P$, reflecting their typically fixed capacity. The expression $\sum_{m \in \mathcal{M}} m\,z^P_{t,l,m}- \mathbf{1}_{t>1}\sum_{m \in \mathcal{M}} m\,z^P_{t-1,l,m}$ calculates the number of new circuits installed on corridor \(l\) in period \(t\). The first term gives the total number of circuits installed by the end of period \(t\), while the second term subtracts those already installed by the end of period \(t{-}1\). The indicator \(\mathbf{1}_{t>1}\) ensures that this subtraction is skipped in the first investment period.
\vspace{-0.7em}
\subsubsection{Network charge revenue constraints}
    \begin{align}
       &VC^{R}_{t,l} = \Psi\sum_{s \in \mathcal{S}}\sum_{b \in \mathcal{B}}z^R_{t,l}\tau^V_{l}\delta_{l,b}(\sum_{k \in \Omega_{t, s,b}^{G}}g^*_{t, s,k,b}+\notag\\[-5pt]
&\sum_{k \in \Omega_{t, s,b}^{D}}d^*_{t, s,k,b}+\sum_{k \in \Omega_{t, s,b}^{W}}p_{t, s,k,b}^{w,sch*}), \quad \forall  t\in\mathcal{T},l\in \mathcal{L}^R,\label{eq:upperg}\end{align}
    \vspace{-15pt}
    \begin{align}
       &VC^{P}_{t,l} = \Psi\sum_{s \in \mathcal{S}}\sum_{b \in \mathcal{B}}\sum_{m \in \mathcal{M}}z^P_{t,l,m}\tau^V_{l}\delta_{l,b} (\sum_{k \in \Omega_{t, s,b}^{G}}g^*_{t, s,k,b}+\notag\\[-5pt]
&\sum_{k \in \Omega_{t, s,b}^{D}}d^*_{t, s,k,b}+\sum_{k \in \Omega_{t, s,b}^{W}}p_{t, s,k,b}^{w,sch*}), \quad \forall  t\in\mathcal{T} ,l \in \mathcal{L}^P,\label{eq:upperh}\end{align}
    \vspace{-15pt}
    \begin{align}
       VC_{t} = \sum_{l\in\mathcal{L}^R} VC^{R}_{t,l} +\sum_{l\in\mathcal{L}^P}VC^{P}_{t,l},\quad \forall  t\in\mathcal{T},\label{eq:upperi}\end{align}
    \vspace{-15pt}
    \begin{align}
       &CC_{t} = \Psi\sum_{s \in \mathcal{S}}\sum_{b \in \mathcal{B}}\tau^C(\sum_{k \in \Omega_{t, s,b}^{G}}g_{t, s,k,b}^{\max}+\sum_{k \in \Omega_{t, s,b}^{D}}d_{t, s,k,b}^{\max}\notag\\[-5pt]
&+\sum_{k \in \Omega_{t, s,b}^{W}}p_{t, s,k,b}^{w,\max}), \quad \forall  t\in\mathcal{T},\label{eq:upperj}\end{align}
    \vspace{-15pt}
    \begin{align}
       \sum_{t\in\mathcal{T}}CC_{t} = \rho^{VC} \sum_{t\in\mathcal{T}}VC_{t}   ,\label{eq:upperk}\end{align}
Revenue from volumetric charges $VC_{t}$ is modelled for reconductored lines in~\eqref{eq:upperg} and parallel lines in~\eqref{eq:upperh}. Capacity-based charge revenue $CC_{t}$, proportional to maximum installed capacity, is modelled in~\eqref{eq:upperj}. Eq.~\eqref{eq:upperk} sets a ratio between total volumetric and capacity charge revenue. \vspace{-0.7em}
\subsubsection{Revenue adequacy constraints}
     \begin{align}
        \sum_{t \in \mathcal{T}} \frac{1}{(1+r)^{t-1}}(MS_t+VC_{t}+CC_{t} -C_t)\geq 0 ,\label{eq:upperl} \end{align}
Revenue adequacy for the network planner is ensured by constraint~\eqref{eq:upperl}, which requires that the total revenue from merchandising surplus $MS_t$ and both forms of network charges $VC_{t}+CC_{t}$ must be sufficient to cover the costs of line investments $C_t$. \vspace{-0.7em}
\subsubsection{Logical constraints for binary investment variables}
     \begin{align}   
         \sum_{t \in \mathcal{T}}\sum_{j \in \mathcal{J}}b_{t,l, j}^{R} \leq 1 , \quad   \forall l \in \mathcal{L}^{R},  \label{eq:upperm} 
\end{align}
\vspace{-15pt}
\begin{align} 
         z^R_{t,l}=\sum_{\hat{t} \leq t}\sum_{j \in \mathcal{J}\setminus\{0\}}b_{\hat{t},l, j}^{R}, \quad  \forall t \in\mathcal{T}, l \in \mathcal{L}^{R},
         \label{eq:uppern} 
\end{align}
\vspace{-15pt}
\begin{align} 
z^P_{t,l,m}=\sum_{c\in\mathcal{C}: I_{c,l}=m}o_{t,c},\quad \forall t\in\mathcal{T}, l\in\mathcal{L}^P, m \in \mathcal{M},
         \label{eq:upperq} 
\end{align}
\vspace{-15pt}
\begin{align}
\sum_{m \in \mathcal{M}} z^P_{t,l,m} \leq 1, \quad \forall t \in \mathcal{T}, l \in \mathcal{L}^P, \label{eq:upperq2}
\end{align}
\vspace{-15pt}
\begin{align}
\sum_{m \in \mathcal{M}} m ( z^P_{t,l,m} - z^P_{t-1,l,m})\geq 0,\hspace{2pt} \forall t \in \mathcal{T} \setminus \{1\}, l \in \mathcal{L}^P, \label{eq:upperr}
\end{align}
\vspace{-15pt}
\begin{align} 
z^R_{t,l} + \sum_{m \in \mathcal{M}} z^P_{t,l,m} \leq 1,\quad\forall t\in\mathcal{T}, l\in\mathcal{L}^P \cap \mathcal{L}^R,
         \label{eq:uppers} 
\end{align}
\vspace{-15pt}
\begin{align}   
    \tau^V_l\geq 0,\tau^C\geq 0, z^R_{t,l}\in\{0,1\},z^P_{t,l,m}\in\{0,1\}.\label{eq:uppert} 
\end{align}
The remaining constraints define logical relationships for the binary line expansion variables and establish feasible domains for the decision variables. Eqs.~\eqref{eq:upperm}--\eqref{eq:uppern} model reconductoring as a lumpy and irreversible capacity expansion. Eq.~\eqref{eq:upperm} ensures that each line $l \in \mathcal{L}^R$ can be reconductored at most once across all periods and capacity indices. Equation~\eqref{eq:uppern} defines the cumulative activation status $z^R_{t,l}$, which reflects whether line $l$ has been reconductored by period $t$ and remains active thereafter. Parallel line expansion is modelled through Eqs.~\eqref{eq:upperq}--\eqref{eq:upperr}. In this formulation, each configuration $c \in \mathcal{C}$ represents a specific parallel expansion plan, where $I_{c,l} = m$ indicates that configuration $c$ installs $m$ parallel circuits on corridor $l$. Eq.~\eqref{eq:upperq} defines the binary variable $z^P_{t,l,m}$, which indicates whether line $l$ has $m$ parallel circuits operating in period $t$. The value of $z^P_{t,l,m}$ is determined by summing over the binary variables $o_{t,c}$, which track the activation of configuration $c$ in period $t$. Eq.~\eqref{eq:upperq2} ensures that for each candidate corridor $l$ and period $t$, only one configuration $m$ is selected. Eq.~\eqref{eq:upperr} enforces a non-decreasing investment policy for parallel lines. Constraint~\eqref{eq:uppers} ensures that for any line $l$ that can either be reconductored or have parallel lines added (i.e., $l \in \mathcal{L}^P \cap \mathcal{L}^R$), only one of these investment strategies can be implemented. Finally, equation~\eqref{eq:uppert} restricts the application of reconductoring and parallel line expansion decisions to their respective eligible sets of lines ($\mathcal{L}^R$ and $\mathcal{L}^P$) and enforces the non-negativity of network charges.
\end{subequations}
\end{fleqn}

\subsection{Lower-level Problem}\label{section2d}
The lower-level joint chance-constrained market clearing problem is specified in Problem~\eqref{eq:lower}. This chance-constrained market-clearing model is similar to that in \cite{8936474}, but excludes unit commitment constraints and incorporates line expansion decisions. All dual variables are reported in each constraint after a colon. The upper-level variables in \(\Xi_{UL}\) serve as parameters in the lower-level problem. Under the SLA reformulation, the lower-level primal and dual variables are given by
\[
\Xi_{LL}^{P,SLA} = \{d_{t,s,k,b} ,g_{t,s,k,b},p^{w,sch}_{t,s,k,b},p^{w,cur}_{t,s,k,b},u_{t,s},v_{t,s,i}\}
\]
\[
\begin{aligned}
\Xi_{LL}^{D,SLA} = \{ & \pi_{t,s},\, \varphi^{D,\min}_{t,s,k,b},\, \varphi^{D,\max}_{t,s,k,b},\, \varphi^{G,\min}_{t,s,k,b},\, \varphi^{G,\max}_{t,s,k,b},\, \varphi^{W,sch}_{t,s,k,b}, \\
& \varphi^{W,cur,\min}_{t,s,k,b},\, \varphi^{W,cur,\max}_{t,s,k,b},\, \mu^{u,\min}_{t,s},\, \mu^{v,\min}_{t,s,i},\, \mu^{1}_{t,s}, \\
& \mu_{t,s,l,i}^{2},\, \mu_{t,s,l,i}^{3},\, \mu_{t,s,l}^{4},\, \mu_{t,s,l}^{5} \}.
\end{aligned}
\]
The lower-level problem is defined as follows: 
\begin{fleqn}
\begin{subequations}\label{eq:lower}
\vspace{-0.7em}
\subsubsection{Objective function}
\begin{align}
&\max\sum_{t \in \mathcal{T}} \sum_{s \in \mathcal{S}}  \sum_{b \in \mathcal{B}}\bigg(\sum_{k \in \Omega_{t, s,b}^{D}} \hat{c}_{t,s, k,b}^{d} d_{t,s, k,b}\notag\\[-6pt]
& -\sum_{k \in \Omega_{t, s,b}^{G}} \hat{c}_{t,s,k,b}^{g} g_{t,s,k,b}-\sum_{k \in \Omega_{t, s,b}^{W}} \hat{c}_{t,s,k,b}^{w} p^{w,sch}_{t,s,k,b} \label{eq:lower_obj}\\[-6pt]
& -\sum_{k \in \Omega_{t, s,b}^{W}} c^{cur}_{t,s,k,b}p^{w,cur}_{t,s,k,b}\bigg)\notag
\end{align}
The objective function of the lower-level problem~\eqref{eq:lower_obj} maximizes the social welfare, considering the actual marginal benefits and costs of market participants, which are influenced by the volumetric network charges $\tau^V_{l}$ as modelled in Eqs.~\eqref{eq:upperb}--\eqref{eq:upperd}. Additionally, the wind curtailment costs are incorporated into the objective function.
\vspace{-0.7em}
\subsubsection{Power balance, generator and consumer constraints}
\begin{align}
&\sum_{b\in\mathcal{B}}(\sum_{k \in \Omega_{t, s,b}^{D}}d_{t, s,k,b}-\sum_{k \in \Omega_{t, s,b}^{G}}g_{t,s,k,b}-\sum_{k \in \Omega_{t, s,b}^{W}}p^{w,sch}_{t,s,k,b})\notag\\
& =0:\pi_{t,s} \quad \forall t\in \mathcal{T}, s\in \mathcal{S}, \hspace{80pt}   \label{eq:lowerb} \end{align}
\vspace{-22pt} 
\begin{align}   
&g_{t, s,k,b}^{ \min }  \leq g_{t, s,k,b} \leq g_{t, s,k,b}^{ \max }:\varphi^{G,\min}_{t,s,k,b},\varphi^{G,\max}_{t,s,k,b}\notag\\ 
&\forall t \in \mathcal{T}, s\in \mathcal{S},  k \in \Omega_{t, s,b}^{G}, b \in \mathcal{B}, \label{eq:lowerc} \end{align}
\vspace{-22pt} 
\begin{align}   
&d_{t, s,k,b}^{ \min } \leq d_{t, s,k,b}\leq d_{t, s,k,b}^{ \max }:\varphi^{D,\min}_{t,s,k,b},\varphi^{D,\max}_{t,s,k,b}\notag\\ 
&\forall t \in \mathcal{T}, s\in \mathcal{S}, k \in \Omega_{t,s, b}^{D}, b \in \mathcal{B}, \label{eq:lowerd} \end{align}
The power balance constraints are modelled in constraint~\eqref{eq:lowerb}. The supply and demand upper and lower limits of generator and consumers are modelled in Eqs.~\eqref{eq:lowerc} and~\eqref{eq:lowerd}, respectively.\vspace{-0.7em}
\subsubsection{Wind power constraints}
\begin{align}   
&p^{w,sch}_{t,s,k,b} = p^{w,fore}_{t,s,k,b} - p^{w,cur}_{t,s,k,b}:\varphi^{W,sch}_{t,s,k,b}\notag\\ 
&\forall t \in \mathcal{T},  s\in \mathcal{S},  k \in \Omega_{t, s,b}^{W}, b \in \mathcal{B}, \label{eq:lowere} \end{align}
\vspace{-22pt} 
\begin{align}   
&0 \leq p^{w,cur}_{t,s,k,b} \leq p^{w,fore}_{t,s,k,b}:\varphi^{W,cur,\min}_{t,s,k,b},\varphi^{W,cur,\max}_{t,s,k,b}   \notag\\ 
&\forall t \in \mathcal{T},  s\in \mathcal{S},  k \in \Omega_{t, s,b}^{W}, b \in \mathcal{B}, \label{eq:lowerf} \end{align}
Following the approach in~\cite{8936474}, the uncertain nature of wind generation is represented as a forecast plus by a random forecasting error: $\tilde{p}^{w,sch}_{t,s,k,b} = p^{w,fore}_{t,s,k,b} + \tilde{e}_{t,s,k,b}$. Consequently, the scheduled wind generation $p^{w,sch}_{t,s,k,b}$, determined by constraint~\eqref{eq:lowere} subject to a specified wind curtailment $p^{w,cur}_{t,s,k,b}$ (constrained by~\eqref{eq:lowerf}), can deviate from the actual wind generation after curtailment. Note that unlike \cite{9479716,6881718}, we consider non-dispatchable wind power to acknowledge that wind energy is a non–dispatchable resource.\vspace{-0.7em}
\subsubsection{Flow limit constraints}
\begin{align}
&\sum_{c\in\mathcal{C}}\sum_{b\in\mathcal{B}} (\sum_{k\in\Omega_{t, s,b}^{G}}s_{b,l,c}o_{t,c} g_{t,s,k,b} +\sum_{k\in\Omega_{t, s,b}^{W}} s_{b,l,c}o_{t,c} (p^{w,sch}_{t,s,k,b} \notag\\[-5pt]
&+ \tilde{e}_{t,s,k,b}) - \sum_{k\in\Omega_{t, s,b}^{D}} s_{b,l,c}o_{t,c} d_{t,s,k,b}) \label{eq:ll:jcc1}\\[-5pt]
&\leq \mathcal{F}^{0}_{l}(1+\sum_{\hat{t} \leq t}\sum_{j\in\mathcal{J}}j b^R_{\hat{t},l,j})+\sum_{m \in \mathcal{M}} m z^P_{t,l,m} \overline{F}^{C}_{l},\notag
\end{align}
\vspace{-15pt} 
\begin{align}
&\sum_{c\in\mathcal{C}}\sum_{b\in\mathcal{B}} (\sum_{k\in\Omega_{t, s,b}^{G}}s_{b,l,c}o_{t,c} g_{t,s,k,b} +\sum_{k\in\Omega_{t, s,b}^{W}} s_{b,l,c}o_{t,c} (p^{w,sch}_{t,s,k,b}\notag\\[-5pt]
& + \tilde{e}_{t,s,k,b}) - \sum_{k\in\Omega_{t, s,b}^{D}} s_{b,l,c}o_{t,c} d_{t,s,k,b}) \label{eq:ll:jcc2}\\[-5pt]
&\geq -\mathcal{F}^{0}_{l}(1+\sum_{\hat{t} \leq t}\sum_{j\in\mathcal{J}}j b^R_{\hat{t},l,j})-\sum_{m \in \mathcal{M}} m z^P_{t,l,m} \overline{F}^{C}_{l}.\notag
\end{align}
\end{subequations}
\end{fleqn}
Due to the uncertain wind forecasting error $\tilde{e}_{t,s,k,b}$ and to ensure operational safety, the limits of transmission lines must be jointly guaranteed with high probability for each market-clearing period $s$ and investment period $t$ and it is thus formulated as an RHS-WDRJCC. For all \( t \in \mathcal{T} \) and \( s \in \mathcal{S} \), the following RHS-WDRJCC holds:
\begin{equation}
\begin{aligned}
& \sup_{\mathbb{P} \in \mathcal{F}^e_N(\theta)} \mathbb{P} \left\{ 
\begin{array}{l}
\eqref{eq:ll:jcc1}, \eqref{eq:ll:jcc2}, \forall l\in\mathcal{L}
\end{array} \right\} \geq 1 - \epsilon.
\end{aligned}
\label{eq:JCC:bilevel}
\end{equation}
In Eqs.~\eqref{eq:ll:jcc1}--\eqref{eq:ll:jcc2}, the right-hand side represents the transmission capacity of line~$l$ at time~$t$, based on upper-level investment decisions. The first term, $\mathcal{F}^{0}_{l}\left(1 + \sum_{\hat{t} \leq t} \sum_{j \in \mathcal{J}} j \, b^R_{\hat{t},l,j} \right)$, models the capacity of line $l$ after reconductoring. The sum $\sum_{\hat{t} \leq t} \sum_{j \in \mathcal{J}} j \, b^R_{\hat{t},l,j}\mathcal{F}^{0}_{l}$ reflects the cumulative upgrade applied by time $t$. The second term, $\sum_{m \in \mathcal{M}} m \, z^P_{t,l,m} \, \overline{F}^{C}_{l}$, accounts for the capacity added through parallel line investments. To model the impact of wind power forecast errors on transmission line flows, the Power Transfer Distribution Factor (PTDF) formulation is applied within a DC load flow model. The RHS-WDRJCC~\eqref{eq:JCC:bilevel} purpose is to satisfy the line thermal limits against uncertain power flows resulting from wind generation variability across all transmission lines. The presence of different network configurations $c$, arising from parallel line expansion decisions, results in varying system susceptance matrices and, consequently, distinct PTDF values $s_{b,l,c}$. To account for the evolving system topology during optimization due to different investment decisions, the PTDF $s_{b,l,c}$ is pre-computed for all configurations $c$. The upper-level decisions $o_{t,c}$ then select the appropriate PTDF for the realized network configuration within the optimization process.

\subsection{Computational Challenges}\label{section:discussion_appro}
It is challenging to solve the proposed bilevel problem ~\eqref{eq:upper},~\eqref{eq:lower} with lower-level RHS-WDRJCC~\eqref{eq:JCC:bilevel}. In fact, various solution techniques have been developed for RHS-WDRJCC problems, and Table~\ref{tab:comparison_methods} provides a comparative overview of these methods along with their applicability to different problem types. While \emph{exact} Mixed-Integer Programming (MIP) formulations for RHS-WDRJCC with valid inequalities and reduced big-M parameters have been developed~\cite{exact_milp_strengthened,2022_OR_exact_DRO}, convex approximations remain essential for large-scale problems or when WDRJCC constraints arise in lower-level market clearing models within bilevel optimization~\cite{zhou2024strengthenedfasterlinearapproximation} (as lower-level convexity is required to derive optimality conditions). Existing approximation methods include the Bonferroni approximation, which is simple but overly conservative~\cite{chen2023approximations,3eaeef4b}; the worst-case CVaR (W-CVaR) approximation, which yields less conservative solutions at the cost of higher computational complexity~\cite{3eaeef4b}; and the Linear Approximation (LA) approach~\cite{chen2023approximations}, which is considered the “best” convex inner-approximation for RHS-WDRJCC and is equivalent to W-CVaR under certain hyperparameters. However, both W-CVaR and LA introduce a large number of ancillary variables, resulting in substantial computational burden, particularly for large-scale problems \cite{zhou2024strengthenedfasterlinearapproximation}.

Our prior work on SFLA~\cite{zhou2024strengthenedfasterlinearapproximation} improves these existing approximations by reducing constraints and adding valid inequalities, significantly decreasing computational burden while maintaining solution quality. However, SFLA has been applied only to bilevel problems where uncertainty coefficients are parameters (Case 1). Unfortunately, SFLA’s efficiency deteriorates in bilevel problems with lower‐level RHS‐WDRJCC when the coefficients of random variables depend on upper‐level decisions (Case 2)—such as investment decisions $o_{t,c}$—as in our formulation~\eqref{eq:JCC:bilevel}. This dependency makes the reduced set $[N]_p$ in SFLA a function of upper-level variables, requiring disjunctive constraints formulated with big-M method, which are then linearized using big-M method again when deriving the Karush-Kuhn-Tucker (KKT) conditions, potentially leading to numerical instability. Additionally, Bonferroni approximation becomes inapplicable in Case 2 (see the formulation in \cite{zhou2024strengthenedfasterlinearapproximation}), as it requires pre-computations of uncertainty terms before solving bilevel problems, which explicitly depend on yet-to-be-determined upper-level decisions.

As shown in the RHS-WDRJCC line thermal constraints~\eqref{eq:JCC:bilevel}, our formulation corresponds to Case 2 where the coefficients of random variables $\tilde{e}_{t,s,k,b}$ incorporate upper-level investment variables $o_{t,c}$. Therefore, computationally efficient for RHS-WDRJCC are needed for our bilevel TEP problem with lower-level RHS-WDRJCC market clearing, without introducing additional numerical challenges and conservativeness.

\begin{table}[t]
\centering
\caption{Comparison of Solution Techniques}
\begin{tabular}{llcccc}
\toprule\midrule
\multirow{2}{*}{Method} 
 & \multirow{2}{*}{Type}
 & \multirow{2}{*}{Problem} 
 & \multirow{2}{*}{\shortstack{Single-\\level}} 
 & \multicolumn{2}{c}{Bilevel} \\ 
\cmidrule{5-6}
 &  &  &  & Case 1 & Case 2  \\
\midrule
Exact \cite{2022_OR_exact_DRO} & Exact & MIP & \checkmark  & \xmark  & \xmark   \\
 ExactS \cite{exact_milp_strengthened} & Exact & MIP  & \checkmark  & \xmark  & \xmark  \\
Bon. \cite{chen2023approximations} & Approx. & LP  & \checkmark  & \checkmark  & \xmark   \\
W-CVaR \cite{mohajerin2018data,chen2023approximations} & Approx. & LP  & \checkmark  & \checkmark  & \checkmark  \\
LA \cite{chen2023approximations} & Approx. & LP  & \checkmark  & \checkmark  & \checkmark   \\
SFLA \cite{zhou2024strengthenedfasterlinearapproximation} & Approx. & LP  & \checkmark  & \checkmark  & $*$ \\
The proposed SLA  & Approx. & LP  & \checkmark  & \checkmark  & \checkmark  \\
\midrule\bottomrule
\end{tabular}
\vspace{1mm}
\begin{minipage}{\linewidth}
\scriptsize 
$\checkmark$ means applicable, \xmark\ means not applicable, $*$ means that this method will bring additional numerical issues.\\
In the bilevel problem with an RHS-WDRJCC lower level, the coefficients of random variables are parameters and independent of upper-level variables (Case 1) or dependent on upper-level variables (Case 2).
\end{minipage}
\label{tab:comparison_methods}
\end{table}

\section{Solution Methodology}\label{section3}
In this section, we present the solution methodology for the bilevel TEP problem with lower-level RHS-WDRJCC dispatch~\eqref{eq:lower}. First, we review uncertainty modelling with Wasserstein DRCCP. Next, we describe existing solution methods for RHS-WDRJCC and then introduce the Strengthened Linear Approximation (SLA) as an efficient and numerically stable approach. Finally, we provide a comprehensive reformulation of the bilevel problem with lower-level RHS-WDRJCC~\eqref{eq:JCC:bilevel}.

\subsection{Uncertainty Modelling}
In this paper, we consider the \textit{Wasserstein} ambiguity set $\mathcal{F}(\theta)$ which is defined as a Wasserstein distance ball of radius $\theta$ around the empirical distribution $\mathbb{P}_N$. The DRCCP over Wasserstein balls has the form
\begin{subequations}
\label{eq:wdrccp}
\begin{align}
    \min_{\bm{x} \in \mathcal{X}} \quad & c(\bm{x}) \label{obj:wdrjcc}  \\[-6pt]
    \mbox{s.t.} \quad & \sup_{\mathbb {P} \in \mathcal{F}_N (\theta)} \mathbb{P}[\bm{\xi} \notin \mathcal{S}(\bm{x})] \leq \epsilon , \label{constr:wdrjcc}
\end{align}
\end{subequations}
where $c(\cdot)$ is the objective function, $\mathcal{X} \subseteq \mathbb{R}^L$ is a compact polyhedron for the decision variables $\bm{x} \in \mathbb{R}^L$, $\mathcal{S}(\bm{x})\subseteq \mathbb{R}^K$ is a decision-dependent safety set, and $\epsilon \in (0,1)$ is a specified risk tolerance. We also assume the ball radius $\theta >0$. Problem~\eqref{eq:wdrccp} minimizes a cost function $c(\bm{x})$ and ~\eqref{constr:wdrjcc} ensures that the worst-case probability of system being unsafe (i.e., $\bm{\xi}$ falling outside a decision-dependent safety set $\mathcal{S}(\bm{x})$) with low probability $\epsilon$ under every distribution resides in the Wasserstein ball $\mathbb{P}\in \mathcal{F}_N (\theta)$. 

Following \cite{2022_OR_exact_DRO} and \cite{exact_milp_strengthened}, we consider the $1$\textit{-Wasserstein distance} based on a general norm $\|\cdot\|$ and two distributions $\mathbb{P}, \mathbb{P}'$:
\begin{align}\label{eq:def:w-distance}
        d_W (\mathbb{P}, \mathbb{P}') \coloneqq \inf\limits_{\Pi \in \mathcal{P}(\mathbb{P}, \mathbb{P}')} \mathbb{E}_{(\bm{\xi}, \bm{\xi'}) \sim \Pi} [\| \bm{\xi} - \bm{\xi'} \|],
\end{align}
where $\mathcal{P}(\mathbb{P}, \mathbb{P}')$ is a set of joint distributions with marginal distributions $\mathbb{P}$ and $\mathbb{P}'$. Then, the Wasserstein ambiguity set is
\begin{equation}
    \mathcal{F}_N(\theta) \coloneqq \left\{ \mathbb{P} \mid  d_W(\mathbb{P}_N, \mathbb{P}) \leq \theta \right\}.
\end{equation}
Given a decision $\bm{x} \in \mathcal{X}$ and a sample $\bm{\xi}_i$, we define the distance from $\bm{\xi}_i$ to the the unsafe set $\bar{\mathcal{S}}(\bm{x})$ (i.e., the closed complement of $\mathcal{S}(\bm{x})$) as:
\begin{equation}
    \mathrm{dist}\left(\bm{\xi}_i, \mathcal{S}(\bm{x})\right) = \inf_{\bm{\xi'} \in \mathbb{R}^K} \left\{ \| \bm{\xi}_i - \bm{\xi'} \| \mid  \bm{\xi'} \notin \mathcal{S}(\bm{x}) \right\}.\label{eq:dist:original}
\end{equation}

This paper focuses on RHS-WDRJCCPs\footnote{We also we say that the chance constraint has left-hand side (LHS) uncertainty if the safety set is defined as follows with $\bm{A}\neq 0$, $\mathcal{S}(\bm{x}) \coloneqq \left\{ \bm{\xi} :  (\bm{b}_p-\bm{A}^\top\bm{x} )^\top\bm{\xi} + d_p - \bm{a}_p^\top \bm{x} \geq 0, \  p\in [P] \right\}.$}, where the safety set is defined as:
\begin{align}
\label{eq:def:safety_set}
    \mathcal{S}(\bm{x}) \coloneqq \left\{ \bm{\xi} :  \bm{b}_p^\top \bm{\xi} + d_p - \bm{a}_p^\top \bm{x} \geq 0, \  p\in [P] \right\},
\end{align}
where $\text{for given } a_p \in \mathbb{R}^K, b_p \in \mathbb{R}^L \text{ and } d_p \in \mathbb{R} \text{ for all } p \in [P]$. The $[P]$ constraints, indexed by $[P] \coloneqq \{1, \cdots, P\}$, must be met jointly with high probability $1 - \epsilon$.

\subsection{Existing Exact and Approximation Schemes}\label{sec:dic_exact_app}
Suppose we have collected $N$ independent and identically distributed (i.i.d.) samples $\{ \bm{\xi}_i \}_{i\in [N]}$ for the random vector $\bm{\xi}$, where the index set is defined as $[N] \coloneqq \{1, \cdots, N\}$. Building on the relation given by~\cite{2022_OR_exact_DRO} (see Equation~12 in~\cite{2022_OR_exact_DRO}) under the safety set defined in~\eqref{eq:def:safety_set}, we can reformulate $ \mathrm{dist}\left(\bm{\xi}_i, \mathcal{S}(\bm{x})\right)$~\eqref{eq:dist:original} to the following analytical expression:
\begin{equation}\label{eq:dist_anal}
\begin{aligned}
    \mathrm{dist}\left(\bm{\xi}_i, \mathcal{S}(\bm{x})\right)
        &= \left(\min_{p\in [P]} \dfrac{\bm{b}_p^\top \bm{\xi}_i + d_p - \bm{a}_p^\top \bm{x}}{\| \bm{b}_p \|_*} \right)^+ \\[-6pt]
        &= \min_{p\in [P]} \left(\dfrac{\bm{b}_p^\top \bm{\xi}_i + d_p - \bm{a}_p^\top \bm{x}}{\| \bm{b}_p \|_*} \right)^+,
\end{aligned}
\end{equation}
where $\|\cdot\|_*$ is the dual norm and $(x)^+=\max\{0,x\}$ for $x\in\mathbb{R}$.

By introducing ancillary variables $s\in \mathbb{R}$ and $\bm{r} \in \mathbb{R}^N$, \cite{2022_OR_exact_DRO} showed that RHS-WDRJCC~\eqref{constr:wdrjcc} can be expressed as the following \emph{exact} formulation:
\begin{subequations}\label{eq:exact1}
    \begin{align}
        &  s\geq 0, \bm{r} \geq \bm{0}, \label{eq:exact1:c1} \\[-5pt]
        & \epsilon N s - \sum\limits_{i=1}^N r_i \geq \theta N, \label{eq:exact1:c2} \\[-5pt]
        & \left( \frac{\bm{b}_p^\top \bm{\xi}_i + d_p - \bm{a}_p^\top \bm{x}}{\| \bm{b}_p \|_*} \right)^+ \geq s-r_i,\forall i\in [N], p \in [P], \label{eq:exact1:c3}
    \end{align}
\end{subequations}
where the $x$-feasible region of \eqref{eq:exact1} is defined as \begin{align}\label{eq:feasi_set_exact1}
    \mathcal{X}_\text{Exact} \coloneqq
    \left\{
        \bm{x} \in \mathcal{X} \mid 
        \exists s, \bm{r} :~\eqref{eq:exact1:c1}\mbox{--}\eqref{eq:exact1:c3}
    \right\},
\end{align}
Note that the minimum operator in \eqref{eq:dist_anal} is eliminated by adding constraints for all $p \in [P]$. \cite{2022_OR_exact_DRO} shows that Problem~\eqref{eq:exact1} can be transformed into a mixed-integer conic program, which is provided in Problem~\eqref{eq:exact_MILP} in Appendix~\ref{section:exact_formulation}. We also provide an overview of existing approximation schemes to RHS-WDRJCC~\eqref{constr:wdrjcc} including LA, W-CVaR and SFLA in Appendix~\ref{section:existing_approx.}. To efficiently solve the bilevel problem with lower-level WDRJCC market clearing, we next propose a SLA scheme for RHS-WDRJCC~\eqref{constr:wdrjcc} to effectively tightens the feasible region of ancillary variables without introducing additional numerical issues.

\subsection{Strengthened Linear Approximation (SLA)}
First, we prove that the exact formulation~\eqref{eq:exact1} can be strengthened by adding valid inequalities. Let $k \coloneqq \lfloor \epsilon N \rfloor$ and consider an index-sorting function $\sigma(\cdot)$ that sorts the set $\{\bm{b}_p^\top \bm{\xi}_i\}_{i \in [N]}$ from smallest to largest as $\{\bm{b}_p^\top \bm{\xi}_{\sigma(i)}\}_{i \in [N]}$. We further denote by $q_p$ the $(k+1)$-th smallest value:
\begin{equation}
    q_p \coloneqq \bm{b}_p^\top \bm{\xi}_{\sigma(k+1)} \label{eq:def:qp}
\end{equation}

\begin{theorem}\label{prop:1}
The feasible set $\mathcal{X}_\text{Exact}$~\eqref{eq:feasi_set_exact1} can be equivalently defined by the following set of strengthened constraints:
\vspace{-0.7em}
\begin{subequations}\label{eq:exact_stren}
    \begin{align}
        & Eqs.~\eqref{eq:exact1:c1}-\eqref{eq:exact1:c3} \label{eq:exact_stren:c1}\\[-5pt]
        & \dfrac{q_p + d_p - \bm{a}_p^\top \bm{x}}{\| \bm{b}_p \|_*} \geq s, &&\hspace{-50pt} \forall p\in [P].\label{eq:exact_stren:c2}
    \end{align}
\end{subequations}
In other words, we have:
\begin{align}\label{eq:feasible:exact_stren}
    \mathcal{X}_{\text{Exact}} =
        \left\{
        \bm{x} \in \mathcal{X} \mid \exists s, \bm{r} :~\eqref{eq:exact_stren:c1}\mbox{--}\eqref{eq:exact_stren:c2}
        \right\}.
\end{align}
\end{theorem}
\begin{IEEEproof}
    Please refer to Appendix~\ref{section:proof_prop1}.
\end{IEEEproof}

\noindent We then replace original distance function $\text{dist}(\bm{\xi}_i, \mathcal{S}(\bm{x}))$ with an inner approximation $\widehat{\text{dist}}(\bm{\xi}_i, \mathcal{S}(\bm{x}))$:
\begin{align}\label{eq:dist_approx}
    \widehat{\text{dist}}(\bm{\xi}_i, \mathcal{S}(\bm{x})) \coloneqq \kappa_i \left(\min_{p \in [P]} \frac{\bm{b}_p^\top \bm{\xi}_i + d_p - \bm{a}_p^\top \bm{x}}{\| \bm{b}_p \|_*} \right)
\end{align} where $\kappa \in[0,1]^N$ is a vector of slope parameters.
Next, we replace the non-convex constraints~\eqref{eq:exact1:c3} in \eqref{eq:exact_stren:c1} with the approximated distance function $\widehat{\mathrm{dist}}(\bm{\xi}_i, \mathcal{S}(\bm{x}))$ introduced in Eq.~\eqref{eq:dist_approx}, while keeping the the strengthening constraints~\eqref{eq:exact_stren:c2}. The SLA of RHS-WDRJCCPs~\eqref{constr:wdrjcc} is:
\begin{subequations}\label{eq:SLA}
    \begin{align}
        & s \geq 0, \bm{r} \geq \bm{0}, \label{eq:SLA:c1} \\[-3pt]
        & \epsilon N s - \sum\limits_{i\in [N]}r_i \geq \theta N, \label{eq:SLA:c2} \\[-3pt]
        & \kappa_i \frac{ \bm{b}_p^\top \bm{\xi}_i + d_p - \bm{a}_p^\top \bm{x}}{\| \bm{b}_p \|_*} \geq s - r_i ,\quad i \in [N],\ p\in [P], \label{eq:SLA:c3}\\[-3pt]
        & \frac{q_p + d_p - \bm{a}_p^\top \bm{x}}{\| \bm{b}_p \|_*} \geq s ,\quad p\in [P], \label{eq:SLA:c4}
    \end{align}
\end{subequations}
where the $x$-feasible region of \eqref{eq:SLA} is defined as \begin{align}\label{eq:feasi_set_SLA}
    \mathcal{X}_\text{SLA}(\bm{\kappa}) \coloneqq
    \left\{
        \bm{x} \in \mathcal{X} \mid 
        \exists s, \bm{r} :~\eqref{eq:SLA:c1}\mbox{--}\eqref{eq:SLA:c4}
    \right\},
\end{align}


For the theoretical properties of the proposed SLA \eqref{eq:SLA}, please refer to Appendix~\ref{section:SLA_property}. We prove that \(\mathcal{X}_\text{SLA}\) forms an inner approximation of \(\mathcal{X}_\text{Exact}\), shares the same feasible region as \(\mathcal{X}_\text{SFLA}\) under certain hyperparameters, and serves as a strengthened formulation of LA. Moreover, \(\mathcal{X}_\text{SLA}\) retains the original feasible region of \(\mathcal{X}_\text{SFLA}\) under specific hyperparameters and can become equivalent to the exact reformulation.

\subsection{Bilevel Reformulation under Different Approximations}
In this section, we will discuss the solution technique for applying the proposed SLA~\eqref{eq:SLA} to the RHS-WDRJCC line thermal constraints~\eqref{eq:JCC:bilevel} in the bilevel TEP problem. 

In the compact form of the safety set~\eqref{eq:def:safety_set}, the uncertainty term \(\bm{b}_p^\top \bm{\xi}_i\) corresponds to the expression \(\sum_{c \in \mathcal{C}} \sum_{b \in \mathcal{B}} \sum_{k \in \Omega_{t, s, b}^{W}} s_{b, l, c} \, o_{t, c} \, \tilde{e}_{t, s, k, b}\) in the WDRJCC flow limit constraints~\eqref{eq:JCC:bilevel}. Here, $\tilde{e}_{t,s,k,b}$ denotes the uncertain wind error, and the upper-level binary variable $o_{t,c}$ is a parameter in the lower-level problem. This entire expression is treated as our random vector $\xi$, supported on a closed convex set $\Xi \subseteq \mathbb{R}^m$, with $\bm{b}_p^\top=\bm{1}$. This treatment is valid because, for the lower-level optimization, $o_{t,c}$ is a fixed parameter determined by the upper-level planner, making the entire expression a function of the uncertain parameter $\tilde{e}_{t,s,k,b}$. In addition, through the definition of $q_p$ in \eqref{eq:def:qp}, we can pre-calculate $Q^{\max}_{t,s,l,c}= -\sum_{b\in\mathcal{B}}\sum_{k\in\Omega_{t, s,b}^{W}} s_{b,l,c}  \tilde{e}_{t,s,k,b}$ for \eqref{eq:ll:jcc1} and $Q^{\min}_{t,s,l,c}= \sum_{b\in\mathcal{B}}\sum_{k\in\Omega_{t, s,b}^{W}} s_{b,l,c}  \tilde{e}_{t,s,k,b}$ for \eqref{eq:ll:jcc2}, for all possible configurations $c\in\mathcal{C}$. In valid inequalities~\eqref{eq:SLA:c4} where $Q^{\max}_{t,s,l,c}$ and $Q^{\min}_{t,s,l,c}$ are involved, we can use the multiplications $\sum_{c\in\mathcal{C}}o_{t,c}Q^{\max}_{t,s,l,c}$ and $\sum_{c\in\mathcal{C}}o_{t,c}Q^{\min}_{t,s,l,c}$ to incorporate the impact of the upper-level investment decisions within the lower-level optimization problem. Appendix~\ref{section:SLA_reformulation1} presents the corresponding SLA approximation~\eqref{eq:SLA} of the RHS-WDRJCC line thermal constraints~\eqref{eq:JCC:bilevel}.

Now, the lower-level RHS-WDRJCC market clearing is approximated as a convex problem using SLA. To solve the bilevel optimization problem, we reformulate it into a single-level problem by replacing the lower-level problem with its Karush–Kuhn–Tucker (KKT) optimality conditions. This single-level reformulation can then be solved using the Gurobi solver for different convex approximations considered. The final single-level problem under the proposed SLA for the RHS-WDRJCC~\eqref{eq:JCC:bilevel} is detailed in Appendix~\ref{section:SLA_reformulation}. For benchmarking purposes, we apply the same reformulation procedure to two existing convex approximations: the LA and W-CVaR approaches \cite{chen2023approximations, mohajerin2018data} (see compact formulation in Appendix~\ref{section:existing_approx.}). Specifically, the RHS-WDRJCC formulation in \eqref{eq:JCC:bilevel} is approximated using these methods, and we reformulate it as a single-level problem using respective KKT optimality conditions. Their corresponding reformulations of the RHS-WDRJCC problem and resulting single-level models are provided in Appendix~\ref{section:benchmark}. 


\section{Case Study}\label{section4}
This case study investigates the proposed bilevel investment planning framework, incorporating three RHS-WDRJCC approximation schemes (SLA, LA and W-CVaR), using a modified Garver's 6-node transmission network~\cite{4435945,XIA2024110617}. Detailed case study settings, such as network topology illustration, procedures to generate maximum generation and demand capacities and bids, wind uncertainty modelling, are provided in Appendix~\ref{section:case_study setting}. The numerical case study is implemented using Python $3.9.13$ and solved with Gurobi $12.0.0$ on a system equipped with a 13th Gen Intel(R) Core(TM) i7-13700K CPU $@$ 3.40 GHz, 16 Core(s) with 32 Gb RAM. The solver is configured with the following tolerances and limits: \texttt{FeasibilityTol=1e-9}, \texttt{OptimalityTol=1e-9}, \texttt{IntFeasTol=1e-9}, \texttt{MIPGap = 0.0001\%}, and a wall-clock time limit of 4 hrs (\texttt{TimeLimit=14,400s}). 

\subsection{Investment and Reliability Analysis}

This section evaluates the proposed model in terms of social welfare, investment costs, and out-of-sample reliability under varying levels of dispatch conservativeness. We consider three Wasserstein radii, $\theta \in \{0.1, 0.2, 0.3\}$, and six risk levels, $\epsilon \in \{0.3, 0.2, 0.1, 0.05, 0.025, 0.01\}$, corresponding to confidence levels/desired reliability of $\{70\%, 80\%, 90\%, 95\%, 97.5\%, 99\%\}$. The out-of-sample reliability is defined as the proportion of the $4,000$ out-of-sample wind samples where line thermal limits~\eqref{eq:ll:jcc1} and~\eqref{eq:ll:jcc2} are jointly satisfied. This reliability is tested for every $t\in\mathcal{T}$ and $s\in\mathcal{S}$. An investment strategy (optimal solutions under desired reliability $1-\epsilon$) is considered reliable if its corresponding out-of-sample reliability exceeds this desired level. To ensure comparability across experiments, a consistent random seed is employed throughout the evaluation process.

\begin{figure}[t] 
  \centering
  \includegraphics[width=\linewidth]{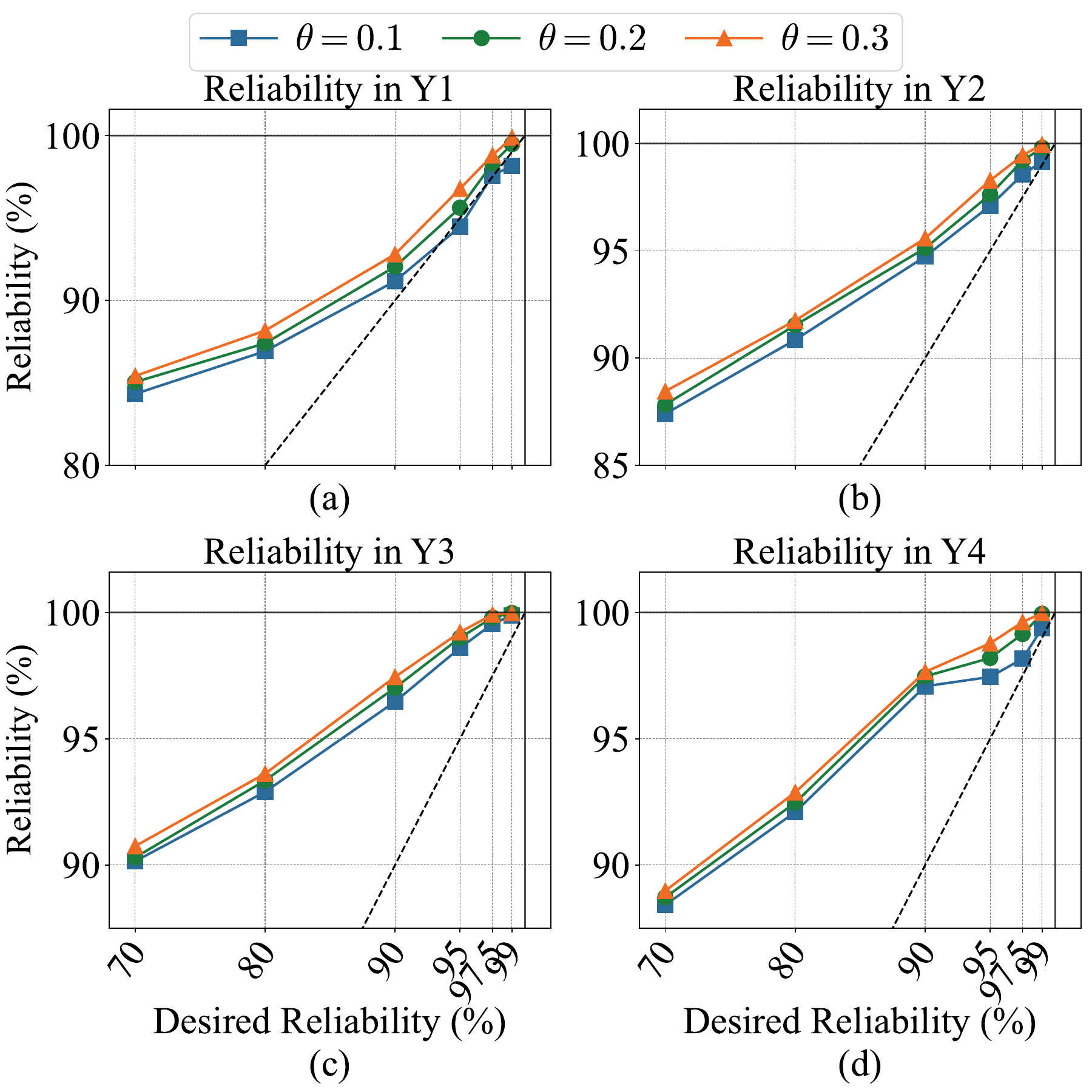}
  \caption{Impact of risk level $\epsilon$ and Wasserstein radius $\theta$ on out-of-sample reliability analysis. The x-axis represents the desired reliability $1-\epsilon$ $\{70\%,80\%,90\%,95\%,97.5\%,99\%\}$, while the y‐axis shows the out‐of‐sample probability that constraints are jointly satisfied. (a)--(d) shows the the out-of-sample reliability results for year 1 to year 4. The black dashed lines represent the desired reliability levels $1-\epsilon$ and the black solid lines represent the maximum $100\%$ reliability levels.}
  \label{fig:investment_relia}
\end{figure}

Fig.~\ref{fig:investment_relia} illustrates the impact of varying confidence levels $1-\epsilon$ and Wasserstein radii $\theta$ on key system performance metrics. The statistical data for out-of-sample reliability and investment decisions are presented in Table~\ref{tab:relia} and Table~\ref{tab:investment} of Appendix~\ref{section:statistical_results}, respectively. A key observation is that more conservative parameter settings in the dispatch problem, characterized by a lower risk level $\epsilon$ and a larger Wasserstein radius $\theta$, generally lead to enhanced out-of-sample reliability, as demonstrated in Table~\ref{tab:relia} and Fig.~\ref{fig:investment_relia}. Furthermore, the desired reliability targets are met across almost all years, with minor shortfalls observed with small radius $\theta = 0.1$ in the first year of out-of-sample evaluations. Specifically, with desired reliability levels of $95\%$ and $99\%$, the corresponding achieved reliabilities in the first year are $94.47\%$ and $98.17\%$, respectively, indicating a slight underachievement. However, out-of-sample reliability can be ensured with larger Wasserstein radii. Overall, the reliability significantly improves from $84.33\%$ under the least conservative dispatch scenario ($\epsilon=0.3$, $\theta=0.01$) to $99.98\%$ under the most conservative scenario ($\epsilon=0.01$, $\theta=0.3$).

As shown in Table~\ref{tab:investment} in Appendix~\ref{section:statistical_results}, increasing dispatch conservativeness tends to raise total transmission expansion investments. Specifically, investment costs increase by $35.69\%$, from \pounds127.50M to \pounds173.00M, when moving from the least to the most conservative dispatch scenarios. Notably, under greater conservativeness, represented by higher values of $1-\epsilon$ and $\theta$, the network planner expands investments by adding a second candidate line in corridor $(2,6)$ and increasing reinforcement of the existing lines $(2,3)$ and $(3,5)$. When the risk level $\epsilon$ is higher, the system allows certain violations of line flow constraints due to wind uncertainty, reducing the need for additional investment. This is because a lower confidence level in the chance constraints corresponds to a lower probability of strictly satisfying line thermal limits, thus reducing the impact of uncertain wind power on economic and operational performance.

In summary, increased dispatch conservativeness, characterized by larger Wasserstein radii and higher desired reliability levels, can lead to greater transmission investment needs. This arises from the requirement to ensure robust dispatch decisions and achieve the targeted out-of-sample reliability under uncertain wind power generation. By explicitly modelling both the probability of constraint violations and the size of the ambiguity set through the WDRJCC framework, system operators and network planners can make informed decisions regarding system reliability and grid investment in the context of increasing renewable energy penetration in TEP.

\subsection{Computing Performance Analysis}
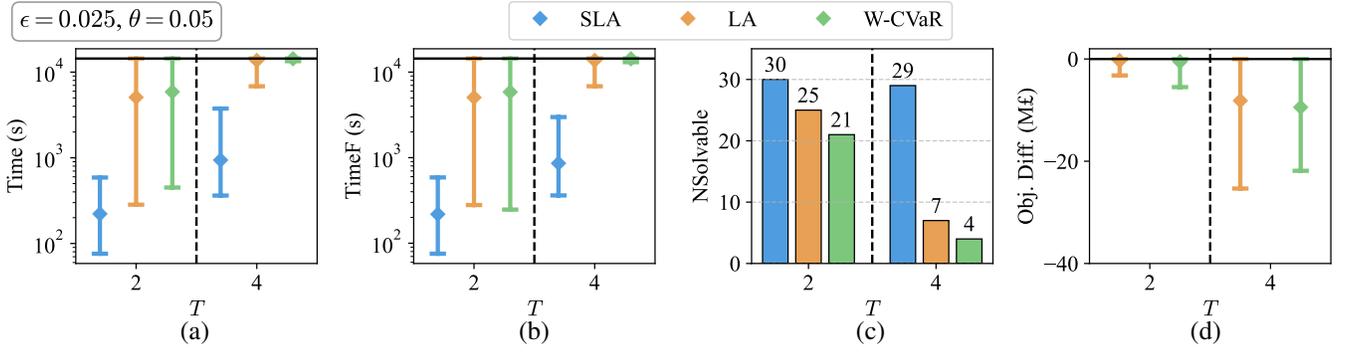
\begin{figure*}[t]
  \centering
  \begin{tikzpicture}
    \node[anchor=south west,inner sep=0] (img)
      {\includegraphics[width=\linewidth]{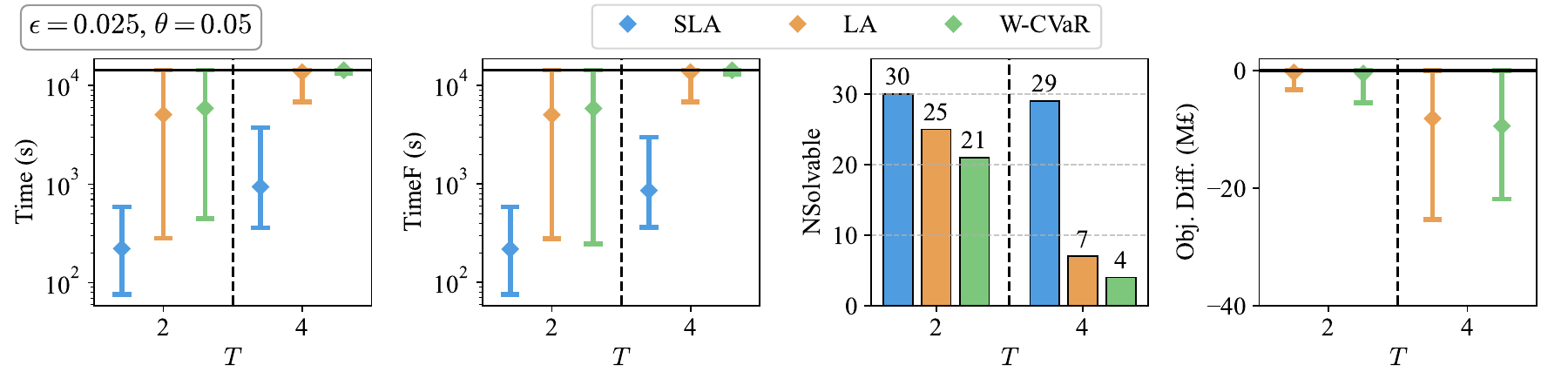}};
    \begin{scope}[x={(img.south east)}, y={(img.north west)}]
      \node[anchor=west] at (0.13,0.02) {(a)};
      \node[anchor=west] at (0.378,0.02) {(b)};
      \node[anchor=west] at (0.625,0.02) {(c)};
      \node[anchor=west] at (0.87,0.02) {(d)};
    \end{scope}
  \end{tikzpicture}
  \caption{Computational performance comparison for the proposed problem under $(\epsilon,\theta) = (0.025,0.05)$, $T\in\{2,4\}$. For \textit{Time }(s) and \textit{TimeF} (s) plots, dots represent the mean values of the $30$ random runs, with error bars representing the $95\%$ percentile interval (from 2.5th to 97.5th) of 30 randomly generated instances. In contrast, error bars in the \textit{Obj. Diff.} (M\pounds) plot represent $95\%$ percentile interval of runs for which the proposed SLA and benchmark are both solvable/feasible. The black horizontal line in the \textit{Time }(s) and \textit{TimeF} (s) plots represents the $14,400$s \texttt{TimeLimit}. The black horizontal line in the \textit{Obj. Diff.} (M\pounds) plots indicates zero difference, while negative values denote a lower objective value (social welfare) achieved by the benchmark method compared to the proposed SLA within the \texttt{TimeLimit}.}
  \label{fig:bilevel_com_results_rep}
\end{figure*}

This section presents a comparative analysis of the computational performance of the proposed SLA against existing approximation schemes, namely the LA ~\cite{chen2023approximations} and W-CVaR~\cite{mohajerin2018data,chen2023approximations}. To ensure the robustness of our findings, we conducted 30 independent random runs for each combination of the following parameter settings: risk level $\epsilon \in \{0.05, 0.025\}$, Wasserstein radius $\theta \in \{0.01, 0.05\}$ and investment period $\mathcal{T}\in\{2,4\}$. Each random run utilizes 4 threads (\texttt{Threads=4}) and a wall-clock time limit of 4 hrs (\texttt{TimeLimit=14,400s}). Other solver setting remains unchanged as discussed beforehand. In each independent run, we generate bid prices for participants from a normal distribution; sample generator and consumer capacities at each node uniformly; draw $N$ wind forecasting‐error samples from the 1,000 training samples; and set the solver seed which can affect tie‐breaking behaviour. To assess the computational performance and solution quality of three approximation methods, we use six evaluation metrics: \emph{TimeF} (s), which records the time to obtain the first comparable high-quality feasible solution; \emph{Time} (s), representing the total solution time; \emph{Nsolvable}, indicating the number of runs out of 30 where a feasible solution was found within the 4-hour time limit; \emph{Bilevel Obj.} (M\pounds), denoting the bilevel objective value that represents the overall social welfare; \emph{Obj. Diff.} (\%), measuring the difference in objective values of benchmark models relative to the proposed SLA, computed only when both models are solvable; and \emph{Reli.} (\%), representing the out-of-sample joint satisfaction rate of the chance constraints. Full definitions for each metric are provided in Appendix~\ref{appendix:evaluation_metrics}. 

Fig.~\ref{fig:bilevel_com_results_rep} presents the distribution of key performance metrics for SLA, LA, and W-CVaR at \((\epsilon,\theta)=(0.025,0.05)\) and Fig.~\ref{fig:bilevel_com_results} in Appendix~\ref{section:statistical_results} illustrates these results across \textit{all} parameter settings. Tables~\ref{tab:computation_time} and~\ref{tab:bilevel_reliability}, also in Appendix~\ref{section:statistical_results}, report the average computation times, out-of-sample reliability, and objective values across \textit{all} parameter settings. As can be observed from Table~\ref{tab:computation_time} and Fig.~\ref{fig:bilevel_com_results_rep}(a)--(b), the proposed SLA demonstrates significant speedup in both \textit{Time} and \textit{TimeF}. For all 30 runs, the proposed SLA achieves an averaged speedup of up to $22.91\times$ and $26.54\times$ in solving to optimality (\textit{Time}) compared to LA and W-CVaR, respectively, and the speedup is up to $23.12\times$ and $26.80\times$ in finding the first comparable high-quality solution (\textit{TimeF}). In runs where all three methods found feasible solutions within the \texttt{TimeLimit}, SLA delivered up to a $15\times$ speedup over LA and W-CVaR in both \textit{Time} and \textit{TimeF}. In addition, for the complex case with $\mathcal{T}=4$, we observe that all 30 random instances under the LA and W-CVaR approaches nearly reach the \texttt{TimeLimit}. In contrast, the proposed SLA demonstrates a significant computational advantage, achieving an average solution time of $941.07s$. Moreover, Fig.~\ref{fig:bilevel_com_results_rep}(c) shows that SLA solves nearly all 30 instances. In contrast, LA and W-CVaR solve only 21–25 cases when \(T=2\) and just 4–7 cases when \(T=4\). This indicates that SLA demonstrates greater effectiveness in finding feasible solutions within the \texttt{TimeLimit}. Furthermore, Fig.~\ref{fig:bilevel_com_results_rep}(d) shows that SLA achieves superior optimality relative to LA and W-CVaR across solvable runs. This improved performance occurs because the benchmark methods often have large optimality gaps, failing to reach optimality within the specified time limit, even when feasible solutions are identified.

\section{Conclusion}\label{section5}

Economic transmission expansion planning (TEP) leverages market signals to support cost-effective investment. However, increasing uncertainty in market operations, particularly from renewable generation variability, requires a method explicitly internalizing uncertainty into market clearing. Joint chance constraints (JCC) provide an effective approach to deal with uncertainties. In addition to renewable uncertainties, the sensitivity of market participants to network charges and the need to ensure cost recovery for the network planner further complicate TEP problems. To address these issues, we propose a bilevel TEP model incorporating a lower-level distributionally robust JCC (DRJCC) market clearing problem. 

Traditional bilevel solution methods depend on lower-level convexity to derive optimality conditions; however, embedding DRJCC at the lower level makes the problem non-convex. Moreover, incorporating joint chance constraints at the lower level with uncertainty terms dependent on upper-level decisions further increases the difficulty of solving the problem. To overcome this challenge, we develop a more efficient and no-more conservative approximation called Strengthened Linear Approximation (SLA) by introducing valid inequalities.

Numerical results demonstrate that increased dispatch conservativeness, characterized by lower risk levels and larger Wasserstein radii, leads to higher transmission investment costs. The achievement of the desired reliability is validated through 4,000 out-of-sample reliability tests. Extensive numerical experiments demonstrate the effectiveness and superiority of the proposed SLA for efficiently solving the proposed bilevel TEP problem. The proposed SLA significantly improves computational efficiency (achieving up to a 26$\times$ speedup) and solution optimality compared to existing linear approximation methods, including linear approximation and worst-case conditional value-at-risk. 

Future research may explore the integrated planning of wind generation and transmission investments, along with the design of associated tariff structures.

\section*{Acknowledgements}
Yuxin Xia's and Yihong Zhou's work are supported by the Engineering Studentship from the University of Edinburgh. Iacopo Savelli's work has been supported by the European Union’s Horizon Europe programme under the Marie Skłodowska-Curie grant agreement No. 101148367.

\bibliographystyle{IEEEtran}
\bibliography{references.bib}

\begin{thebibliography}{10}
\providecommand{\url}[1]{#1}
\csname url@samestyle\endcsname
\providecommand{\newblock}{\relax}
\providecommand{\bibinfo}[2]{#2}
\providecommand{\BIBentrySTDinterwordspacing}{\spaceskip=0pt\relax}
\providecommand{\BIBentryALTinterwordstretchfactor}{4}
\providecommand{\BIBentryALTinterwordspacing}{\spaceskip=\fontdimen2\font plus
\BIBentryALTinterwordstretchfactor\fontdimen3\font minus \fontdimen4\font\relax}
\providecommand{\BIBforeignlanguage}[2]{{%
\expandafter\ifx\csname l@#1\endcsname\relax
\typeout{** WARNING: IEEEtran.bst: No hyphenation pattern has been}%
\typeout{** loaded for the language `#1'. Using the pattern for}%
\typeout{** the default language instead.}%
\else
\language=\csname l@#1\endcsname
\fi
#2}}
\providecommand{\BIBdecl}{\relax}
\BIBdecl

\bibitem{6345475}
Y.~Gu, M.~Ni, and R.~Bo, ``Transmission expansion planning considering economic and reliability criteria,'' in \emph{2012 IEEE Power and Energy Society General Meeting}, 2012, pp. 1--8.

\bibitem{9242289}
B.~K. Poolla, A.~R. Hota, S.~Bolognani, D.~S. Callaway, and A.~Cherukuri, ``Wasserstein distributionally robust look-ahead economic dispatch,'' \emph{IEEE Transactions on Power Systems}, vol.~36, no.~3, pp. 2010--2022, 2021.

\bibitem{SAVELLI2020113979}
\BIBentryALTinterwordspacing
I.~Savelli, A.~{De Paola}, and F.~Li, ``Ex-ante dynamic network tariffs for transmission cost recovery,'' \emph{Applied Energy}, vol. 258, p. 113979, 2020. [Online]. Available: \url{https://www.sciencedirect.com/science/article/pii/S0306261919316666}
\BIBentrySTDinterwordspacing

\bibitem{2411207121}
\BIBentryALTinterwordspacing
E.~Chojkiewicz, U.~Paliwal, N.~Abhyankar, C.~Baker, R.~O’Connell, D.~Callaway, and A.~Phadke, ``Accelerating transmission capacity expansion by using advanced conductors in existing right-of-way,'' \emph{Proceedings of the National Academy of Sciences}, vol. 121, no.~40, p. e2411207121, 2024. [Online]. Available: \url{https://www.pnas.org/doi/abs/10.1073/pnas.2411207121}
\BIBentrySTDinterwordspacing

\bibitem{8275041}
S.~Dehghan, N.~Amjady, and A.~J. Conejo, ``A multistage robust transmission expansion planning model based on mixed binary linear decision rules—part i,'' \emph{IEEE Transactions on Power Systems}, vol.~33, no.~5, pp. 5341--5350, 2018.

\bibitem{10607938}
S.~Mahmoudi, B.~Alizadeh, and S.~Dehghan, ``Distributionally robust chance-constrained transmission expansion planning using a distributed solution method,'' \emph{IEEE Transactions on Network Science and Engineering}, vol.~11, no.~6, pp. 6431--6447, 2024.

\bibitem{zhou2024strengthenedfasterlinearapproximation}
\BIBentryALTinterwordspacing
Y.~Zhou, Y.~Xia, H.~Yang, and T.~Morstyn, ``Strengthened and faster linear approximation to joint chance constraints with wasserstein ambiguity,'' 2024. [Online]. Available: \url{https://arxiv.org/abs/2412.12992}
\BIBentrySTDinterwordspacing

\bibitem{Zymler2013}
\BIBentryALTinterwordspacing
S.~Zymler, D.~Kuhn, and B.~Rustem, ``Distributionally robust joint chance constraints with second-order moment information,'' \emph{Mathematical Programming}, vol. 137, no.~1, pp. 167--198, 2013. [Online]. Available: \url{https://doi.org/10.1007/s10107-011-0494-7}
\BIBentrySTDinterwordspacing

\bibitem{chen2023approximations}
Z.~Chen, D.~Kuhn, and W.~Wiesemann, ``On approximations of data-driven chance constrained programs over wasserstein balls,'' \emph{Operations Research Letters}, vol.~51, no.~3, pp. 226--233, 2023.

\bibitem{mohajerin2018data}
P.~Mohajerin~Esfahani and D.~Kuhn, ``Data-driven distributionally robust optimization using the wasserstein metric: performance guarantees and tractable reformulations,'' \emph{Mathematical Programming}, vol. 171, no.~1, pp. 115--166, 2018.

\bibitem{ZHAN2022107417}
\BIBentryALTinterwordspacing
S.~Zhan, P.~Hou, G.~Yang, and J.~Hu, ``Distributionally robust chance-constrained flexibility planning for integrated energy system,'' \emph{International Journal of Electrical Power \& Energy Systems}, vol. 135, p. 107417, 2022. [Online]. Available: \url{https://www.sciencedirect.com/science/article/pii/S0142061521006566}
\BIBentrySTDinterwordspacing

\bibitem{8933104}
F.~Pourahmadi, J.~Kazempour, C.~Ordoudis, P.~Pinson, and S.~H. Hosseini, ``Distributionally robust chance-constrained generation expansion planning,'' \emph{IEEE Transactions on Power Systems}, vol.~35, no.~4, pp. 2888--2903, 2020.

\bibitem{8447238}
D.~Alvarado, A.~Moreira, R.~Moreno, and G.~Strbac, ``Transmission network investment with distributed energy resources and distributionally robust security,'' \emph{IEEE Transactions on Power Systems}, vol.~34, no.~6, pp. 5157--5168, 2019.

\bibitem{9028125}
A.~Velloso, D.~Pozo, and A.~Street, ``Distributionally robust transmission expansion planning: A multi-scale uncertainty approach,'' \emph{IEEE Transactions on Power Systems}, vol.~35, no.~5, pp. 3353--3365, 2020.

\bibitem{8442871}
D.~Pozo, A.~Street, and A.~Velloso, ``An ambiguity-averse model for planning the transmission grid under uncertainty on renewable distributed generation,'' in \emph{2018 Power Systems Computation Conference (PSCC)}, 2018, pp. 1--7.

\bibitem{8294298}
C.~Duan, W.~Fang, L.~Jiang, L.~Yao, and J.~Liu, ``Distributionally robust chance-constrained approximate ac-opf with wasserstein metric,'' \emph{IEEE Transactions on Power Systems}, vol.~33, no.~5, pp. 4924--4936, 2018.

\bibitem{9026959}
A.~Zhou, M.~Yang, M.~Wang, and Y.~Zhang, ``A linear programming approximation of distributionally robust chance-constrained dispatch with wasserstein distance,'' \emph{IEEE Transactions on Power Systems}, vol.~35, no.~5, pp. 3366--3377, 2020.

\bibitem{9961917}
B.~Chen, T.~Liu, X.~Liu, C.~He, L.~Nan, L.~Wu, X.~Su, and J.~Zhang, ``A wasserstein distance-based distributionally robust chance-constrained clustered generation expansion planning considering flexible resource investments,'' \emph{IEEE Transactions on Power Systems}, vol.~38, no.~6, pp. 5635--5647, 2023.

\bibitem{Xie2021}
\BIBentryALTinterwordspacing
W.~Xie, ``On distributionally robust chance constrained programs with wasserstein distance,'' \emph{Mathematical Programming}, vol. 186, no.~1, pp. 115--155, 2021. [Online]. Available: \url{https://doi.org/10.1007/s10107-019-01445-5}
\BIBentrySTDinterwordspacing

\bibitem{XIA2025124721}
\BIBentryALTinterwordspacing
Y.~Xia, I.~Savelli, and T.~Morstyn, ``Integrating local market operations into transmission investment: A tri-level optimization approach,'' \emph{Applied Energy}, vol. 378, p. 124721, 2025. [Online]. Available: \url{https://www.sciencedirect.com/science/article/pii/S0306261924021044}
\BIBentrySTDinterwordspacing

\bibitem{8031353}
I.~Savelli, A.~Giannitrapani, S.~Paoletti, and A.~Vicino, ``An optimization model for the electricity market clearing problem with uniform purchase price and zonal selling prices,'' \emph{IEEE Transactions on Power Systems}, vol.~33, no.~3, pp. 2864--2873, 2018.

\bibitem{BECK2023401}
\BIBentryALTinterwordspacing
Y.~Beck, I.~Ljubić, and M.~Schmidt, ``A survey on bilevel optimization under uncertainty,'' \emph{European Journal of Operational Research}, vol. 311, no.~2, pp. 401--426, 2023. [Online]. Available: \url{https://www.sciencedirect.com/science/article/pii/S0377221723000073}
\BIBentrySTDinterwordspacing

\bibitem{10038580}
S.~Ramyar, M.~Tanaka, A.~L. Liu, and Y.~Chen, ``Endogenous risk management of prosumers by distributionally robust chance-constrained optimization,'' \emph{IEEE Transactions on Energy Markets, Policy and Regulation}, vol.~1, no.~1, pp. 48--59, 2023.

\bibitem{8937824}
Y.~Dvorkin, ``A chance-constrained stochastic electricity market,'' \emph{IEEE Transactions on Power Systems}, vol.~35, no.~4, pp. 2993--3003, 2020.

\bibitem{Heitsch2022}
\BIBentryALTinterwordspacing
H.~Heitsch, R.~Henrion, T.~Kleinert \emph{et~al.}, ``On convex lower-level black-box constraints in bilevel optimization with an application to gas market models with chance constraints,'' \emph{Journal of Global Optimization}, vol.~84, no.~3, pp. 651--685, 2022. [Online]. Available: \url{https://doi.org/10.1007/s10898-022-01161-z}
\BIBentrySTDinterwordspacing

\bibitem{2022_OR_exact_DRO}
Z.~Chen, D.~Kuhn, and W.~Wiesemann, ``Data-driven chance constrained programs over wasserstein balls,'' \emph{Operations Research}, 2022.

\bibitem{exact_milp_strengthened}
N.~Ho-Nguyen, F.~K{\i}l{\i}n{\c{c}}-Karzan, S.~K{\"u}{\c{c}}{\"u}kyavuz, and D.~Lee, ``Distributionally robust chance-constrained programs with right-hand side uncertainty under wasserstein ambiguity,'' \emph{Mathematical Programming}, 2021.

\bibitem{SAVELLI2025101937}
\BIBentryALTinterwordspacing
I.~Savelli, D.~Howey, and T.~Morstyn, ``Locating large-scale energy storage: spillover effects, carbon emissions, and balancing costs across italy,'' \emph{Utilities Policy}, vol.~95, p. 101937, 2025. [Online]. Available: \url{https://www.sciencedirect.com/science/article/pii/S0957178725000529}
\BIBentrySTDinterwordspacing

\bibitem{8936474}
Y.~Yang, W.~Wu, B.~Wang, and M.~Li, ``Analytical reformulation for stochastic unit commitment considering wind power uncertainty with gaussian mixture model,'' \emph{IEEE Transactions on Power Systems}, vol.~35, no.~4, pp. 2769--2782, 2020.

\bibitem{9479716}
C.~Ning and F.~You, ``Deep learning based distributionally robust joint chance constrained economic dispatch under wind power uncertainty,'' \emph{IEEE Transactions on Power Systems}, vol.~37, no.~1, pp. 191--203, 2022.

\bibitem{6881718}
F.~Qiu and J.~Wang, ``Chance-constrained transmission switching with guaranteed wind power utilization,'' \emph{IEEE Transactions on Power Systems}, vol.~30, no.~3, pp. 1270--1278, 2015.

\bibitem{3eaeef4b}
\BIBentryALTinterwordspacing
W.~Chen, M.~Sim, J.~Sun, and C.-P. Teo, ``From cvar to uncertainty set: Implications in joint chance-constrained optimization,'' \emph{Operations Research}, vol.~58, no.~2, pp. 470--485, 2010. [Online]. Available: \url{http://www.jstor.org/stable/40605930}
\BIBentrySTDinterwordspacing

\bibitem{4435945}
S.~de~la Torre, A.~J. Conejo, and J.~Contreras, ``Transmission expansion planning in electricity markets,'' \emph{IEEE Transactions on Power Systems}, vol.~23, no.~1, pp. 238--248, 2008.

\bibitem{XIA2024110617}
\BIBentryALTinterwordspacing
Y.~Xia, I.~Savelli, and T.~Morstyn, ``An incentive regulation approach for balancing stakeholder interests in transmission merchant investment,'' \emph{Electric Power Systems Research}, vol. 234, p. 110617, 2024. [Online]. Available: \url{https://www.sciencedirect.com/science/article/pii/S0378779624005030}
\BIBentrySTDinterwordspacing

\bibitem{8036231}
E.~Nasrolahpour, J.~Kazempour, H.~Zareipour, and W.~D. Rosehart, ``A bilevel model for participation of a storage system in energy and reserve markets,'' \emph{IEEE Transactions on Sustainable Energy}, vol.~9, no.~2, pp. 582--598, 2018.

\bibitem{gefc2012}
T.~Hong, Glider, and P.~Pinson, ``{Global Energy Forecasting Competition 2012 - Wind Forecasting},'' \url{https://kaggle.com/competitions/GEF2012-wind-forecasting}, 2012, accessed: 2024-10-11.

\end{thebibliography}

\appendices

\nomenclature[S]{$\mathcal{L}$}{Set of transmission lines.}
\nomenclature[S]{$\mathcal{L}^R$}{Set of transmission lines for reconductoring.}
\nomenclature[S]{$\mathcal{L}^P$}{Set of transmission lines for parallel expansion.}
\nomenclature[S]{$\mathcal{J}$}{Set of lumpy capacity indices for reconductoring.}
\nomenclature[S]{$\mathcal{B}$}{Set of transmission network nodes.}
\nomenclature[S]{$\mathcal{T}$}{Set of investment periods.}
\nomenclature[S]{$\mathcal{S}$}{Set of operation periods.}
\nomenclature[S]{$\mathcal{C}$}{Set of configurations.}
\nomenclature[S]{$\mathcal{M}$}{Set indexing the number of candidate lines.}
\nomenclature[S]{$\Omega_{t,s,b}^{G}$}{Set of generators in investment period $t$, operation period $s$ in node $b$.}
\nomenclature[S]{$\Omega_{t,s,b}^{D}$}{Set of consumers in investment period $t$, operation period $s$ in node $b$.}
\nomenclature[S]{$\Omega_{t,s,b}^{W}$}{Set of wind farms in investment period $t$, operation period $s$ in node $b$.}
\nomenclature[S]{$\overline{\mathcal{F}}_{l,j}$}{Lumpy capacity expansion for reconductoring line $l$ with capacity index $j$.}

\nomenclature[P]{$r$}{The discount rate.}
\nomenclature[P]{$\Psi$}{Scaling factor representing the number of actual operational periods covered by each representative market-clearing period in an investment period.}
\nomenclature[P]{$c^d_{t,s,k,b}$}{Demand bid price (willingness-to-pay) for consumer $k$ in node $b$ at year $t$ and period $s$ without network costs.}
\nomenclature[P]{$c^{cur}_{t,s,k,b}$}{Wind curtailment cost for wind farm $k$ in node $b$ at year $t$ and period $s$.}
\nomenclature[P]{$c^p_{t,s,k,b}$}{Supply bid price (marginal cost) for generator $k$ in node $b$ at year $t$ and period $s$ without network costs.}
\nomenclature[P]{$K_l^{fix}$}{Fixed cost of expanding line $l$.}
\nomenclature[P]{$K_l^{var}$}{Variable cost of expanding line $l$.}
\nomenclature[P]{$\mathcal{F}^{0}_{l}$}{Existing capacity on transmission line $l$, $l \in \mathcal{L}$.}
\nomenclature[P]{$\overline{F}^{C}_{l}$}{Candidate line capacity on transmission line $l$, $l \in \mathcal{L}^R$.}
\nomenclature[P]{$\delta_{l,b}$}{Allocation factor for cost distribution among market participants.}
\nomenclature[P]{$d_{t,s,k,b}^{\min}$}{Minimum quantity of power demanded by consumer $k$ in investment period $t$, operation period $s$ in node $b$.}
\nomenclature[P]{$d_{t,s,k,b}^{\max}$}{Maximum quantity of power demanded by consumer $k$ in investment period $t$, operation period $s$ in node $b$.}
\nomenclature[P]{$g_{t,s,k,b}^{\min}$}{Minimum quantity of power produced by generator $k$ in investment period $t$, operation period $s$ in node $b$.}
\nomenclature[P]{$g_{t,s,k,b}^{\max}$}{Maximum quantity of power produced by generator $k$ in investment period $t$, operation period $s$ in node $b$.}
\nomenclature[P]{$p_{t,s,k,b}^{w,\max}$}{Maximum quantity of power produced by wind farm $k$ in investment period $t$, operation period $s$ in node $b$.}
\nomenclature[P]{$s_{b,l,c}$}{Power transfer distribution factor that measures the increase of flow on line $l$ given a power injection at bus $b$ under configuration $c$.}
\nomenclature[P]{$M$}{The big-M value.}
\nomenclature[P]{$\epsilon$}{The risk level of the RHS-WDRJCC.} 
\nomenclature[P]{$N$}{The number of samples (historical data) for random variables.}
\nomenclature[P]{$\kappa_i$}{The tunable hyperparameter in the proposed SLA and LA methods.}
\nomenclature[P]{$\theta$}{The Wasserstein radius.}
\nomenclature[P]{$Q^{\max}_{t,s,l,c}$}{The corresponding $q_p$ of the proposed SLA in the bilevel problem.}
\nomenclature[P]{$Q^{\min}_{t,s,l,c}$}{The corresponding $q_p$ of the proposed SLA in the bilevel problem.}
\nomenclature[P]{$\rho^{VC}$}{The capacity-to-volumetric revenue ratio.}
\nomenclature[P]{$w^{\max}_{t,s,l}$}{The tunable hyperparameters in the W-CVaR method.}
\nomenclature[P]{$w^{\min}_{t,s,l}$}{The tunable hyperparameters in the W-CVaR method.}
\nomenclature[P]{$I_{c,l}$}{Matrix element indicating the number of parallel circuits installed on corridor $l \in \mathcal{L}$ under configuration $c \in \mathcal{C}$.}

\nomenclature[R]{$\tilde{e}_{t,s,k,b}$}{Wind forecast error of wind farm $k$ at the investment period $t$, operation period $s$ in node $b$.}
\nomenclature[R]{$\tilde{p}^{w,sch}_{t,s,k,b}$}{Actual power output of wind farm $k$ at the investment period $t$, operation period $s$ in node $b$.}

\nomenclature[V]{$\hat{c}^d_{t,s,k,b}$}{Demand bid price (willingness-to-pay) for consumer $k$ in node $b$ at year $t$ and period $s$ with network costs.}
\nomenclature[V]{$\hat{c}^p_{t,s,k,b}$}{Supply bid price (marginal cost) for generator $k$ in node $b$ at year $t$ and period $s$ with network costs.}
\nomenclature[V]{$\hat{c}^w_{t,s,k,b}$}{Supply bid price (marginal cost) for wind farm $k$ in node $b$ at year $t$ and period $s$ with network costs.}
\nomenclature[V]{$z^R_{t,l}$}{Binary variable equal to one if line $l$ has been reconductored in investment period $t$, and zero otherwise, $l \in \mathcal{L}^R$.}
\nomenclature[V]{$b_{t, l, j}^{R}$}{Binary variable equal to one if the lumpy investment for reconductoring in additional capacity $j\mathcal{F}^0_{l}$ for line $l$ is made in investment period $t$, and zero otherwise, $l \in \mathcal{L}^R$.}
\nomenclature[V]{$z^P_{t,l,m}$}{Binary variable equal to one if corridor $l$ has been expanded by adding $m$ parallel circuits in investment period $t$, and zero otherwise, $l \in \mathcal{L}^P$.}
\nomenclature[V]{$o_{t,c}$}{Binary variable equal to one if configuration $c$ is activated in investment period $t$, and zero otherwise.}
\nomenclature[V]{$\tau^C, \tau^V_l$}{Transmission capacity-based network charges and volumetric-based network charges applied on line $l$.}
\nomenclature[V]{$g_{t,s,k,b}$}{Scheduled power of generator $k$ in investment period $t$, operation period $s$ in node $b$.}
\nomenclature[V]{$d_{t,s,k,b}$}{Scheduled power of consumer $k$ in investment period $t$, operation period $s$ in node $b$.}
\nomenclature[V]{$p^{w,sch}_{t,s,k,b}$}{Scheduled wind power of wind farm $k$ in investment period $t$, operation period $s$ in node $b$.}
\nomenclature[V]{$p^{w,cur}_{t,s,k,b}$}{Scheduled wind curtailment of wind farm $k$ in investment period $t$, operation period $s$ in node $b$.}
\nomenclature[V]{$\pi_{t,s,b}$}{Day-ahead market-clearing price (as a dual variable) at transmission node $b$ in investment period $t$, operation period $s$.}
\nomenclature[V, 8]{$\alpha_i, \beta, \tau$}{Auxiliary variables for W-CVaR.}
\nomenclature[V]{$\eta$}{Auxiliary variables to replace the nonlinear terms.}
\nomenclature[V]{$u,v_{i}$}{Auxiliary variables for the proposed SFLA and the existing LA method.}
\nomenclature[V]{$\psi,\mu$}{Dual variables corresponding to inequality day-ahead lower-level constraints.}

\printnomenclature

\section{Model Assumptions}\label{section:assumption}
Throughout this paper, the following assumptions are made:
\begin{enumerate}
    \item Reinforced/Reconductored lines are considered to be a subset of the existing transmission lines, denoted as $\mathcal{L}^R \subseteq \mathcal{L}^E$. For reconductoring projects, with the exception of resistance, conductor properties such as reactance and susceptance remain consistent across different conductor types with the same diameters~\cite{2411207121}.
    \item Transmission lines that are currently non-existent or unconnected are not eligible for re-conductoring but can be considered for candidate line expansion.
    \item For existing corridors $l$, both candidate line expansion and reconductoring represent potential investment strategies (i.e., $\mathcal{L}^P \cap \mathcal{L}^R$ may be non-empty). However, it is assumed that only one of these strategies can be implemented for any given existing corridor.
\end{enumerate}

\section{Exact Formulation} \label{section:exact_formulation}

\begin{lemma}[Lemma 1 in \cite{exact_milp_strengthened}]
\label{lemma:1}
    For any fixed $\bm{x} \in \mathcal{X}_\text{Exact}$, there exists ($\bm{r}$, $s$) such that $s$ is equal to the $(k+1)$-th smallest value amongst $\{ \mathrm{dist}(\bm{\xi}_i, \mathcal{S}(\bm{x})) \}_{i\in [N]}$ and constraints~\eqref{eq:exact1:c1}--\eqref{eq:exact1:c3} are satisfied.
\end{lemma}

\cite{2022_OR_exact_DRO} shows that Problem~\eqref{eq:exact1} can be transformed into the following mixed-integer conic program where $M$ is a suitably large (but finite) positive constant:
\begin{fleqn}
    \begin{subequations}\label{eq:exact_MILP}
    \begin{align}
        &\min_{\bm{x} \in \mathcal{X},z,s,r} \quad c(\bm{x}) \label{obj:exact_MILP}  \\
        & \mbox{s.t.} \\
        &\bm{z} \in \{0,1\}^N, s \geq 0, \bm{r} \geq \bm{0}, \label{con:exact_MILP:c1}\\
        & \epsilon N s - \sum\limits_{i=1}^N r_i \geq \theta N, \label{con:exact_MILP:c4}\\
        & \dfrac{\bm{b}_p^\top \bm{\xi}_i + d_p - \bm{a}_p^\top \bm{x}}{\| \bm{b}_p \|_*} + M z_i \geq s-r_i, \forall i\in [N], p\in[P],  \label{con:exact_MILP:c2}\\ 
        & M(1-z_i) \geq s-r_i, \quad \forall i\in [N]. \label{con:exact_MILP:c3}
    \end{align}
\end{subequations}
\end{fleqn}

\section{Existing Approximation Schemes} \label{section:existing_approx.}
In this section, we provide an overview of three existing inner approximation schemes for the RHS-WDRJCCP~\eqref{constr:wdrjcc}. We refer the readers to \cite{zhou2024strengthenedfasterlinearapproximation} for detailed explanation.
\subsection{Worst-case CVaR (W-CVaR) \cite{chen2023approximations}}
\cite{chen2023approximations} shows that the RHS-WDRJCCP~\eqref{constr:wdrjcc} can be approximated by the following worst-case CVaR constrained programs:
\begin{subequations}
\begin{align}
    &  \bm{\alpha} \geq \bm{0}, \beta \in \mathbb{R}, \tau \in \mathbb{R}, \label{eq:CVaR:c1}\\
    & \tau + \frac{1}{\epsilon} (\theta \beta + \frac{1}{N} \sum_{i\in [N]} \alpha_i) \leq 0, \label{eq:CVaR:c2}\\
    & \alpha_i \geq w_p (\bm{a}_p^\top \bm{x} - \bm{b}_p^\top\bm{\xi}_i - d_p) - \tau,\ i\in[N], p\in [P], \label{eq:CVaR:c3}\\
    & \beta \geq w_p \| \bm{b}_p \|_*,\ \forall p \in [P],\label{eq:CVaR:c4}
\end{align}\label{eq:CVaR}
\end{subequations}
where $\bm{w} \in \mathbb{R}^P$ is a tunable hyperparameter that affects the performance of the worst-case CVaR. The $x$-feasible region of \eqref{eq:CVaR} is defined as \begin{align}\label{eq:feasi_set_wcvar}
    \mathcal{X}_\text{W-CVaR} (\bm{w}) \coloneqq
    \left\{
        \bm{x} \in \mathcal{X} \mid 
        \exists s, \bm{r} :~\eqref{eq:CVaR:c1}\mbox{--}\eqref{eq:CVaR:c4}
    \right\}.
\end{align}
\subsection{Linear Approximation (LA) \cite{zhou2024strengthenedfasterlinearapproximation}}
With the approximated distance function $\widehat{\mathrm{dist}}(\bm{\xi}_i, \mathcal{S}(\bm{x}))$ introduced in Eq.~\eqref{eq:dist_approx}, the LA \cite{chen2023approximations} of RHS-WDRJCCPs~\eqref{constr:wdrjcc} is formulated as 
\begin{subequations}\label{eq:LA}
    \begin{align}
        & s \geq 0, \bm{r} \geq \bm{0}, \label{eq:LA:c1} \\
        & \epsilon N s - \sum\limits_{i\in [N]}r_i \geq \theta N, \label{eq:LA:c2} \\
        & \kappa_i \frac{ \bm{b}_p^\top \bm{\xi}_i + d_p - \bm{a}_p^\top \bm{x}}{\| \bm{b}_p \|_*} \geq s - r_i ,\quad i \in [N],\ p\in [P], \label{eq:LA:c3}
    \end{align}
\end{subequations}
where $\kappa_i \in [0, 1]$ is a vector of slope parameters $\kappa \in[0,1]^N$. The $x$-feasible region of \eqref{eq:LA} is defined as \begin{align}\label{eq:feasi_set_LA}
    \mathcal{X}_\text{LA}(\bm{\kappa}) \coloneqq
    \left\{
        \bm{x} \in \mathcal{X} \mid 
        \exists s, \bm{r} :~\eqref{eq:LA:c1}\mbox{--}\eqref{eq:LA:c3}
    \right\}.
\end{align}

\begin{lemma} [Lemma 3 in \cite{zhou2024strengthenedfasterlinearapproximation}]\label{lemma:2}
    For any $\bm{x} \in \mathcal{X}_\text{LA}(\bm{\kappa})$, there exists ($\bm{r}$, $s$) such that $s$ is equal to the $(k+1)$-th smallest value amongst $\{ \widehat{\mathrm{dist}}(\bm{\xi}_i, \mathcal{S}(\bm{x}))\}_{i\in [N]}$ and the constraints in~\eqref{eq:feasi_set_LA}, namely~\eqref{eq:LA:c1}--\eqref{eq:LA:c3}, are satisfied.
\end{lemma}
\subsection{Strengthened and Faster Linear Approximation (SFLA) \cite{zhou2024strengthenedfasterlinearapproximation}}
Given the definition of $q_p$ as the $(k+1)$-th smallest value of $\bm{b}_p^\top \bm{\xi}_{i}$ in~\eqref{eq:def:qp}, the index set $[N]_p$ with only $k$ elements is introduced with this $q_p$ in SFLA \cite{zhou2024strengthenedfasterlinearapproximation} where:
\begin{equation}
   [N]_p \coloneqq \left\{ i\in [N] \mid  \bm{b}_p^\top \bm{\xi}_i < q_p \right\},\quad \forall p\in [P]. \label{eq:def:npset}
\end{equation}
With the approximated distance function $\widehat{\mathrm{dist}}(\bm{\xi}_i, \mathcal{S}(\bm{x}))$ introduced in Eq.~\eqref{eq:dist_approx}, the SFLA of RHS-WDRJCCPs~\eqref{constr:wdrjcc} proposed in \cite{zhou2024strengthenedfasterlinearapproximation} with reduced $[N]_p$ set is:
\begin{subequations}\label{eq:SFLA}
    \begin{align}
        & s \geq 0, \bm{r} \geq \bm{0}, \label{eq:SFLA:c1} \\
        & \epsilon N s - \sum\limits_{i\in [N]}r_i \geq \theta N, \label{eq:SFLA:c2} \\
        & \kappa_i \frac{ \bm{b}_p^\top \bm{\xi}_i + d_p - \bm{a}_p^\top \bm{x}}{\| \bm{b}_p \|_*} \geq s - r_i ,\quad i \in [N]_p\color{black},\ p\in [P], \label{eq:SFLA:c3}\\
        & \frac{q_p + d_p - \bm{a}_p^\top \bm{x}}{\| \bm{b}_p \|_*} \geq s ,\quad p\in [P],  \label{eq:SFLA:c4}
    \end{align}
\end{subequations}
where the $x$-feasible region of \eqref{eq:SFLA} is defined as \begin{align}\label{eq:feasi_set_SFLA}
    \mathcal{X}_\text{SFLA}(\bm{\kappa}) \coloneqq
    \left\{
        \bm{x} \in \mathcal{X} \mid 
        \exists s, \bm{r} :~\eqref{eq:SFLA:c1}\mbox{--}\eqref{eq:SFLA:c4}
    \right\},
\end{align}
This reduced set $[N]_p$ in SFLA may introduce numerical issues for our bilevel problem because it is based on the ranking of $\{\bm{b}_p^\top \bm{\xi}_i\}_{i \in [N]}$. Specifically, when $\bm{b}_p$ depends on the leader's variables, additional disjunctive constraints (typically modelled via big-M constraints) are required to represent the bindingness of constraint~\eqref{eq:SFLA:c3}. These disjunctive constraints in the lower-level problem are then reformulated as complementarity conditions and subsequently linearized, typically using big-M constraints again. This dual use of big-M parameters can lead to numerical instability.

\section{proof of Theorem~\ref{prop:1}}\label{section:proof_prop1}
\begin{IEEEproof}
it is equivalent to show that 
\begin{equation}
\begin{aligned}
    &\left\{
        \bm{x} \in \mathcal{X} \mid 
        \exists s, \bm{r} :~\eqref{eq:exact1:c1}\mbox{--}\eqref{eq:exact1:c3}
    \right\} \\
    &\qquad\qquad\qquad\qquad= \left\{
        \bm{x} \in \mathcal{X} \mid \exists s, \bm{r} :~\eqref{eq:exact_stren:c1}\mbox{--}\eqref{eq:exact_stren:c2}
        \right\}
\end{aligned}
\end{equation}
By~\eqref{eq:exact_stren:c1}, we know that the set on the right-hand side is contained in the set on the
left-hand side. Let us show that the reverse direction also holds. For notational convenience, for every $p\in [P]$, we let $f_{i,p}(\bm{x})\coloneqq \frac{ \bm{b}_p^\top \bm{\xi}_i + d_p - \bm{a}_p^\top \bm{x}}{\| \bm{b}_p \|_*}$ and $f^*_{p}(\bm{x})\coloneqq \frac{q_p + d_p - \bm{a}_p^\top \bm{x}}{\| \bm{b}_p \|_*}$. Therefore, the distance function~\eqref{eq:dist_anal} is reformulated as 
\begin{align} \label{eq:dist_ref}
    \mathrm{dist}\left(\bm{\xi}_i, \mathcal{S}(\bm{x})\right)
        = \left(\min_{p\in [P]} f_{i,p}(\bm{x}) \right)^+ = \min_{p\in [P]} \left(f_{i,p}(\bm{x})\right)^+,
\end{align}
Without loss of generality, assume that we have an ordering $\mathrm{dist}(\bm{\xi}_1, \mathcal{S}(\bm{x}))\leq \cdots \leq \mathrm{dist}(\bm{\xi}_N, \mathcal{S}(\bm{x}))$. Then denote the ($k+1$)-th smallest distance value as $d^*(\bm{x}) \coloneqq \mathrm{dist}(\bm{\xi}_{k+1}, \mathcal{S}(\bm{x}))$. By Lemma~\ref{lemma:1}, take $(s,\bm{r})$ to be a solution that satisfies $s = d^*(\bm{x})$, $r_i = \max\{0, s - \text{dist}(\bm{\xi}_i, S(\bm{x}))\}$.

Now we consider two cases:
\begin{enumerate}
    \item if $\min_{p\in [P]} f_{k+1,p}(\bm{x})\leq 0$, then $s = d^*(\bm{x})=\mathrm{dist}(\bm{\xi}_{k+1}, \mathcal{S}(\bm{x}))=\left(\min_{p\in [P]} f_{k+1,p}(\bm{x}) \right)^+=0$. However, given that the radius $\theta >0$ and $\bm{r}\geq 0$, constraint~\eqref{eq:exact1:c2} cannot be satisfied if $s =0$. Thus, we cannot have $\min_{p\in [P]} f_{k+1,p}(\bm{x})\leq 0$.
    \item if $\min_{p\in [P]} f_{k+1,p}(\bm{x})> 0$, due to the definition in~\eqref{eq:dist_ref}, we have 
\begin{align*}
    0 < \min_{p\in[P]} f_{k+1,p}(\bm{x}) = d^*(\bm{x}),
\end{align*}
Recall the ordering $\mathrm{dist}(\bm{\xi}_{1}, \mathcal{S}(\bm{x}))\leq \cdots \leq \mathrm{dist}(\bm{\xi}_N, \mathcal{S}(\bm{x}))$, for $i\geq k+1$, we have 
\begin{align*}
     d^*(\bm{x}) \leq \text{dist}(\bm{\xi}_i, S(\bm{x})) = \min_{p'\in[P]} f_{i,p'}(\bm{x}) \leq f_{i,p}(\bm{x}).
\end{align*}
It means that there are at least $N-k$ indices $i$ that satisfy $d^*(\bm{x})\leq f_{i,p}(\bm{x})$. Recall the definition of $q_p$ is the $(k+1)$-th smallest among $\{\bm{b}_p^\top \bm{\xi}_{\sigma(i)}\}_{i \in [N]}$, i.e., $q_p \coloneqq \bm{b}_p^\top \bm{\xi}_{\sigma(k+1)}$, $f_p^*$ is the $(k+1)$-th smallest value amongst $\{f_{i,p}(\bm{x})\}_{i\in[N]}$. Therefore, we have
\begin{align*}
     s=d^*(\bm{x}) \leq f_p^*(\bm{x})=\frac{q_p + d_p - \bm{a}_p^\top \bm{x}}{\| \bm{b}_p \|_*},
\end{align*}
\end{enumerate}
Thus, the claim follows.\end{IEEEproof}

\section{Theoretical Properties of SLA}\label{section:SLA_property}
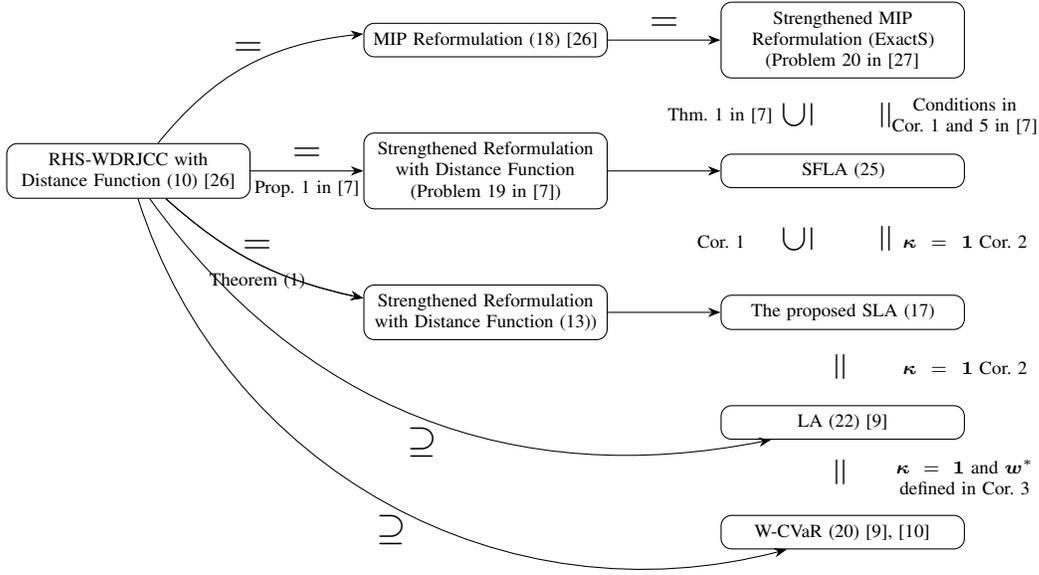
\begin{figure*}
   \centering
   \begin{tikzpicture}[ node distance=1cm and 1.5cm, font=\scriptsize,
                         every node/.style={draw, text width=3cm, align=center, rounded corners},
                         arrow/.style={-{Stealth}}
                       ]
       \node (rhs) {RHS-WDRJCC with Distance Function~\eqref{eq:exact1}~\cite{2022_OR_exact_DRO}};
       \node (strengthenedR) [right=of rhs] {Strengthened Reformulation with Distance Function (Problem 19 in \cite{zhou2024strengthenedfasterlinearapproximation})};
       \node (strengthenedLA) [below=of strengthenedR] {Strengthened Reformulation with Distance Function \eqref{eq:exact_stren})};
       \node (milp) [above=of strengthenedR] {MIP Reformulation~\eqref{eq:exact_MILP}~\cite{2022_OR_exact_DRO}};
       \node (strengthenedM) [right=of milp] {Strengthened MIP Reformulation (ExactS) (Problem 20 in~\cite{exact_milp_strengthened}};
       \node (sfla) [right=of strengthenedR] {SFLA~\eqref{eq:feasi_set_SFLA}};
       \node (sla) [right=of strengthenedLA] {The proposed SLA~\eqref{eq:feasi_set_SLA}};
       \node (la) [below=of sla] {LA~\eqref{eq:feasi_set_LA}~\cite{chen2023approximations}};
       \node (wcvar) [below=of la] {W-CVaR~\eqref{eq:feasi_set_wcvar}~\cite{mohajerin2018data,chen2023approximations}};
       
       \draw[arrow, bend left=23] (rhs) to node[midway, above, draw=none] {{\Large $=$}} (milp);
       \draw[arrow] (rhs) to node[midway, above, draw=none] {{\Large $=$}} (strengthenedR);
       \path[arrow] (rhs) to node[midway, below, draw=none] {Prop.~1 in \cite{zhou2024strengthenedfasterlinearapproximation}} (strengthenedR);
       \draw[arrow, bend right=15] (rhs) to node[midway, above, draw=none] {{\Large $=$}} (strengthenedLA);
       \draw[arrow, bend right=15] (rhs) to node[midway, below, draw=none] {Theorem~\eqref{prop:1} } (strengthenedLA);
       \draw[arrow, bend right=33] (rhs) to node[midway, below, draw=none] {\Large$\supseteq$} (la);
       \draw[arrow, bend right=44] (rhs) to node[midway, below, draw=none] {\Large$\supseteq$} (wcvar);
       \draw[arrow] (strengthenedR) -- (sfla);
       \draw[arrow] (strengthenedLA) -- (sla);
       \draw[arrow] (milp) to node[midway, above, draw=none] {\Large$=$} (strengthenedM);
       \node at (9.5,-4) [draw = none] {\rotatebox{90}{\Large$=$}};
       \node at (9.5,-2.6) [draw = none] {\rotatebox{90}{\Large$=$}};
       \node at (10.1,-0.9) [draw = none] {\rotatebox{90}{\Large$=$}};
       \node at (10.1,0.75) [draw = none] {\rotatebox{90}{\Large$=$}};
       \path[arrow] (wcvar) to node[right, draw=none] {$\bm{\kappa} = \bm{1}$ and ${\bm{w}}^*$ defined in Cor.~\ref{cor:wcvar}} (la);
       \path[arrow] (la) to node[right, draw=none] {$\bm{\kappa} = \bm{1}$ Cor.~\ref{cor:2}} (sla);
       \path[arrow] (la) to node[right, draw=none] {} (sla);
       \path[arrow] (sfla) to node[left, draw=none] {Cor.~\ref{cor:1}} (sla);
       \path[arrow] (sfla) to node[right, draw=none] {$\bm{\kappa} = \bm{1}$ Cor.~\ref{cor:2}} (sla);
       \path[arrow] (sfla) to node[right, draw=none] {Conditions in Cor.~1 and~5 in \cite{zhou2024strengthenedfasterlinearapproximation}} (strengthenedM);
       \node at (8.9,-0.9) [draw = none] {\rotatebox{90}{\Large$\subseteq$}};
       \node at (8.9,0.75) [draw = none] {\rotatebox{90}{\Large$\subseteq$}};
       \path[arrow] (sfla) to node[left, draw=none] {Thm.~1 in \cite{zhou2024strengthenedfasterlinearapproximation}} (strengthenedM);
   \end{tikzpicture}
   \caption{The derivation flow and the formulation comparisons of the $\bm{x}$-feasible region defined by different reformulation of RHS-WDRJCC~\eqref{constr:wdrjcc}.}
\label{fig:flowchart}
\end{figure*}
In this section, we will provide the theoretical properties of the proposed SLA set, $\mathcal{X}_\text{SLA}$. In addition, Fig.~\ref{fig:flowchart}
displays the derivation process and a comparison of different formulations in the $x$-feasible region. We first demonstrate that the proposed SLA provides an inner approximation to the exact RHS-WDRJCC~\eqref{eq:exact1}:
\begin{theorem}\label{theorem:1}
We have $ \mathcal{X}_\text{SLA}(\bm{\kappa})\subseteq\mathcal{X}_\text{Exact}$ for every $\kappa \in[0,1]^N$.
\end{theorem}
\begin{IEEEproof}
    Given a vector of slope parameter $\kappa \in[0,1]^N$, we have
    \begin{align*}
        \left( \frac{\bm{b}_p^\top \bm{\xi}_i + d_p - \bm{a}_p^\top \bm{x}}{\| \bm{b}_p \|_*} \right)^+\geq \kappa_i \frac{ \bm{b}_p^\top \bm{\xi}_i + d_p - \bm{a}_p^\top \bm{x}}{\| \bm{b}_p \|_*}. 
    \end{align*} Thus, the claim follows.
\end{IEEEproof}
We now verify that the inequalities~\eqref{eq:exact_stren:c2} can be used to strengthen the formulation of LA~\eqref{eq:LA}. 
\begin{proposition}\label{prop:2}
We have $\mathcal{X}_\text{LA}(\bm{\kappa}) =\mathcal{X}_\text{SLA}(\bm{\kappa})$.
\end{proposition}
\begin{IEEEproof}
    it is sufficient to show that 
\begin{equation}
\begin{aligned}
    &\left\{
        \bm{x} \in \mathcal{X} \mid 
        \exists s, \bm{r} :~\eqref{eq:LA:c1}\mbox{--}\eqref{eq:LA:c3}
    \right\} \\ &\qquad\qquad\qquad\qquad= \left\{
        \bm{x} \in \mathcal{X} \mid \exists s, \bm{r} :~\eqref{eq:SLA:c1}\mbox{--}\eqref{eq:SLA:c4}
        \right\}
\end{aligned}
\end{equation}
By~\eqref{eq:SLA:c1}--\eqref{eq:SLA:c3}, we know that the set on the right-hand side is contained in the set on the
left-hand side. Let us show that the reverse direction also holds. For notational convenience, for every $p\in [P]$, we let $f_{i,p}(\bm{x})\coloneqq \frac{ \bm{b}_p^\top \bm{\xi}_i + d_p - \bm{a}_p^\top \bm{x}}{\| \bm{b}_p \|_*}$ and $f^*_{p}(\bm{x})\coloneqq \frac{q_p + d_p - \bm{a}_p^\top \bm{x}}{\| \bm{b}_p \|_*}$. Therefore, the distance function~\eqref{eq:dist_approx} is reformulated as 
\begin{align} \label{eq:dist_ref_hat}
    \widehat{\text{dist}}\left(\bm{\xi}_i, \mathcal{S}(\bm{x})\right)
        = \kappa_i\left(\min_{p\in [P]} f_{i,p}(\bm{x}) \right)^+ = \kappa_i\min_{p\in [P]} \left(f_{i,p}(\bm{x})\right)^+,
\end{align}
Without loss of generality, assume that we have an ordering $\widehat{\text{dist}}(\bm{\xi}_1, \mathcal{S}(\bm{x}))\leq \cdots \leq \widehat{\text{dist}}(\bm{\xi}_N, \mathcal{S}(\bm{x}))$. Then denote the ($k+1$)-th smallest distance value as $\hat{d}^*(\bm{x}) \coloneqq \widehat{\text{dist}}(\bm{\xi}_{k+1}, \mathcal{S}(\bm{x}))$. By Lemma~\ref{lemma:2}, take $(s,\bm{r})$ to be a solution that satisfies $s = \hat{d}^*(\bm{x})$, $r_i = \max\{0, s - \widehat{\text{dist}}(\bm{\xi}_i, S(\bm{x}))\}$. Similar to the proof of Theorem~\ref{prop:1} in Appendix~\ref{section:proof_prop1}, we cannot have $\min_{p\in [P]} f_{k+1,p}(\bm{x})\leq 0$. Then for any $i\geq k+1$, $p\in [P]$, due to the definition in~\eqref{eq:dist_ref_hat}, we have 
\begin{align*}
   \hat{d}^*(\bm{x}) \leq \widehat{\text{dist}}(\bm{\xi}_i, S(\bm{x})) = \kappa_i\min_{p'\in[P]} f_{i,p'}(\bm{x}) \leq \kappa_i f_{i,p}(\bm{x}).
\end{align*}
It means that there are at least $N-k$ indices $i$ that satisfy $\hat{d}^*(\bm{x})\leq \kappa_i f_{i,p}(\bm{x})$. Recall the definition of $q_p$ is the $(k+1)$-th smallest among $\{\bm{b}_p^\top \bm{\xi}_{\sigma(i)}\}_{i \in [N]}$, i.e., $q_p \coloneqq \bm{b}_p^\top \bm{\xi}_{\sigma(k+1)}$, $f_p^*$ is the $(k+1)$-th smallest value amongst $\{f_{i,p}(\bm{x})\}_{i\in[N]}$. Given that $\kappa \in[0,1]^N$, we have
\begin{align*}
     s=\hat{d}^*(\bm{x}) \leq \kappa_i f_p^*(\bm{x})\leq f_p^*(\bm{x}) =\frac{q_p + d_p - \bm{a}_p^\top \bm{x}}{\| \bm{b}_p \|_*} ,
\end{align*}
Thus, the claim follows.
\end{IEEEproof}

\noindent We then present the following corollaries to establish the relationships between different approximation schemes.

\begin{corollary}
\label{cor:1}
    The following relationship holds: 
    \[
    \mathcal{X}_\text{LA}(\bm{\kappa})=\mathcal{X}_\text{SLA}(\bm{\kappa})\subseteq\mathcal{X}_\text{SFLA}(\bm{\kappa}).
    \]
\end{corollary}
\begin{IEEEproof}
    Based on Prop.~\ref{prop:2}, we have $\mathcal{X}_\text{LA}(\bm{\kappa})=\mathcal{X}_\text{SLA}(\bm{\kappa})$. We now demonstrate that $\mathcal{X}_\text{SLA}(\bm{\kappa})\subseteq\mathcal{X}_\text{SFLA}(\bm{\kappa})$. This follows from Theorem 2 in~\cite{zhou2024strengthenedfasterlinearapproximation}, which establishes that $\mathcal{X}_\text{LA}(\bm{\kappa})\subseteq\mathcal{X}_\text{SFLA}(\bm{\kappa})$, coupled with the previously established equality $\mathcal{X}_\text{LA}(\bm{\kappa})=\mathcal{X}_\text{SLA}(\bm{\kappa})$. Furthermore, a direct comparison of the constraints reveals that~\eqref{eq:LA:c3} in SLA introduces $NP-\lfloor \epsilon N \rfloor P$ additional valid inequalities compared to~\eqref{eq:SFLA:c3} in SFLA.\end{IEEEproof}

\begin{corollary}
\label{cor:2}
    Under the specific parameter setting of $\bm{\kappa} = \bm{1}$, the feasible regions of the LA, SLA, and SFLA are equivalent:
    \[
    \mathcal{X}_\text{LA}(\bm{1})=\mathcal{X}_\text{SLA}(\bm{1})=\mathcal{X}_\text{SFLA}(\bm{1}).
    \]
\end{corollary}
\begin{IEEEproof}
    The equivalence of $\mathcal{X}_\text{LA}(\bm{1})=\mathcal{X}_\text{SFLA}(\bm{1})$ is demonstrated in Corollary 3 of \cite{zhou2024strengthenedfasterlinearapproximation}. This, combined with the established equality $\mathcal{X}_\text{LA}(\bm{\kappa})=\mathcal{X}_\text{SLA}(\bm{\kappa})$, proves the corollary.\end{IEEEproof}

\noindent Furthermore, the proposed SLA formulation, $\mathcal{X}_\text{SLA}(\bm{\kappa})$, becomes equivalent to W-CVaR in~\eqref{eq:CVaR} under a specific choice of $\bm{w}$ and $\bm{\kappa} = \bm{1}$:

\begin{corollary}\label{cor:wcvar}
    For a particular weighting vector $\bm{w}^*$, the following equivalence holds:
    \[
    \mathcal{X}_\text{SLA}(\bm{1}) = \mathcal{X}_\text{W-CVaR}(\bm{w}^*),
    \]
    where $\bm{w}^*=[w^*_p]_{p\in [P]}$ and $w^*_p = \frac{\| \bm{b}_p \|_*^{-1}}{\sum_{l \in [P]} \| \bm{b}_l \|_*^{-1}}$.
\end{corollary}
\begin{IEEEproof}
    The result follows from \cite{chen2023approximations}, where this equivalence is established for $\mathcal{X}_\text{LA}(\bm{\kappa})$~\eqref{eq:feasi_set_LA} from an optimal value perspective. Given Prop.~\ref{prop:2}, the same conclusion extends to the proposed SLA formulation, $\mathcal{X}_\text{SLA}(\bm{1})$.\end{IEEEproof}

\noindent Finally, under certain conditions, the proposed SLA formulation $\mathcal{X}_\text{SLA}(\bm{1})$ becomes exact when $\bm{\kappa} = \bm{1}$:

\begin{corollary}\label{cor:5}
    The equivalence $\mathcal{X}_\text{SLA}(\bm{1}) = \mathcal{X}_\text{Exact}$ holds under either of the following conditions:
    \begin{enumerate}[label={(\arabic*)}, leftmargin=3em]
        \item $\bm{\xi}_i \in \mathcal{S}(\bm{x})$ for all $i \in [N]$ and $\bm{x} \in \mathcal{X}_\text{Exact}$.
        \item $\epsilon \leq 1/N$.
    \end{enumerate}
\end{corollary}
\begin{IEEEproof}
    This follows from Corollary 5 presented in~\cite{zhou2024strengthenedfasterlinearapproximation}, combined with the previously established equivalence $\mathcal{X}_\text{SLA}(\bm{1})=\mathcal{X}_\text{SFLA}(\bm{1})$ under the parameter setting $\bm{\kappa} = \bm{1}$in Corollary~\ref{cor:2}.\end{IEEEproof}

\section{Reformulation of the bilevel problem under SLA}\label{section:SLA_reformulation}
This section presents the SLA-based reformulation of the bilevel problem. We begin with the SLA reformulation of the RHS-WDRJCC line flow constraints~\eqref{eq:JCC:bilevel}, followed by the KKT conditions. Finally, we introduce the linearization techniques applied to the resulting single-level problem and show the final problem formulation.
\subsection{SLA Reformulation}\label{section:SLA_reformulation1}

For ease of reading, we denote the following term for every $t \in \mathcal{T},s\in \mathcal{S},l\in\mathcal{L}$ as 
\begin{IEEEeqnarray}{rCl}
\Lambda_{t,s,l} := \sum_{c\in\mathcal{C}}\sum_{b\in\mathcal{B}} s_{b,l,c} o_{t,c} \bigg(
&& \sum_{k\in\Omega_{t,s,b}^{G}} g_{t,s,k,b} 
+ \sum_{k\in\Omega_{t,s,b}^{W}} p^{w,sch}_{t,s,k,b} \nonumber \\
&& {} - \sum_{k\in\Omega_{t,s,b}^{D}} d_{t,s,k,b} \bigg)
\label{eq:Xi}
\end{IEEEeqnarray}

The RHS-WDRJCC~\eqref{eq:JCC:bilevel} under the proposed SLA~\eqref{eq:SLA} for every $t \in \mathcal{T},s\in \mathcal{S}$ is 
\begin{subequations}
\label{eq:bilevel:SLA}
\begin{align}
& -u_{t,s} \leq 0,  -v_{t,s,i} \leq 0, : \mu^{u,\min}_{t,s},\mu^{v,\min}_{t,s,i} \quad \forall i\in[N], \label{eq:SLAa}\\
& -(\epsilon N u_{t,s} - \sum\limits_{i\in [N]}v_{t,s,i})\leq -\theta N,  :\mu^{1}_{t,s} \label{eq:SLAb}  \\
&u_{t,s}-v_{t,s,i} + \kappa_i\Lambda_{t,s,l}\leq \kappa_i (\mathcal{F}^{0}_{l}(1+\sum_{\hat{t} \leq t}\sum_{j\in\mathcal{J}}j b^R_{\hat{t},l,j})+\notag\\
&\sum_{m \in \mathcal{M}} m z^P_{t,l,m} \overline{F}^{C}_{l})-\kappa_i \sum_{c\in\mathcal{C}}\sum_{b\in\mathcal{B}}\sum_{k\in\Omega_{t, s,b}^{W}} s_{b,l,c}o_{t,c} \tilde{e}_{t,s,k,b} \notag\\
&:\mu_{t,s,l,i}^{2}\quad \forall l\in\mathcal{L},\forall i\in [N],\label{eq:SLAc}\\
&  u_{t,s}-v_{t,s,i} - \kappa_i\Lambda_{t,s,l} \leq \kappa_i (\mathcal{F}^{0}_{l}(1+\sum_{\hat{t} \leq t}\sum_{j\in\mathcal{J}}j b^R_{\hat{t},l,j})+\notag\\
&\sum_{m \in \mathcal{M}} m z^P_{t,l,m} \overline{F}^{C}_{l})+\kappa_i\sum_{c\in\mathcal{C}}\sum_{b\in\mathcal{B}}\sum_{k\in\Omega_{t, s,b}^{W}} s_{b,l,c}o_{t,c} \tilde{e}_{t,s,k,b}\notag\\
&:\mu_{t,s,l,i}^{3} \quad\forall l\in\mathcal{L},\forall i\in [N],\label{eq:SLAd}\\
& u_{t,s}+ \Lambda_{t,s,l}\leq  (\mathcal{F}^{0}_{l}(1+\sum_{\hat{t} \leq t}\sum_{j\in\mathcal{J}}j b^R_{\hat{t},l,j})+\sum_{m \in \mathcal{M}} m z^P_{t,l,m} \overline{F}^{C}_{l}) \notag\\
&+   \sum_{c\in\mathcal{C}}  o_{t,c} Q^{\max}_{t,s,l,c},:\mu_{t,s,l}^{4}\quad  \forall l\in\mathcal{L}, \label{eq:SLAe}\\
& u_{t,s}- \Lambda_{t,s,l}\leq  (\mathcal{F}^{0}_{l}(1+\sum_{\hat{t} \leq t}\sum_{j\in\mathcal{J}}j b^R_{\hat{t},l,j})+\sum_{m \in \mathcal{M}} m z^P_{t,l,m} \overline{F}^{C}_{l}) \notag\\
&+   \sum_{c\in\mathcal{C}}   o_{t,c} Q^{\min}_{t,s,l,c},:\mu_{t,s,l}^{5}\quad  \forall l\in\mathcal{L}. \label{eq:SLAf}
\end{align}
\end{subequations}
\subsection{KKT Conditions under SLA}

The proposed SLA, along with SFLA, LA and W-CVaR approximations, leads to convex reformulation of the RHS-WDRJCC in the lower-level problem. Consequently, by substituting the lower level problem with its KKT conditions, the bilevel problem can be effectively transformed into a single-level Mathematical Program with Equilibrium Constraints (MPEC). Firstly, we define the following term under SLA approximation as
\begin{equation}
\begin{aligned}
& \Lambda^{SLA}_{t,s,b} := \sum_{i\in[N]} \kappa_i \sum_{c\in\mathcal{C}}\sum_{l\in\mathcal{L}}s_{b,l,c}o_{t,c}(\mu^2_{t,s,l,i}-\mu^3_{t,s,l,i})\\
& +  \sum_{c\in\mathcal{C}}\sum_{l\in\mathcal{L}}s_{b,l,c}o_{t,c}(\mu^4_{t,s,l}-  \mu^5_{t,s,l} )
\end{aligned}
\label{eq:Xi_SLA_original}
\end{equation}
then KKT conditions under SLA are
\begin{subequations}
\begin{align}
& Eqs.~\eqref{eq:lowerb},\eqref{eq:lowere} \label{KKT:SLAa},\\
& \hat{c}_{t,s,k,b}^{g} -\pi_{t,s} + \varphi^{G,\max}_{t,s,k,b}-\varphi^{G,\min}_{t,s,k,b}+ \Lambda^{SLA}_{t,s,b} =0,\notag\\
&\forall t \in \mathcal{T},  s\in \mathcal{S},  k \in \Omega_{t, s,b}^{G}, b \in \mathcal{B}, \label{KKT:SLAb}\\
& - \hat{c}_{t,s,k,b}^{d} + \pi_{t,s} +  \varphi^{D,\max}_{t,s,k,b}-\varphi^{D,\min}_{t,s,k,b}-\Lambda^{SLA}_{t,s,b} =0,\notag\\
&\forall t \in \mathcal{T},  s\in \mathcal{S},  k \in \Omega_{t, s,b}^{D},b \in \mathcal{B}, \label{KKT:SLAc}\\
&  \hat{c}_{t,s,k,b}^{w}-\pi_{t,s} + \varphi_{t,s,k,b}^{W,sch}+ \Lambda^{SLA}_{t,s,b} =0,\notag\\
&\forall t \in \mathcal{T},  s\in \mathcal{S},  k \in \Omega_{t, s,b}^{W}, b \in \mathcal{B},\label{KKT:SLAd}\\
& c_{t,s,k,b}^{cur} + \varphi_{t,s,k,b}^{W,sch} +\varphi^{W,cur,\max}_{t,s,k,b}-\varphi^{W,cur,\min}_{t,s,k,b} =0 ,\notag\\
&\forall t \in \mathcal{T},  s\in \mathcal{S}, k \in \Omega_{t, s,b}^{W}, b \in \mathcal{B},\label{KKT:SLAe}\\
& -\mu^{v,\min}_{t,s,i}+ \mu^1_{t,s} -\sum_{l\in\mathcal{L}} \mu^2_{t,s,l,i} -\sum_{l\in\mathcal{L}} \mu^3_{t,s,l,i} = 0, \notag\\
&\forall  t \in \mathcal{T}, s \in \mathcal{S}, i\in [N],\label{KKT:SLAf}\\
& -\mu^{u,\min}_{t,s}-\epsilon N \mu^1_{t,s} + \sum_{l\in\mathcal{L}}(\sum_{i\in [N]}\mu^2_{t,s,l,i}+\sum_{i\in [N]}\mu^3_{t,s,l,i}\notag\\
&+\mu^4_{t,s,l} +\mu^5_{t,s,l}) = 0,\quad \forall t \in \mathcal{T}, s \in \mathcal{S},\label{KKT:SLAg}\\
& 0 \leq (g_{t,s,k,b}-g^{\min}_{t,s,k,b})\perp \varphi^{G,\min}_{t,s,k,b} \geq 0,\notag\\
&\forall t\in \mathcal{T},  s\in \mathcal{S},  k \in \Omega_{t, s,b}^{G}, b \in \mathcal{B},\label{KKT:SLAh}\\
& 0 \leq (g^{\max}_{t,s,k,b}-g_{t,s,k,b})\perp \varphi^{G,\max}_{t,s,k,b} \geq 0,\notag\\
&\forall t\in \mathcal{T},  s\in \mathcal{S},  k \in \Omega_{t, s,b}^{G}, b \in \mathcal{B},\label{KKT:SLAi}\\
& 0 \leq (d_{t,s,k,b}-d^{\min}_{t,s,k,b})\perp \varphi^{D,\min}_{t,s,k,b} \geq 0 ,\notag\\
&\forall t \in \mathcal{T},  s\in \mathcal{S},  k \in \Omega_{t, s,b}^{D},b \in \mathcal{B}, \label{KKT:SLAj}\\
& 0 \leq (d^{\max}_{t,s,k,b}-d_{t,s,k,b})\perp \varphi^{D,\max}_{t,s,k,b} \geq 0 ,\notag\\
&\forall t\in \mathcal{T},  s\in \mathcal{S},  k \in \Omega_{t, s,b}^{D},b \in \mathcal{B},  \label{KKT:SLAk}\\
& 0 \leq p^{w,cur}_{t,s,k,b}\perp \varphi^{W,cur,\min}_{t,s,k,b} \geq 0 ,\notag\\
&\forall t \in \mathcal{T},  s\in \mathcal{S},  k \in \Omega_{t, s,b}^{W},b \in \mathcal{B}, \label{KKT:SLAl}\\
& 0 \leq (p^{w,fore}_{t,s,k,b}-p^{w,cur}_{t,s,k,b})\perp \varphi^{W,cur,\max}_{t,s,k,b} \geq 0,\notag\\
&\forall t\in \mathcal{T},  s\in \mathcal{S},  k \in \Omega_{t, s,b}^{W},b \in \mathcal{B},  \label{KKT:SLAm}\\
&0 \leq u_{t,s} \perp \mu^{u,\min}_{t,s} \geq 0,  \quad \forall t \in \mathcal{T},  s\in \mathcal{S},\label{KKT:SLAn}\\
&0 \leq v_{t,s,i} \perp \mu^{v,\min}_{t,s,i} \geq 0 ,\quad \forall t\in \mathcal{T},  s\in \mathcal{S},i \in [N],\label{KKT:SLAo}\\
&0 \leq (-\theta N +\epsilon N u_{t,s} - \sum\limits_{i\in [N]}v_{t,s,i}) \perp \mu^1_{t,s} \geq 0 ,\notag\\
&\forall t\in \mathcal{T},  s\in \mathcal{S},\label{KKT:SLAp}\\
&0 \leq\kappa_i  (\mathcal{F}^{0}_{l}(1+\sum_{\hat{t} \leq t}\sum_{j\in\mathcal{J}}j b^R_{\hat{t},l,j})+\sum_{m \in \mathcal{M}} m z^P_{t,l,m} \overline{F}^{C}_{l})\notag\\
&-\kappa_i  \sum_{c\in\mathcal{C}}\sum_{b\in\mathcal{B}}\sum_{k\in\Omega_{t, s,b}^{W}} s_{b,l,c}o_{t,c} \tilde{e}_{t,s,k,b} \notag\\
&-  (u_{t,s}-v_{t,s,i}) - \kappa_i \Lambda_{t,s,l} \perp \mu_{t,s,l,i}^{2} \geq 0, \notag\\
& \forall t \in \mathcal{T}, s\in \mathcal{S},l\in\mathcal{L},i\in [N],\label{KKT:SLAq}\\
&0 \leq\kappa_i  (\mathcal{F}^{0}_{l}(1+\sum_{\hat{t} \leq t}\sum_{j\in\mathcal{J}}j b^R_{\hat{t},l,j})+\sum_{m \in \mathcal{M}} m z^P_{t,l,m} \overline{F}^{C}_{l})\notag\\
&+\kappa_i \sum_{c\in\mathcal{C}} \sum_{b\in\mathcal{B}}\sum_{k\in\Omega_{t, s,b}^{W}} s_{b,l,c}o_{t,c} \tilde{e}_{t,s,k,b} \notag\\
&-  (u_{t,s}-v_{t,s,i}) + \kappa_i \Lambda_{t,s,l} \perp \mu_{t,s,l,i}^{3} \geq 0,\notag\\
& \forall t \in \mathcal{T}, s\in\mathcal{S},l\in\mathcal{L},i\in [N],\label{KKT:SLAr}\\
& 0 \leq  (\mathcal{F}^{0}_{l}(1+\sum_{\hat{t} \leq t}\sum_{j\in\mathcal{J}}j b^R_{\hat{t},l,j})+\sum_{m \in \mathcal{M}} m z^P_{t,l,m} \overline{F}^{C}_{l})\notag\\
& +   \sum_{c\in\mathcal{C}}  o_{t,c} Q^{\max}_{t,s,l,c} -  u_{t,s}- \Lambda_{t,s,l}\perp \mu_{t,s,l}^{4}\geq 0, \notag\\
&\forall t \in \mathcal{T},s \in \mathcal{S},l\in\mathcal{L},\label{KKT:SLAs}\\
& 0 \leq   (\mathcal{F}^{0}_{l}(1+\sum_{\hat{t} \leq t}\sum_{j\in\mathcal{J}}j b^R_{\hat{t},l,j})+\sum_{m \in \mathcal{M}} m z^P_{t,l,m} \overline{F}^{C}_{l})\notag\\
&+ \sum_{c\in\mathcal{C}}  o_{t,c} Q^{\min}_{t,s,l,c} - u_{t,s}+ \Lambda_{t,s,l}\perp\mu_{t,s,l}^{5}\geq 0, \notag\\
&\forall t \in\mathcal{T},s\in\mathcal{S},l\in\mathcal{L}.\label{KKT:SLAt}
\end{align}
\label{KKT:SLA}
\end{subequations}
\subsection{Reformulation}
The resulting single-level optimization model is a MPEC. This model, however, is non-linear. We have four types of nonlinearities and they are in the form of
\begin{enumerate}
    \item complementarity conditions~\eqref{KKT:SLAh}--\eqref{KKT:SLAt} in KKT conditions;
    \item $\pi_{t,s,b}d_{t, s,k,b},
 \pi_{t,s,b}g_{t, s,k,b},\pi_{t,s,b}p_{t, s,k,b}^{w,sch}$ in calculation of merchandising surplus~\eqref{eq:uppere}
    \item $o_{t,c} \mu^2_{t,s,l,i}$, $o_{t,c} \mu^3_{t,s,l,i}$, $o_{t,c} \mu^4_{t,s,l}$ and $o_{t,c} \mu^5_{t,s,l}$ in $\Lambda^{SLA}_{t,s,b}$~\eqref{eq:Xi_SLA_original}, $o_{t,c} g_{t,s,k,b}$, $o_{t,c} p^{w,sch}_{t,s,k,b}$ and $o_{t,c} d_{t,s,k,b}$ in $\Lambda_{t,s,b}$~\eqref{eq:Xi};
    \item $z^R_{t,l}\tau^V_l$ and $z^P_{t,l,m}\tau^V_l$ in calculation of actual bids~\eqref{eq:upperb}--\eqref{eq:upperd}; $z^R_{t,l}\tau^V_l g_{t,s,k,b}$, $z^R_{t,l}\tau^V_l d_{t,s,k,b}$, $z^R_{t,l}\tau^V_l p^{w,sch}_{t,s,k,b}$, $z^P_{t,l,m}\tau^V_l g_{t,s,k,b}$, $z^P_{t,l,m}\tau^V_l d_{t,s,k,b}$, $z^P_{t,l,m}\tau^V_l p^{w,sch}_{t,s,k,b}$ in calculation of volumetric-based network charge revenue~\eqref{eq:upperg} and~\eqref{eq:upperh}.
\end{enumerate}

For the first source of non-linearity which is complementarity conditions~\eqref{KKT:SLAh}--\eqref{KKT:SLAt} within the KKT conditions, each complementarity condition in the form of $0\leq a \perp b \geq 0$ can be readily linearized using a Big-M approach \cite{8036231}. The resulting mixed-integer conditions are $a \geq 0$, $b \geq 0$, $a\leq M*u$ and $b\geq M*(1-u)$. Note that $u$ is an auxiliary binary variable and $M$ is a large enough positive value.

For the second type of non-linearity, we first calculate the LMP $\pi_{t,s,b}$: 
\begin{flalign}
&  \pi_{t,s,b} = \pi_{t,s} - \Lambda^{SLA}_{t,s,b} \label{SLA:LMP}
\end{flalign}
Through dual constraints and complementary slackness, and the definition of actual marginal benefits and costs~\eqref{eq:upperb}--\eqref{eq:upperd}, we have:
\begin{equation}
\begin{aligned}
&\pi_{t,s,b} d_{t, s,k,b}= (\hat{c}_{t,s,k,b}^{d} -  \varphi^{D,\max}_{t,s,k,b}+\varphi^{D,\min}_{t,s,k,b} )d_{t, s,k,b}\\
&=\hat{c}_{t,s,k,b}^{d}d_{t, s,k,b}-  \varphi^{D,\max}_{t,s,k,b}d^{\max}_{t, s,k,b}+\varphi^{D,\min}_{t,s,k,b}d^{\min}_{t, s,k,b} \\
&=-(\sum_{l\in\mathcal{L}^R}z^R_{t,l}\tau^V_{l}\delta_{l,b}+\sum_{l\in\mathcal{L}^P}z^P_{t,l}\tau^V_{l}\delta_{l,b})d_{t, s,k,b}\\
&\quad+c_{t,s, k,b}^{d}d_{t, s,k,b}-  \varphi^{D,\max}_{t,s,k,b}d^{\max}_{t, s,k,b}+\varphi^{D,\min}_{t,s,k,b}d^{\min}_{t, s,k,b}
\end{aligned}  
\label{eq:ms1}
\end{equation}

\begin{equation}
\begin{aligned}
&\pi_{t,s,b} g_{t, s,k,b}= (\hat{c}_{t,s,k,b}^{g} +  \varphi^{G,\max}_{t,s,k,b}-\varphi^{G,\min}_{t,s,k,b})g_{t, s,k,b}\\
&=\hat{c}_{t,s,k,b}^{g}g_{t, s,k,b}+  \varphi^{G,\max}_{t,s,k,b}g^{\max}_{t, s,k,b}-\varphi^{G,\min}_{t,s,k,b}g^{\min}_{t, s,k,b}\\
&=(\sum_{l\in\mathcal{L}^R}z^R_{t,l}\tau^V_{l}\delta_{l,b}+\sum_{l\in\mathcal{L}^P}z^P_{t,l}\tau^V_{l}\delta_{l,b})g_{t, s,k,b}\\
&\quad+c_{t,s, k,b}^{g}g_{t, s,k,b}+  \varphi^{G,\max}_{t,s,k,b}g^{\max}_{t, s,k,b}-\varphi^{G,\min}_{t,s,k,b}g^{\min}_{t, s,k,b}
\end{aligned}  
\label{eq:ms2}
\end{equation}

\begin{equation}
\begin{aligned}
&\pi_{t,s,b} p^{w,sch}_{t, s,k,b}= (\hat{c}_{t,s,k,b}^{w}+\varphi_{t,s,k,b}^{W,sch})p^{w,sch}_{t, s,k,b}\\
&= \hat{c}_{t,s,k,b}^{w}p^{w,sch}_{t, s,k,b}+\varphi_{t,s,k,b}^{W,sch}p^{w,fore}_{t,s,k,b} - \varphi_{t,s,k,b}^{W,sch}p^{w,cur}_{t,s,k,b}\\
&= \hat{c}_{t,s,k,b}^{w}p^{w,sch}_{t, s,k,b}+\varphi_{t,s,k,b}^{W,sch}p^{w,fore}_{t,s,k,b} \\
&\quad+(c_{t,s,k,b}^{cur}  +\varphi^{W,cur,\max}_{t,s,k,b}-\varphi^{W,cur,\min}_{t,s,k,b})p^{w,cur}_{t,s,k,b}\\
&= \hat{c}_{t,s,k,b}^{w}p^{w,sch}_{t, s,k,b}+\varphi_{t,s,k,b}^{W,sch}p^{w,fore}_{t,s,k,b}  \\
&\quad+c_{t,s,k,b}^{cur}p^{w,cur}_{t,s,k,b} +\varphi^{W,cur,\max}_{t,s,k,b}p^{w,fore}_{t,s,k,b}\\
& =(\sum_{l\in\mathcal{L}^R}z^R_{t,l}\tau^V_{l}\delta_{l,b}+\sum_{l\in\mathcal{L}^P}z^P_{t,l}\tau^V_{l}\delta_{l,b})p^{w,sch}_{t, s,k,b} \\
&\quad+\varphi_{t,s,k,b}^{W,sch}p^{w,fore}_{t,s,k,b} +c_{t,s,k,b}^{cur}p^{w,cur}_{t,s,k,b}  \\
&\quad+\varphi^{W,cur,\max}_{t,s,k,b}p^{w,fore}_{t,s,k,b}
\end{aligned}  
\label{eq:ms3}
\end{equation}

\begin{proposition}
The merchandising surplus decreases as a result of the changes in the market curve induced by the imposed volumetric network charges, but this reduction is offset by the increase in the volumetric-based network charge revenue $VC$ collected by the network planner. 
\end{proposition}

For the third type of nonlinearity, we have the following bilinear terms $o_{t,c} \mu^2_{t,s,l,i}$, $o_{t,c} \mu^3_{t,s,l,i}$, $o_{t,c} \mu^4_{t,s,l}$ and $o_{t,c} \mu^5_{t,s,l}$ in $\Lambda^{SLA}_{t,s,b}$. We will introduce auxiliary variables $\eta^{o,2}_{t,s,l,c,i}$, 
$\eta^{o,3}_{t,s,l,c,i}$, 
$\eta^{o,4}_{t,s,l,c}$ 
and $\eta^{o,5}_{t,s,l,c}$ to replace these bilinear terms with big-M constraints, $\forall t \in \mathcal{T},s \in \mathcal{S},l \in\mathcal{L},c \in\mathcal{C}$:
\begin{subequations}
\begin{align}
&  0 \leq \eta^{o,2}_{t,s,l,c,i}\leq M o_{t,c},\quad \forall i\in [N],\label{linearization:oa}\\
& 0 \leq \mu^2_{t,s,l,i}-\eta^{o,2}_{t,s,l,c,i} \leq M(1-o_{t,c} ),\quad \forall i\in [N], \label{linearization:ob}\\
&  0 \leq \eta^{o,3}_{t,s,l,c,i}\leq M o_{t,c},\quad \forall i\in [N],\label{linearization:oc}\\
& 0 \leq \mu^3_{t,s,l,i}-\eta^{o,3}_{t,s,l,c,i} \leq M(1-o_{t,c} ),\quad \forall i\in [N], \label{linearization:od}\\
&  0 \leq \eta^{o,4}_{t,s,l,c}\leq M o_{t,c},\label{linearization:oe}\\
& 0 \leq \mu^4_{t,s,l}-\eta^{o,4}_{t,s,l,c} \leq M(1-o_{t,c} ), \label{linearization:of}\\
&  0 \leq \eta^{o,5}_{t,s,l,c}\leq M o_{t,c},\label{linearization:og}\\
& 0 \leq \mu^5_{t,s,l}-\eta^{o,5}_{t,s,l,c} \leq M(1-o_{t,c} ). \label{linearization:oh}
\end{align}
\label{linearization:o}
\end{subequations}
Then $\Lambda^{SLA}_{t,s,b}$ with auxiliary variables $\eta^{o,2}_{t,s,l,c,i}$, 
$\eta^{o,3}_{t,s,l,c,i}$, 
$\eta^{o,4}_{t,s,l,c}$ 
and $\eta^{o,5}_{t,s,l,c}$ $\forall t\in \mathcal{T},  s\in \mathcal{S},b\in\mathcal{B}$ becomes
\begin{equation}
\begin{aligned}
& \Lambda^{SLA}_{t,s,b} := \sum_{i\in[N]} \kappa_i \sum_{c\in\mathcal{C}}\sum_{l\in\mathcal{L}}s_{b,l,c}(\eta^{o,2}_{t,s,l,c,i}- \eta^{o,3}_{t,s,l,c,i})\\
& +  \sum_{c\in\mathcal{C}}\sum_{l\in\mathcal{L}}s_{b,l,c}(\eta^{o,4}_{t,s,l,c}-  \eta^{o,5}_{t,s,l,c} )
\end{aligned}
\label{eq:SLAterm}
\end{equation}

Similarly, for bilinear terms $o_{t,c} g_{t,s,k,b}$, $o_{t,c} p^{w,sch}_{t,s,k,b}$ and $o_{t,c} d_{t,s,k,b}$ in $\Lambda_{t,s,l}$, we will introduce the auxiliary variable $\eta^{o,g}_{t,s,k,b,c}$, 
$\eta^{o,w}_{t,s,k,b,c}$, 
$\eta^{o,d}_{t,s,k,b,c}$ to replace these bilinear terms with big-M constraints, $\forall t \in \mathcal{T},  s\in \mathcal{S}, b \in \mathcal{B}, c \in \mathcal{C} $:
\begin{subequations}
\begin{align}
&  0 \leq \eta^{o,g}_{t,s,k,b,c}\leq M o_{t,c},\quad  \forall k \in \Omega_{t, s,b}^{G}, \label{linearization:og1}\\
& 0 \leq g_{t,s,k,b}-\eta^{o,g}_{t,s,k,b,c} \leq M(1-o_{t,c} ),\quad \forall k \in \Omega_{t, s,b}^{G}, \label{linearization:og2}\\
&  0 \leq \eta^{o,d}_{t,s,k,b,c}\leq M o_{t,c},\quad \forall  k \in \Omega_{t, s,b}^{D}, \label{linearization:od1}\\
& 0 \leq d_{t,s,k,b}-\eta^{o,d}_{t,s,k,b,c} \leq M(1-o_{t,c} ),\quad \forall  k \in \Omega_{t, s,b}^{D},  \label{linearization:od2}\\
&  0 \leq \eta^{o,w}_{t,s,k,b,c}\leq M o_{t,c},\quad \forall  k \in \Omega_{t, s,b}^{W},  \label{linearization:ow1}\\
& 0 \leq p^{w,sch}_{t,s,k,b}-\eta^{o,w}_{t,s,k,b,c} \leq M(1-o_{t,c} ),\quad \forall k \in \Omega_{t, s,b}^{W}.  \label{linearization:ow2}
\end{align}
\label{linearization:ogwd}
\end{subequations}
Then $\Lambda_{t,s,l}$ with auxiliary variables $\eta^{o,g}_{t,s,k,b,c}$, 
$\eta^{o,w}_{t,s,k,b,c}$, 
$\eta^{o,d}_{t,s,k,b,c}$ $ \forall t\in \mathcal{T},  s\in \mathcal{S},l\in\mathcal{L}$ becomes
\begin{equation}
\begin{aligned}
&\Lambda_{t,s,l} := \sum_{c\in\mathcal{C}}\sum_{b\in\mathcal{B}} (\sum_{k\in\Omega_{t, s,b}^{G}}s_{b,l,c}\eta^{o,g}_{t,s,k,b,c} \\
&+\sum_{k\in\Omega_{t, s,b}^{W}} s_{b,l,c} \eta^{o,w}_{t,s,k,b,c} - \sum_{k\in\Omega_{t, s,b}^{D}} s_{b,l,c}\eta^{o,d}_{t,s,k,b,c}). 
\end{aligned}
\label{eq:SLAXiterm}
\end{equation}

The last type of non-linearity is in the calculation of capacity-based network charge revenue and bids updates, we have bilinear terms $z^R_{t,l}\tau^V_l$ and $z^P_{t,l,m}\tau^V_l$. We then introduce auxiliary non-negative continuous variables $\eta^{R,V}_{t,l} =z^R_{t,l}\tau^V_l$ and $\eta^{P,V}_{t,l,m} =z^P_{t,l,m}\tau^V_l$ and the following constraints $\forall t\in\mathcal{T}$:
\begin{subequations}
\begin{align}
& 0 \leq \eta^{R,V}_{t,l}\leq M z^R_{t,l},\quad \forall l \in\mathcal{L}^R,\\
& 0 \leq \tau^V_l-\eta^{R,V}_{t,l} \leq M(1-z^R_{t,l} ), \quad \forall l \in\mathcal{L}^R, \\
& 0 \leq \eta^{P,V}_{t,l,m}\leq M z^P_{t,l,m}, \quad \forall l \in\mathcal{L}^P,m\in\mathcal{M}\\
& 0 \leq \tau^V_l-\eta^{P,V}_{t,l,m} \leq M(1-z^P_{t,l,m} ), \forall l \in\mathcal{L}^P,m\in\mathcal{M}.
\end{align}  
\label{linearization:tau}
\end{subequations}

It is important to note that the calculation of volume-based network charge revenue in equations~\eqref{eq:upperg} and~\eqref{eq:upperh} involves a non-convex bilinear term, arising from the product of auxiliary variables $\eta^{R,V}_{t,l}$ and $\eta^{P,V}_{t,l}$ with the corresponding cleared quantities $g_{t,s,k,b},d_{t,s,k,b},p^{w,sch}_{t,s,k,b}$. The Gurobi solver can handle the resulting mixed-integer programming with bilinear constraints, which ensures global optimality. Alternatively, linearization techniques such as binary expansion of tariffs or McCormick relaxations may be applied, but are left for future work.

\subsection{Final Single-level Problem under SLA}
The final mixed-integer problem with bilinear constraints under SLA is 
\begin{fleqn}
\begin{subequations}
\begin{align}
    &\max_{\Xi^{SL}_{SLA}} \sum_{t \in \mathcal{T}} \frac{1}{(1+r)^{t-1}}\bigg[\sum_{s \in \mathcal{S}}\sum_{b \in \mathcal{B}}\Psi \bigg(\sum_{k \in \Omega_{t, s,b}^{D}} c_{t,s, k,b}^{d} d_{t,s, k,b} \notag\\
    &-\sum_{k \in \Omega_{t, s,b}^{G}} c_{t,s,k,b}^{g} g_{t,s,k,b}-\sum_{k \in \Omega_{t, s,b}^{W}} c^{cur}_{t,s,k,b}p^{w,cur}_{t,s,k,b}\bigg)-C_t \bigg] \label{eq:MIQCPSLA_obj}\end{align}
    s.t.
    \begin{align}
       &\hat{c}_{t,s, k,b}^{d}=c_{t,s, k,b}^{d}-(\sum_{l\in\mathcal{L}^R}\eta^{R,V}_{t,l}\delta_{l,b}+\sum_{l\in\mathcal{L}^P}\sum_{m\in\mathcal{M}}\eta^{P,V}_{t,l,m}\delta_{l,b}), \notag\\
    & \forall t\in \mathcal{T},  s\in \mathcal{S},  k \in \Omega_{t, s,b}^{D}, b \in \mathcal{B},\label{eq:MIQCPSLAb} \end{align}
     \vspace{-15pt}
     \begin{align}
       &\hat{c}_{t,s, k,b}^{g}=c_{t,s, k,b}^{g}+(\sum_{l\in\mathcal{L}^R}\eta^{R,V}_{t,l}\delta_{l,b}+\sum_{l\in\mathcal{L}^P}\sum_{m\in\mathcal{M}}\eta^{P,V}_{t,l,m}\delta_{l,b}), \notag\\
    & \forall t\in \mathcal{T},  s\in \mathcal{S},  k \in \Omega_{t, s,b}^{G}, b \in \mathcal{B}, \label{eq:MIQCPSLAc} \end{align}
    \vspace{-15pt}
    \begin{align}
       &\hat{c}_{t,s, k,b}^{w}=(\sum_{l\in\mathcal{L}^R}\eta^{R,V}_{t,l}\delta_{l,b}+\sum_{l\in\mathcal{L}^P}\sum_{m\in\mathcal{M}}\eta^{P,V}_{t,l,m}\delta_{l,b}), \notag\\
    &  \forall t\in \mathcal{T},  s\in \mathcal{S},  k \in \Omega_{t, s,b}^{W}, b \in \mathcal{B},  \label{eq:MIQCPSLAd} \end{align}
    \vspace{-15pt}
    \begin{align}
       &VC^{R}_{t,l} = \Psi\sum_{s \in \mathcal{S}}\sum_{b \in \mathcal{B}}(\sum_{k \in \Omega_{t, s,b}^{G}}\delta_{l,b}\eta^{R,V}_{t,l}g_{t, s,k,b} + \notag\\
    &\sum_{k \in \Omega_{t, s,b}^{D}}\delta_{l,b}\eta^{R,V}_{t,l}d_{t, s,k,b} + \sum_{k \in \Omega_{t, s,b}^{W}}\delta_{l,b}\eta^{R,V}_{t,l}p^{w,sch}_{t, s,k,b}),\notag\\
    & \forall  t\in\mathcal{T}, l\in\mathcal{L}^R,\label{eq:MIQCPSLAe}\end{align}
    \vspace{-15pt}
    \begin{align}
           &VC^{P}_{t,l} = \Psi\sum_{s \in \mathcal{S}}\sum_{b \in \mathcal{B}}\sum_{m \in \mathcal{M}}(\sum_{k \in \Omega_{t, s,b}^{G}}\delta_{l,b}\eta^{P,V}_{t,l,m}g_{t, s,k,b} \notag\\
    &+ \sum_{k \in \Omega_{t, s,b}^{D}}\delta_{l,b}\eta^{P,V}_{t,l,m}d_{t, s,k,b} + \sum_{k \in \Omega_{t, s,b}^{W}}\delta_{l,b}\eta^{P,V}_{t,l,m}p^{w,sch}_{t, s,k,b}), \notag\\
    &\forall  t\in\mathcal{T},l\in\mathcal{L}^P, \label{eq:MIQCPSLAf}\end{align}
    \vspace{-15pt}
     \begin{align}
    &\sum_{t \in \mathcal{T}} \frac{1}{(1+r)^{t-1}}(\Psi\sum_{s \in \mathcal{S}} \sum_{b \in \mathcal{B}}(\sum_{k \in \Omega_{t, s,b}^{D}}(c_{t,s,k,b}^{d}d_{t, s,k,b}\notag\\
    &-  \varphi^{D,\max}_{t,s,k,b}d^{\max}_{t, s,k,b}+\varphi^{D,\min}_{t,s,k,b}d^{\min}_{t, s,k,b}) \notag\\
    &- \sum_{k \in \Omega_{t, s,b}^{G}} (c_{t,s,k,b}^{g}g_{t, s,k,b}+  \varphi^{G,\max}_{t,s,k,b}g^{\max}_{t, s,k,b}-\varphi^{G,\min}_{t,s,k,b}g^{\min}_{t, s,k,b}) \notag\\
    &- \sum_{k \in \Omega_{t, s,b}^{W}} (\varphi_{t,s,k,b}^{W,sch}p^{w,fore}_{t,s,k,b}  +\varphi^{W,cur,\max}_{t,s,k,b}p^{w,fore}_{t,s,k,b}\notag\\
    &+c_{t,s,k,b}^{cur}p^{w,cur}_{t,s,k,b} ))+ CC_t -C_t)\geq 0,\label{eq:MIQCPSLAg}
\end{align}
\vspace{-15pt} 
\begin{align}
    Eqs.~\eqref{eq:upperf},~\eqref{eq:upperi}-\eqref{eq:upperk},~\eqref{eq:upperm}-\eqref{eq:uppert} ,\label{eq:MIQCPSLAh}
\end{align}\vspace{-15pt} 
\begin{align}
Eqs.~\eqref{KKT:SLA} \text{ with linearized complementarity condition}, \label{eq:MIQCPSLAi}
\end{align}\vspace{-15pt} 
\begin{align}
    Eqs.~\eqref{linearization:o}-\eqref{linearization:tau}.\label{eq:MIQCPSLAj}
\end{align}
    \label{eq:MIQCPSLA}
\end{subequations}
\end{fleqn}
The variables of the single-level problem under SLA include 
$$\Xi^{SL}_{SLA} = \{b_{t,l, j}^{R},z^R_{t,l},o_{t,c},z^P_{t,l,m},\tau^V_l,\tau^C,\hat{c}_{t,s, k,b}^{d},\hat{c}_{t,s, k,b}^{g},$$
$$\hat{c}_{t,s, k,b}^{w},d_{t,s,k,b} ,g_{t,s,k,b},p^{w,sch}_{t,s,k,b},p^{w,cur}_{t,s,k,b},u_{t,s},v_{t,s,i},$$
$$\pi_{t,s},\varphi^{D,\min}_{t,s,k,b},\varphi^{D,\max}_{t,s,k,b},\varphi^{G,\min}_{t,s,k,b},\varphi^{G,\max}_{t,s,k,b},\varphi^{W,sch}_{t,s,k,b},\varphi^{W,cur,\min}_{t,s,k,b},$$
$$\varphi^{W,cur,\max}_{t,s,k,b},\mu^{u,\min}_{t,s},\mu^{v,\min}_{t,s,i},\mu^{1}_{t,s},\mu_{t,s,l,i}^{2},\mu_{t,s,l,i}^{3},\mu_{t,s,l}^{4},\mu_{t,s,l}^{5},$$
$$u^{\varphi^{D,\min}}_{t,s,k,b},u^{\varphi^{D,\max}}_{t,s,k,b},u^{\varphi^{G,\min}}_{t,s,k,b},u^{\varphi^{G,\max}}_{t,s,k,b},u^{\varphi^{W,cur,\min}}_{t,s,k,b},u^{\varphi^{W,cur,\max}}_{t,s,k,b},$$
$$u^{\mu^{u,\min}}_{t,s},u^{\mu^{v,\min}}_{t,s,i},u^{\mu^{1}}_{t,s},u^{\mu^{2}}_{t,s,l,i},u^{\mu^{3}}_{t,s,l,i},u^{\mu^{4}}_{t,s,l},u^{\mu^{5}}_{t,s,l},\eta^{R,V}_{t,l},\eta^{P,V}_{t,l,m},$$
$$\eta^{o,2}_{t,s,l,c,i},\eta^{o,3}_{t,s,l,c,i},\eta^{o,4}_{t,s,l,c},\eta^{o,5}_{t,s,l,c}, \eta^{o,g}_{t,s,k,b,c}, \eta^{o,w}_{t,s,k,b,c}, \eta^{o,d}_{t,s,k,b,c}\}$$.

\section{Benchmark models}\label{section:benchmark}
This section presents the reformulation of the bilevel problem under two different linear approximation schemes (i.e., LA and W-CVaR) for the RHS-WDRJCC~\eqref{eq:JCC:bilevel}.

\subsection{LA}
The LA scheme~\eqref{eq:LA} with respect to the RHS-WDRJCC~\eqref{eq:JCC:bilevel} is 
\begin{align}
&Eqs.~\eqref{eq:SLAa}-\eqref{eq:SLAd}\label{eq:bilevel:LA}
\end{align}
We denote the term as 
\begin{equation}
\begin{aligned}
& \Lambda^{LA}_{t,s,b} := \sum_{i\in[N]} \kappa_i \sum_{c\in\mathcal{C}}\sum_{l\in\mathcal{L}}s_{b,l,c}o_{t,c}(\mu^2_{t,s,l,i}- \mu^3_{t,s,l,i})
\end{aligned}\label{eq:Xi_LA_original}
\end{equation}
then KKT conditions under LA are
\begin{subequations}
\begin{align}
& \hat{c}_{t,s,k,b}^{g} -\pi_{t,s} + \varphi^{G,\max}_{t,s,k,b}-\varphi^{G,\min}_{t,s,k,b}+ \Lambda^{LA}_{t,s,b} =0,\notag\\
&\forall t \in \mathcal{T},  s\in \mathcal{S},  k \in \Omega_{t, s,b}^{G}, b \in \mathcal{B}, \label{KKT:LAa}\\
& - \hat{c}_{t,s,k,b}^{d} + \pi_{t,s} +  \varphi^{D,\max}_{t,s,k,b}-\varphi^{D,\min}_{t,s,k,b}-\Lambda^{LA}_{t,s,b} =0,\notag\\
&\forall t \in \mathcal{T}, s\in \mathcal{S}, k \in \Omega_{t, s,b}^{D}, b \in \mathcal{B}, \label{KKT:LAb}\\
&  \hat{c}_{t,s,k,b}^{w}-\pi_{t,s} + \varphi_{t,s,k,b}^{W,sch}+ \Lambda^{LA}_{t,s,b} =0,\notag\\
&\forall t \in \mathcal{T}, s\in \mathcal{S},  k \in \Omega_{t, s,b}^{W}, b \in \mathcal{B},\label{KKT:LAc}\\
& -\mu^{u,\min}_{t,s}-\epsilon N \mu^1_{t,s} + \sum_{l\in\mathcal{L}}\sum_{i\in [N]}(\mu^2_{t,s,l,i}+\mu^3_{t,s,l,i}) = 0,\notag\\
&\forall t \in \mathcal{T}, s\in \mathcal{S}, \label{KKT:LAd}\\
&Eqs.~\eqref{KKT:SLAa},~\eqref{KKT:SLAe}-\eqref{KKT:SLAf},~\eqref{KKT:SLAh}-\eqref{KKT:SLAr} \label{KKT:LAf}
\end{align}
\label{KKT:LA}
\end{subequations}
The LMP $\pi_{t,s,b}$ under LA is 
\begin{flalign}
    & \pi_{t,s,b} = \pi_{t,s}- \Lambda^{LA}_{t,s,b}.\label{LA:LMP}
\end{flalign}
Consequently, the merchandising surplus under the LA is identical to that under the SLA. This equality arises because the terms $\Lambda^{LA}_{t,s,b}$ cancel each other in the calculation of the merchandising surplus, as detailed in Eqs.~\eqref{eq:ms1}--\eqref{eq:ms3}.

To linearize the single-level MPEC problem, ~\eqref{eq:Xi_LA_original} is reformulated as 
\begin{align}
& \Lambda^{SLA}_{t,s,b} := \sum_{i\in[N]} \kappa_i \sum_{c\in\mathcal{C}}\sum_{l\in\mathcal{L}}s_{b,l,c}(\eta^{o,2}_{t,s,l,c,i}- \eta^{o,3}_{t,s,l,c,i})\label{eq:LAterm_updated}
\end{align} using big-Ms with auxiliary variables $\eta^{o,2}_{t,s,l,c,i}$ and $\eta^{o,3}_{t,s,l,c,i}$.


The final mixed-integer problem with bilinear constraints under LA approximation is to replace these following constraints in Problem~\eqref{eq:MIQCPSLA}:
\begin{itemize}
    \item replace~\eqref{eq:MIQCPSLAi} by~\eqref{KKT:LA} with linearized complementarity condition;
    \item replace~\eqref{eq:MIQCPSLAj} by ~\eqref{linearization:oa}--\eqref{linearization:od},~\eqref{linearization:ogwd}--~\eqref{linearization:tau},~\eqref{eq:LAterm_updated}.
\end{itemize}
The variable set for the single-level problem under LA $\Xi^{SL}_{LA}$ is to delete $\mu_{t,s,l}^{4},\mu_{t,s,l}^{5},\eta^{o,4}_{t,s,l,c},\eta^{o,5}_{t,s,l,c},u^{\mu^{4}}_{t,s,l},u^{\mu^{5}}_{t,s,l}$ in $\Xi^{SL}_{SLA}$.

\subsection{W-CVaR}
The W-CVaR scheme~\eqref{eq:CVaR} with respect to the RHS-WDRJCC~\eqref{eq:JCC:bilevel} is:
\begin{subequations}
\begin{align}
& -\alpha_{t,s,i} \leq 0: \mu^{\alpha,\min}_{t,s,i} \geq 0 \quad \forall t \in \mathcal{T},s\in \mathcal{S},i\in[N],\label{eq:WCVARa}\\
& \beta_{t,s} \in \mathbb{R}, \tau_{t,s} \in \mathbb{R}\quad \forall t \in \mathcal{T},s\in \mathcal{S},\label{eq:WCVARa}\\
& \tau_{t,s} + \frac{1}{\epsilon} (\theta \beta_{t,s} + \frac{1}{N} \sum_{i\in [N]} \alpha_{t,s,i}) \leq 0:\mu^{1}_{t,s} \notag\\
&\forall t \in \mathcal{T},s\in \mathcal{S},\label{eq:WCVARb} \\
& -\tau_{t,s}-\alpha_{t,s,i} + w^{\max}_{t,s,l}\Lambda_{t,s,l}\leq \notag\\
&w^{\max}_{t,s,l} (-\sum_{c\in\mathcal{C}}\sum_{b\in\mathcal{B}}\sum_{k\in\Omega_{t, s,b}^{W}} s_{b,l,c}o_{t,c}  \tilde{e}_{t,s,k,b} \notag\\
&+\mathcal{F}^{0}_{l}(1+\sum_{\hat{t} \leq t}\sum_{j\in\mathcal{J}}j b^R_{\hat{t},l,j})+\sum_{m \in \mathcal{M}} m z^P_{t,l,m} \overline{F}^{C}_{l}):\mu_{t,s,l,i}^{2}\notag\\
& \forall t \in \mathcal{T},s\in \mathcal{S},l \in \mathcal{L},i\in[N], \label{eq:WCVARc}\\
& -\tau_{t,s}-\alpha_{t,s,i} - w^{\min}_{t,s,l}\Lambda_{t,s,l}\leq \notag\\
&w^{\min}_{t,s,l} (\sum_{c\in\mathcal{C}}\sum_{b\in\mathcal{B}}\sum_{k\in\Omega_{t, s,b}^{W}} s_{b,l,c}o_{t,c}  \tilde{e}_{t,s,k,b}\notag\\
&+\mathcal{F}^{0}_{l}(1+\sum_{\hat{t} \leq t}\sum_{j\in\mathcal{J}}j b^R_{\hat{t},l,j})+\sum_{m \in \mathcal{M}} m z^P_{t,l,m} \overline{F}^{C}_{l}):\mu_{t,s,l,i}^{3}\notag\\
&\forall t \in \mathcal{T},s\in \mathcal{S},l \in \mathcal{L},i\in[N], \label{eq:WCVARd}\\
&  w^{\max}_{t,s,l} - \beta_{t,s} \leq 0 :\mu_{t,s,l}^{4} \quad  \forall t \in \mathcal{T},s\in \mathcal{S},l\in \mathcal{L} ,\label{eq:WCVARe}\\
&  w^{\min}_{t,s,l} - \beta_{t,s} \leq 0 :\mu_{t,s,l}^{5}  \quad  \forall t\in \mathcal{T},s\in \mathcal{S},l\in \mathcal{L} .\label{eq:WCVARf}
\end{align}
\end{subequations}

We denote the term as 
\begin{equation}
\begin{aligned}
&\Lambda^{WCVaR}_{t,s,b} := \\ &\sum_{i\in[N]}\sum_{c\in\mathcal{C}}\sum_{l\in\mathcal{L}}s_{b,l,c}o_{t,c}(w^{\max}_{t,s,l} \mu^2_{t,s,l,i} - w^{\min}_{t,s,l}\mu^3_{t,s,l,i} )
\end{aligned}
\end{equation}
then KKT conditions under W-CVaR are
\begin{subequations}
\begin{align}
& Eqs.~\eqref{eq:lowerb},\eqref{eq:lowere}, \label{KKT:WCVARa}\\
& \hat{c}_{t,s,k,b}^{g} -\pi_{t,s} + \varphi^{G,\max}_{t,s,k,b}-\varphi^{G,\min}_{t,s,k,b}+ \Lambda^{WCVaR}_{t,s,b}=0,\notag\\
&\forall t \in \mathcal{T},  s\in \mathcal{S},  k \in \Omega_{t, s,b}^{G}, b \in \mathcal{B}, \label{KKT:WCVARb}\\
& - \hat{c}_{t,s,k,b}^{d} + \pi_{t,s} +  \varphi^{D,\max}_{t,s,k,b}-\varphi^{D,\min}_{t,s,k,b}- \Lambda^{WCVaR}_{t,s,b}=0,\notag\\
&\forall t \in \mathcal{T},  s\in \mathcal{S},  k \in \Omega_{t, s,b}^{D}, b \in \mathcal{B}, \label{KKT:WCVARc}\\
&   \hat{c}_{t,s,k,b}^{w}-\pi_{t,s} + \varphi_{t,s,k,b}^{W,sch}+ \Lambda^{WCVaR}_{t,s,b}=0,\notag\\
&\forall t \in \mathcal{T},  s\in \mathcal{S},  k \in \Omega_{t, s,b}^{W}, b \in \mathcal{B},\label{KKT:WCVARd}\\
& c_{t,s,k,b}^{cur} + \varphi_{t,s,k,b}^{W,sch} +\varphi^{W,cur,\max}_{t,s,k,b}-\varphi^{W,cur,\min}_{t,s,k,b} =0,\notag\\
&\forall t \in \mathcal{T},  s\in \mathcal{S},  k \in \Omega_{t, s,b}^{W}, b \in \mathcal{B}, \label{KKT:WCVARe}\\
& -\mu^{\alpha,\min}_{i,t,s}+ \frac{1}{\epsilon N}\mu^1_{t,s} -\sum_{l\in\mathcal{L}}(\mu^2_{t,s,l,i}+\mu^3_{t,s,l,i}) = 0,\notag\\
&\forall t \in \mathcal{T},s \in \mathcal{S}, i\in[N], \label{KKT:WCVARf}\\
& \frac{\theta}{\epsilon} \mu^1_{t,s} -\sum_{l\in\mathcal{L}}(\mu^4_{t,s,l}+\mu^5_{t,s,l}) = 0,\quad\forall t \in \mathcal{T},s \in\mathcal{S},\label{KKT:WCVARg}\\
& \mu^1_{t,s} -\sum_{i\in [N]}\sum_{l\in\mathcal{L}}(\mu^2_{t,s,l,i}+\mu^3_{t,s,l,i})= 0,\notag\\
&\forall t \in \mathcal{T}, s\in\mathcal{S},  \label{KKT:WCVARh}\\
& Eqs.~\eqref{KKT:SLAh}-\eqref{KKT:SLAm},\label{KKT:WCVARn}\\
&0 \leq \alpha_{t,s,i} \perp \mu^{\alpha,\min}_{t,s,i} \geq 0 ,\quad \forall t \in \mathcal{T},s\in \mathcal{S},i\in[N],\label{KKT:WCVARo}\\
&0 \leq -(\tau_{t,s} + \frac{1}{\epsilon} (\theta \beta_{t,s} + \frac{1}{N} \sum_{i\in [N]} \alpha_{t,s,i})) \perp \mu^{1}_{t,s} \geq 0 ,\notag\\
& \forall t\in \mathcal{T},s \in \mathcal{S},  \label{KKT:WCVARp}\\
&0 \leq  \big(w^{\max}_{t,s,l} (-\sum_{c\in\mathcal{C}}\sum_{b\in\mathcal{B}}\sum_{k\in\Omega_{t, s,b}^{W}} s_{b,l,c}o_{t,c}  \tilde{e}_{t,s,k,b} \notag\\
&+\mathcal{F}^{0}_{l}(1+\sum_{\hat{t} \leq t}\sum_{j\in\mathcal{J}}j b^R_{\hat{t},l,j})+\sum_{m \in \mathcal{M}} m z^P_{t,l,m} \overline{F}^{C}_{l})\notag\\
&+\tau_{t,s}+\alpha_{t,s,i} - w^{\max}_{t,s,l}\Lambda_{t,s,l}\big) \perp \mu_{t,s,l,i}^{2} \geq 0,\notag\\
& \forall t\in \mathcal{T},s \in \mathcal{S},l\in\mathcal{L},i\in[N], \label{KKT:WCVARq}\\
&0 \leq  \big(w^{\min}_{t,s,l} (\sum_{c\in\mathcal{C}}\sum_{b\in\mathcal{B}}\sum_{k\in\Omega_{t, s,b}^{W}} s_{b,l,c}o_{t,c}  \tilde{e}_{t,s,k,b} \notag\\
&+\mathcal{F}^{0}_{l}(1+\sum_{\hat{t} \leq t}\sum_{j\in\mathcal{J}}j b^R_{\hat{t},l,j})+\sum_{m \in \mathcal{M}} m z^P_{t,l,m} \overline{F}^{C}_{l})\notag\\
&+\tau_{t,s}+\alpha_{t,s,i}+ w^{\min}_{t,s,l}\Lambda_{t,s,l}\big) \perp \mu_{t,s,l,i}^{3} \geq 0,\notag\\
&\forall t\in \mathcal{T},s \in \mathcal{S},l\in\mathcal{L},i\in[N],\label{KKT:WCVARr}\\
& 0 \leq  -w^{\max}_{t,s,l} + \beta_{t,s}  \perp \mu_{t,s,l}^{4}\geq 0,\forall t\in \mathcal{T},s \in \mathcal{S},l \in\mathcal{L},\label{KKT:WCVARs}\\
& 0 \leq  -  w^{\min}_{t,s,l} + \beta_{t,s}  \perp\mu_{t,s,l}^{5}\geq 0,\forall t \in\mathcal{T}, s\in\mathcal{S},l \in\mathcal{L}. \label{KKT:WCVARt}
\end{align}
\label{KKT:WCVAR}
\end{subequations}
The LMP $\pi_{t,s,b}$ under W-CVaR is 
\begin{flalign}
    & \pi_{t,s,b} = \pi_{t,s}- \Lambda^{WCVaR}_{t,s,b}\label{WCVAR:LMP}
\end{flalign}
Consequently, the merchandising surplus under the W-CVaR is identical to that under the SLA. This equality arises because the terms $\Lambda^{WCVaR}_{t,s,b}$ cancel each other in the calculation of the merchandising surplus, as detailed in Eqs.~\eqref{eq:ms1}--\eqref{eq:ms3}.

We will also introduce auxiliary variables $\eta^{o,2}_{t,s,l,c,i},\eta^{o,3}_{t,s,l,c,i}$ for $o_{t,c} \mu^2_{t,s,l,i}$ and $o_{t,c} \mu^3_{t,s,l,i}$ with constraints~\eqref{linearization:oa}--\eqref{linearization:od}. Then $\Lambda^{WCVaR}_{t,s,b}$ becomes 
\begin{equation}
\begin{aligned}
&\Lambda^{WCVaR}_{t,s,b}\\
&:= \sum_{i\in[N]}  \sum_{c\in\mathcal{C}}\sum_{l\in\mathcal{L}}s_{b,l,c}(w^{\max}_{t,s,l}\eta^{o,2}_{t,s,l,c,i} -w^{\min}_{t,s,l} \eta^{o,3}_{t,s,l,c,i} ) 
\end{aligned}  \label{eq:WCVARterm} 
\end{equation}

The final mixed-integer problem with bilinear constraints under W-CVaR approximation is to replace these following constraints in Problem~\eqref{eq:MIQCPSLA}:
\begin{itemize}
    \item replace~\eqref{eq:MIQCPSLAi} by~\eqref{KKT:WCVAR} with linearized complementarity condition 
    \item replace~\eqref{eq:MIQCPSLAj} by ~\eqref{linearization:oa}--\eqref{linearization:od},~\eqref{linearization:ogwd}--~\eqref{linearization:tau},~\eqref{eq:WCVARterm} .
\end{itemize}
The variable set for the single-level problem under W-CVaR $\Xi^{SL}_{WCVaR}$ is to replace $u_{t,s},v_{t,s,i},\mu^{u,\min}_{t,s},\mu^{v,\min}_{t,s,i},u^{\mu^{u,\min}}_{t,s},u^{\mu^{v,\min}}_{t,s,i},\eta^{o,4}_{t,s,l,c},\eta^{o,5}_{t,s,l,c}$ in $\Xi^{SL}_{SLA}$ with $\alpha_{t,s,i}, \beta_{t,s}, \tau_{t,s}, \mu^{\alpha,\min}_{t,s,i},u^{\mu^{\alpha,\min}}_{t,s,i}$.

\section{Case Study Setting}\label{section:case_study setting}
\begin{figure}[htbp] \label{fig:6node_network}
  \centering
  \includegraphics[width=\linewidth]{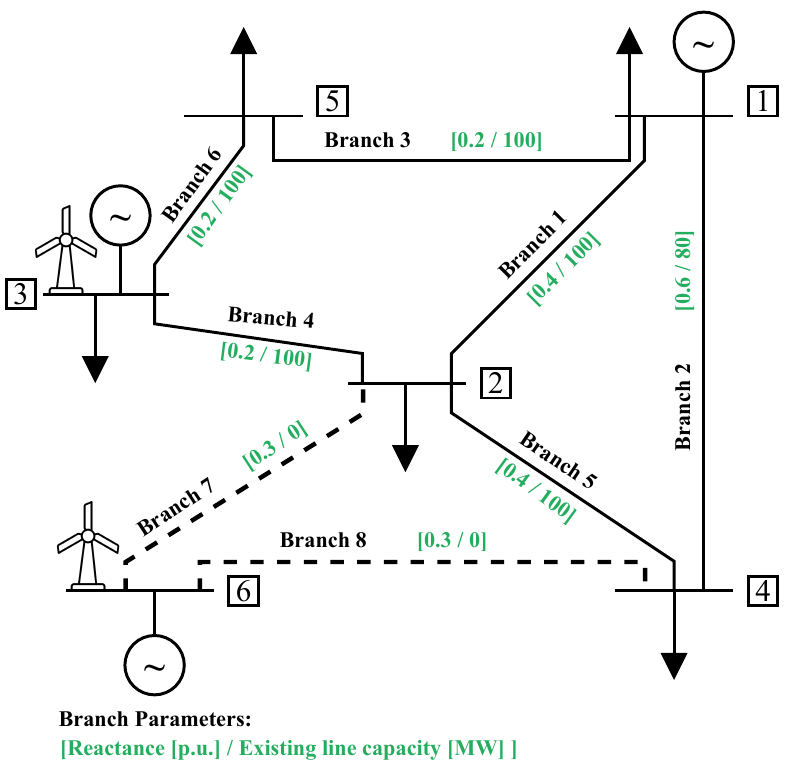}
  \caption{Topology of the modified Garver’s 6-node transmission network.}
\end{figure}

\begin{figure}[htbp] \label{fig:histogram}
  \centering
  \includegraphics[width=\linewidth]{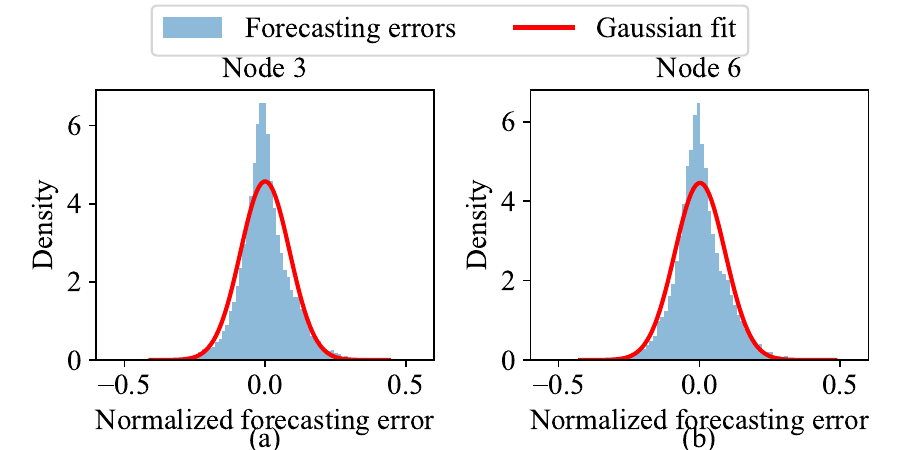}
  \caption{Histograms of normalized wind power forecasting errors at Node 3 and Node 6, overlaid with fitted Gaussian distributions.}
\end{figure}
 
This case study evaluates the proposed bilevel investment planning framework using a modified Garver 6-node transmission network, as shown in Fig.~\ref{fig:6node_network}. The potential addition of new transmission lines is considered for the corridors connecting nodes $(2,6)$ and $(4,6)$. Each candidate line has a reactance of $0.3$ p.u., a transmission capacity of $100 MW$, and an investment cost of \pounds30M, and we assume any corridor can accommodate at most three lines \cite{4435945}.  Reconductoring is considered for existing lines in corridors $(2,3)$ and $(3,5)$, with a fixed investment cost $K^{fix}_l$ of \pounds1M and a variable cost $K^{var}_l$ of \pounds0.1M/MW of added capacity. The investment cost for constructing a new parallel transmission line is set approximately three times higher than that of reconductoring an existing line~\cite{2411207121}. Discrete capacity increase factors for reconductoring are considered, denoted by the set $J=[0\%, 5\%, 10\%, \cdots, 200\%]$. 

Conventional generators are located at nodes $\{1, 3, 6\}$, while consumers are situated at nodes $\{1, 2, 3, 4, 5\}$. At each node where either generation or consumption is present, it is assumed there are five individual market participants. Demand bids at each node are sourced from Table III in~\cite{4435945}. Given the presence of more generators compared to~\cite{4435945}, generator bids are sampled from a normal distribution with mean values taken from the average generator bids reported in Table III of~\cite{4435945}, and standard deviations set to $10\%$ of the corresponding means to introduce variability. 

The referenced total generation capacities at nodes 1, 3, and 6 are $150$ MW, $360$ MW, and $600$ MW, respectively in~\cite{4435945}. To introduce variability in the generation capacity at these nodes, the capacity of each generator participant is sampled from a uniform distribution between $0.5$ and $1.5$ times the average referenced total capacity at the respective node. The referenced total demand capacities at nodes 1, 2, 3, 4, and 5 are $80$ MW, $240$ MW, $40$ MW, $160$ MW, and $240$ MW, respectively~\cite{4435945}. Since we have flexible demand and wind generation, we introduce variability in the demand by sampling the base demand capacity for each consumer participant from a uniform distribution ranging between 1 and 2 times the average referenced total demand at the respective node.

Two wind farms exist within the system, located at node 3 and node 6, with a combined installed wind capacity of $380$ MW. This capacity corresponds to a $50\%$ penetration level relative to the total referenced system load. The wind farm at node 6 has a larger installed capacity ($80\%$ of the total wind capacity), while the wind farm at node 3 accounts for the remaining $20\%$. This configuration is designed to create a scenario where a substantial portion of generation is initially located at the initially isolated node 6, likely incentivizing transmission expansion plans that construct lines connecting node 6 to the rest of the system. The wind curtailment cost is set at \pounds60/MWh.

Wind power data is handled following the methodology outlined in our previous work~\cite{zhou2024strengthenedfasterlinearapproximation}, utilizing publicly available Bronze-medal wind forecasts from the Global Energy Forecasting Competition 2012 – Wind Forecasting~\cite{gefc2012}, corresponding to the second wind farm. Forecasting errors are estimated by computing the difference between these predictions and the actual ground-truth values from~\cite{gefc2012}, yielding 156 samples of next 48-hour wind power forecasting errors. We adopt the approach in~\cite{zhou2024strengthenedfasterlinearapproximation} to construct a temporal correlation matrix and model the marginal distribution of forecasting errors. To increase the sample size while preserving both marginal distributions and temporal correlations, we employ the Gaussian kernel density estimate (KDE) for marginals and a Gaussian copula to generate new temporally correlated samples. Fig.~\ref{fig:histogram} illustrates the empirical distributions of normalized wind power forecasting errors at Node 3 and Node 6, along with fitted Gaussian curves.

The transmission investment planning problem is examined under a single operational scenario ($\mathcal{S} = 1$) for each year, with each investment period comprising 8,760 hourly operational intervals. We consider a planning horizon of four investment periods ($\mathcal{T} = 4$) and a discount rate of $5\%$. The total system demand is assumed to increase by $5\%$ annually. The network cost allocation factor $\delta_{l,b}$ is set to 1 for all lines and buses as we assume a postage stamp network cost allocation approach~\cite{SAVELLI2020113979}. For all norm calculations in the RHS-WDRJCC formulations, we adopt the L2 norm which has the property of being self-dual. Optimal hyper-parameters are set to $\bm{\kappa}=1$ for the proposed SLA and LA, and $\bm{w}_p=1/P$ for the W-CVaR (see Cor.~\ref{cor:wcvar}). We have 50 training samples ($N=50$) for the RHS-WDRJCC problems.

\section{Statistical Results}\label{section:statistical_results}

\begin{table*}[h]
\centering
\caption{Out-of-sample reliability results across investment years (Y1-Y4). Red values indicate instances where achieved reliability falls below the desired reliability.}
\begin{tabular}{l|c|c|c}
\toprule
\hline
$1-\epsilon$ (\%)  & $\theta=0.1$ & $\theta=0.2$ & $\theta=0.3$ \\
\hline
$70$ & $84.33$, $87.40$, $90.15$, $88.42$  & $85.05$, $87.83$, $90.30$, $88.72$ &  $85.42$, $88.45$, $90.75$, $88.98$ \\
$80$  & $86.92$, $90.85$, $92.90$, $92.10$ & $87.40$, $91.53$, $93.35$, $92.47$ & $88.17$, $91.75$, $93.62$, $92.88$ \\
$90$  & $91.17$, $94.73$, $96.47$, $97.08$ & $92.05$, $95.12$, $97.02$, $97.47$ & $92.80$, $95.58$, $97.45$, $97.65$ \\
$95$ & $\textcolor{red}{\textbf{94.47}}$, $97.10$, $98.60$, $97.45$ &
$95.62$, $97.62$, $99.00$, $98.20$ &
$96.78$, $98.28$, $99.22$, $98.78$ \\
$97.5$ &  $97.58$, $98.55$, $99.55$, $98.17$ &
$98.30$, $99.20$, $99.80$, $99.15$ &
$98.80$, $99.45$, $99.92$, $99.62$ \\
$99$ & $\textcolor{red}{\textbf{98.17}}$, $99.15$, $99.88$, $99.38$ &
$99.48$, $99.80$, $99.98$, $99.95$ &
$99.88$, $99.95$, $99.98$, $99.98$ \\
\hline\bottomrule
\end{tabular}\label{tab:relia}
\end{table*}

\begin{table*}[h]
\centering
\caption{Investment results showing corridor expansions (2,3), (3,5), (2,6), (4,6) in MW, followed by total MW. Expansions without brackets occurred in Year 1; those with brackets indicate the specific expansion year.}
\begin{tabular}{l|c|c|c}
\toprule
\hline
$1-\epsilon$ (\%)  & $\theta=0.1$ & $\theta=0.2$ & $\theta=0.3$ \\
\hline
$70$ & $0$, $65$, $100$, $300$:\hspace{0.5em}\hspace{0.5em}\hspace{0.5em}\hspace{0.7em} $465$ & $0$, $65$, $100$, $300$:\hspace{0.5em}\hspace{0.5em}\hspace{0.5em}\hspace{0.7em} $465$ &  $0$, $65$, $100$, $300$:\hspace{0.5em}\hspace{0.5em}\hspace{0.5em}\hspace{0.7em} $465$ \\
$80$ & $0$, $75$, $200$, $300$:\hspace{0.5em}\hspace{0.5em}\hspace{0.5em}\hspace{0.7em} $575$ & $0$, $75$, $200$, $300$:\hspace{0.5em}\hspace{0.5em}\hspace{0.5em}\hspace{0.7em} $575$ & $0$, $80$, $200$, $300$:\hspace{0.5em}\hspace{0.5em}\hspace{0.5em}\hspace{0.7em} $580$  \\
$90$ & $0$, $85$, $200$, $300$:\hspace{0.5em}\hspace{0.5em}\hspace{0.5em}\hspace{0.7em} $585$  & $0$, $85$, $200$, $300$:\hspace{0.5em}\hspace{0.5em}\hspace{0.5em}\hspace{0.7em} $585$ &  $0$, $85$, $200$, $300$:\hspace{0.5em}\hspace{0.5em}\hspace{0.5em}\hspace{0.7em} $585$ \\
$95$ & $0$, $80$, $200$, $300$:\hspace{0.5em}\hspace{0.5em}\hspace{0.5em}\hspace{0.7em} $580$ & $0$, $85$, $200$, $300$:\hspace{0.5em}\hspace{0.5em}\hspace{0.5em}\hspace{0.7em} $585$  & $40(2)$, $100$, $200$, $300$: $640$ \\
$97.5$ & $45(2)$, $110$, $200$, $300$: $655$ & $50(2)$, $110$, $200$, $300$: $660$  &  $55$, $110$, $200$, $300$:\hspace{0.5em}\hspace{0.7em} $665$ \\
$99$ & $50(2)$, $115$, $200$, $300$: $665$ & $65$, $120$, $200$, $300$:\hspace{0.5em}\hspace{0.7em} $685$  &  $80$, $130$, $200$, $300$:\hspace{0.5em}\hspace{0.7em} $710$ \\
\hline \bottomrule
\end{tabular} \label{tab:investment}
\end{table*}

\section{Evaluation metrics}\label{appendix:evaluation_metrics}
The exclusive evaluation metrics used in this case study include:
\begin{enumerate}
    \item \emph{TimeF} (s): 
    This metric records the time (in seconds) to obtain the first ``comparable high-quality'' solution, defined as the solution with an objective within the \texttt{MIPGap} of the final objective achieved by the proposed SLA (whether optimal or the best found within the \texttt{TimeLimit}). If SLA does not yield a feasible solution within the \texttt{TimeLimit}, the \emph{TimeF} of the proposed SLA is set to the \texttt{TimeLimit}, and the \emph{TimeF} for other benchmarks is set to the time when the first feasible solution is found (if any, otherwise set to the \texttt{TimeLimit}). This metric offers a practical perspective in cases where the solver identifies a solution close to optimal early but requires extra time for verifying optimality.
    \item \emph{Time} (s): This measures the total time (in seconds) spent to solve the problem to optimality (within \texttt{MIPGap}), and is set to the \texttt{TimeLimit} if no optimal solution is found within the \texttt{TimeLimit}. This metric complements the previous \emph{TimeF} by reporting the total time to reach optimality.
    \item \emph{Nsolvable} (s): The number of instances (out of 30) for which a feasible solution was found within the \texttt{TimeLimit}. This metric reflects the practical solvability of each method under \texttt{TimeLimit}.
    \item \emph{Bilevel Obj.} (M\pounds): Objective value obtained by solving the bilevel problem within the \texttt{TimeLimit}.
    \item \emph{Obj. Diff.} (\%): This represents the difference between the objective values achieved by the benchmark methods and SLA, considering only cases where both methods returned feasible solutions within the \texttt{TimeLimit}. A negative value indicates that the benchmark methods achieve lower social welfare (i.e., worse optimality) than SLA, while a positive value indicates that SLA performs worse than the benchmark methods. We only keep solvable runs where the benchmark method and the proposed SLA identified feasible solutions within the \texttt{TimeLimit}. 
    \item \emph{Reli.} (\%): This is the out-of-sample joint satisfaction rate of the RHS-WDRJCC~\eqref{eq:ll:jcc1}--\eqref{eq:ll:jcc2}. This represents the averaged value of reliability of $t\in\mathcal{T}, s\in\mathcal{S}$, which is calculated on $4,000$ testing samples.
\end{enumerate}
In addition, we classify the runs based on results into two groups:
\begin{enumerate}
    \item \emph{All Runs (AR)}: All 30 random runs.
    \item \emph{Comparable Runs (CR)}: Only those runs in which \textit{all} three methods obtained a feasible solution within the \texttt{TimeLimit}.
\end{enumerate}
\begin{figure*}
    \centering    \includegraphics[width=\linewidth]{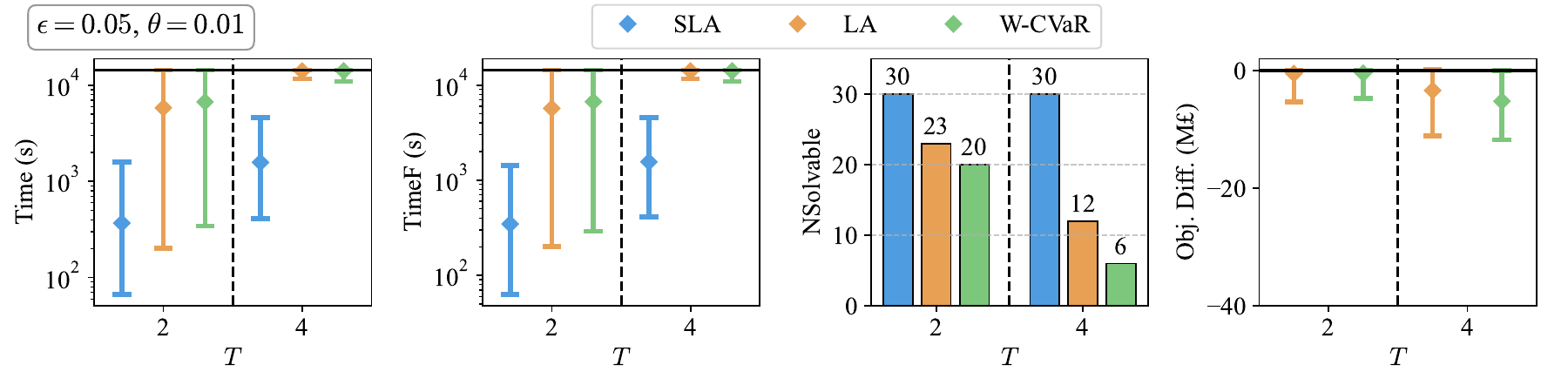}
    \begin{minipage}{\textwidth}
    \hspace{1em}
    \begin{minipage}{\textwidth}
        \centering
        \vspace{-3em}  
        \subfloat[]{\phantom{\rule{0.249\textwidth}{0pt}}}
        \subfloat[]{\phantom{\rule{0.249\textwidth}{0pt}}}
        \subfloat[\label{fig:n_infeasible_eps0.05_theta0.01}]{\phantom{\rule{0.249\textwidth}{0pt}}}
        \subfloat[\label{fig:profit_diff_eps0.05_theta0.01}]{\phantom{\rule{0.249\textwidth}{0pt}}}
    \end{minipage}
    \end{minipage}
    \includegraphics[width=\linewidth]{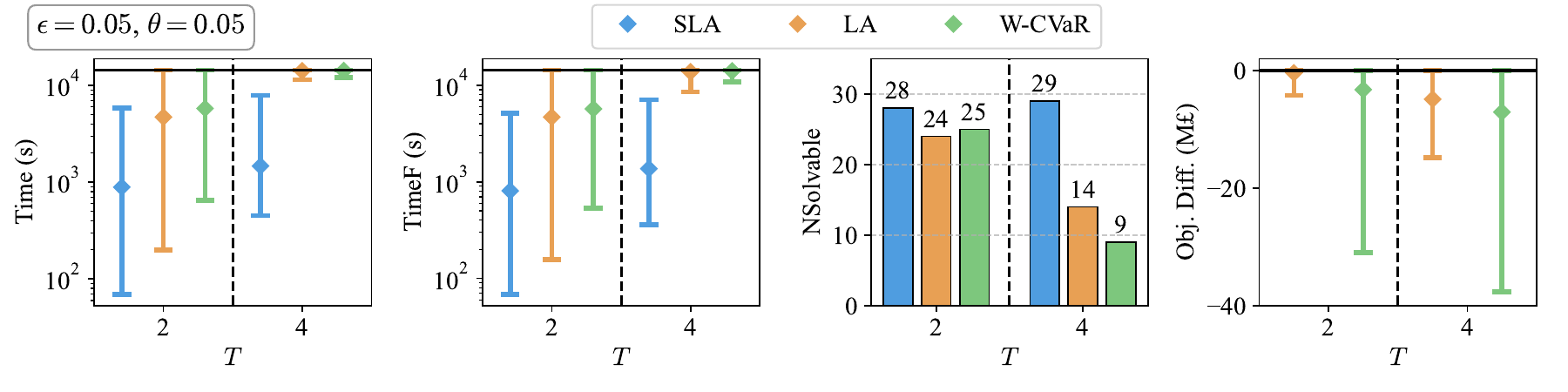}
    \begin{minipage}{\textwidth}
    \hspace{1em}
    \begin{minipage}{\textwidth}
        \centering
        \vspace{-3em}  
        \subfloat[]{\phantom{\rule{0.249\textwidth}{0pt}}}
        \subfloat[]{\phantom{\rule{0.249\textwidth}{0pt}}}
        \subfloat[\label{fig:n_infeasible_eps0.05_theta0.05}]{\phantom{\rule{0.249\textwidth}{0pt}}}
        \subfloat[\label{fig:profit_diff_eps0.05_theta0.05}]{\phantom{\rule{0.249\textwidth}{0pt}}}
    \end{minipage}
    \end{minipage}
    \includegraphics[width=\linewidth]{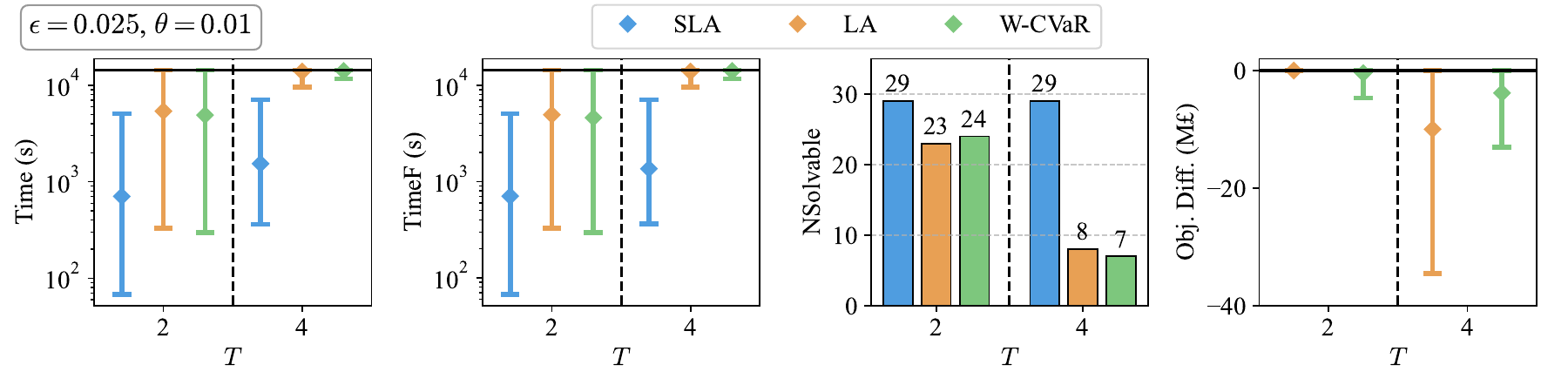}
    \begin{minipage}{\textwidth}
    \hspace{1em}
    \begin{minipage}{\textwidth}
        \centering
        \vspace{-3em}  
        \subfloat[]{\phantom{\rule{0.249\textwidth}{0pt}}}
        \subfloat[]{\phantom{\rule{0.249\textwidth}{0pt}}}
        \subfloat[\label{fig:n_infeasible_eps0.025_theta0.01}]{\phantom{\rule{0.249\textwidth}{0pt}}}
        \subfloat[\label{fig:profit_diff_eps0.025_theta0.01}]{\phantom{\rule{0.249\textwidth}{0pt}}}
    \end{minipage}
    \end{minipage}
    \includegraphics[width=\linewidth]{fig/bilevel_eps0.025_theta0.05_N50_FGap1e-06_Tcompare.pdf}
    \begin{minipage}{\textwidth}
    \hspace{1em}
    \begin{minipage}{\textwidth}
        \vspace{-3em}  
        \subfloat[]{\phantom{\rule{0.249\textwidth}{0pt}}}
        \subfloat[]{\phantom{\rule{0.249\textwidth}{0pt}}}
        \subfloat[\label{fig:n_infeasible_eps0.025_theta0.05}]{\phantom{\rule{0.249\textwidth}{0pt}}}
        \subfloat[\label{fig:profit_diff_eps0.025_theta0.05}]{\phantom{\rule{0.249\textwidth}{0pt}}}
    \end{minipage}
    \end{minipage}
    \vspace{-3em}
    \caption{ 
    Computational performance comparison for the proposed bilevel problem under different parameter combinations of risk level ($\epsilon \in \{0.05,0.025\}$), radius ($\theta \in \{0.01,0.05\}$), and the number of investment years ($T\in\{2,4\}$). Subplots (e)--(p) replicate (a)--(b) but for varying $(\epsilon, \theta)$ combinations. For \textit{Time }(s) and \textit{TimeF} (s) plots, dots represent the mean values of the $30$ random runs, with error bars representing the $95\%$ percentile interval (from 2.5th to 97.5th) of 30 randomly generated instances. In contrast, error bars in the \textit{Obj. Diff.} M\pounds plot represent $95\%$ percentile interval of runs for which the proposed SLA and benchmark are both solvable. The black horizontal line in the \textit{Time }(s) and \textit{TimeF} (s) plots represents the $14,400$s \texttt{TimeLimit}. The black horizontal line in the \textit{Obj. Diff.} (M\pounds) plots indicates zero difference, while negative values denote a lower objective value achieved by the benchmark method compared to the proposed SLA within the \texttt{TimeLimit}.}
    \label{fig:bilevel_com_results}
\end{figure*}

\begin{table*}[h]
\centering
\caption{Comparison of Average Computation Time}
\begin{tabular}{ccc|ccc|ccc}
\toprule
& & & \multicolumn{3}{c|}{\textit{Time} (s)} & \multicolumn{3}{c}{\textit{TimeF} (s)} \\
\cmidrule{4-9}
$\epsilon$ & $\theta$ & $T$ & SLA & LA & W-CVaR  & SLA   & LA   & W-CVaR  \\
\midrule
\multicolumn{9}{c}{All Runs (AR)} \\
\midrule
$0.025$ & $0.01$ & $2$ & $705.18$ & $5400.62$ \hspace{0.5em}\hspace{0.5em}$(7.66\times)$ & $4906.58$ \hspace{0.5em}\hspace{0.5em}$(6.96\times)$ & $704.37$ & $4936.37$ \hspace{0.5em}\hspace{0.5em}$(7.01\times)$ & $4602.94$ \hspace{0.5em}\hspace{0.5em}$(6.53\times)$ \\
$0.025$ & $0.01$ & $4$ & $1540.79$ & $13892.25$\hspace{0.5em} $(9.02\times)$ & $14078.34$ \hspace{0.5em}$(9.14\times)$ & $1357.72$ & $13888.25$ $(10.23\times)$ & $14077.57$ $(10.37\times)$ \\
$0.025$ & $0.05$ & $2$ & $221.72$ & $5080.64$ \hspace{0.5em}$(22.91\times)$ & $5885.24$ \hspace{0.5em}$(26.54\times)$ & $218.40$ & $5048.82$ \hspace{0.5em}$(23.12\times)$ & $5852.70$ \hspace{0.5em}$(26.80\times)$ \\
$0.025$ & $0.05$ & $4$ & $941.07$ & $13698.96$ $(14.56\times)$ & $14308.04$ $(15.20\times)$ & $859.55$ & $13698.86$ $(15.94\times)$ & $14290.98$ $(16.63\times)$ \\
$0.05$ & $0.01$ & $2$ & $368.66$ & $5824.94$ \hspace{0.5em}$(15.80\times)$ & $6710.46$ \hspace{0.5em}$(18.20\times)$ & $346.57$ & $5695.30$\hspace{0.5em} $(16.43\times)$ & $6698.97$ \hspace{0.5em}$(19.33\times)$ \\
$0.05$ & $0.01$ & $4$ & $1573.64$ & $14063.45$ \hspace{0.5em}$(8.94\times)$ & $14076.75$ \hspace{0.5em}$(8.95\times)$ & $1557.27$ & $14041.17$ \hspace{0.5em}$(9.02\times)$ & $14065.87$ \hspace{0.5em}$(9.03\times)$ \\
$0.05$ & $0.05$ & $2$ & $890.77$ & $4706.49$ \hspace{0.5em}\hspace{0.5em}$(5.28\times)$ & $5783.30$ \hspace{0.5em}\hspace{0.5em}$(6.49\times)$ & $806.90$ & $4689.39$ \hspace{0.5em}\hspace{0.5em}$(5.81\times)$ & $5702.43$\hspace{0.5em}\hspace{0.5em} $(7.07\times)$ \\
$0.05$ & $0.05$ & $4$ & $1472.31$ & $14136.31$ \hspace{0.5em}$(9.60\times)$ & $14223.04$\hspace{0.5em} $(9.66\times)$ & $1367.79$ & $13792.38$ $(10.08\times)$ & $14079.15$ $(10.29\times)$ \\
\midrule
\multicolumn{9}{c}{Comparable Runs (CR)} \\
\midrule
$0.025$ & $0.01$ & $2$ & $180.78$ & $2933.28$ \hspace{0.5em}$(16.23\times)$ & $2025.26$ \hspace{0.5em}$(11.20\times)$ & $179.83$ & $2199.61$\hspace{0.5em} $(12.23\times)$ & $2022.85$ \hspace{0.5em}$(11.25\times)$ \\
$0.025$ & $0.01$ & $4$ & $7549.98$ & $9670.59$ \hspace{0.5em}\hspace{0.5em}$(1.28\times)$ & $10909.98$ \hspace{0.5em}$(1.45\times)$ & $7549.00$ & $9670.00$ \hspace{0.5em}\hspace{0.5em}$(1.28\times)$ & $10909.00$ \hspace{0.5em}$(1.45\times)$ \\
$0.025$ & $0.05$ & $2$ & $202.80$ & $2890.62$\hspace{0.5em}\hspace{0.5em}$(14.25\times)$ & $1870.42$ \hspace{0.5em}\hspace{0.5em}$(9.22\times)$ & $197.76$ & $2844.76$\hspace{0.5em} $(14.38\times)$ & $1814.35$ \hspace{0.5em}\hspace{0.5em}$(9.17\times)$ \\
$0.025$ & $0.05$ & $4$ & $937.78$ & $12896.38$ $(13.75\times)$ & $14400.67$ $(15.36\times)$ & $937.00$ & $12896.38$ $(13.76\times)$ & $14400.00$ $(15.37\times)$ \\
$0.05$ & $0.01$ & $2$ & $336.08$ & $2530.71$ \hspace{0.5em}\hspace{0.5em}$(7.53\times)$ & $2913.00$ \hspace{0.5em}\hspace{0.5em}$(8.67\times)$ & $332.69$ & $2528.44$ \hspace{0.5em}\hspace{0.5em}$(7.60\times)$ & $2892.56$\hspace{0.5em}\hspace{0.5em} $(8.69\times)$ \\
$0.05$ & $0.01$ & $4$ & $1164.58$ & $13477.25$ $(11.57\times)$ & $12668.60$ $(10.88\times)$ & $1164.00$ & $13318.75$ $(11.44\times)$ & $12587.25$ $(10.81\times)$ \\
$0.05$ & $0.05$ & $2$ & $307.39$ & $2070.68$\hspace{0.5em}\hspace{0.5em} $(6.74\times)$ & $4000.46$\hspace{0.5em} $(13.01\times)$ & $305.68$ & $2049.98$ \hspace{0.5em}\hspace{0.5em}$(6.71\times)$ & $3970.94$ \hspace{0.5em}$(12.99\times)$ \\
$0.05$ & $0.05$ & $4$ & $2880.33$ & $13564.70$\hspace{0.5em} $(4.71\times)$ & $14361.08$ \hspace{0.5em}$(4.99\times)$ & $2797.60$ & $13546.60$ \hspace{0.5em}$(4.84\times)$ & $14360.69$ \hspace{0.5em}$(5.13\times)$ \\
\bottomrule
\end{tabular}
\label{tab:computation_time}
\end{table*}

\begin{table*}[h]
\centering
\caption{Comparison of Average Objective Values and Out-of-Sample Reliability}
\begin{tabular}{ccc|ccc|ccc}
\toprule \hline
& & & \multicolumn{3}{c}{\textit{Bilevel Obj.} (M£)} & \multicolumn{3}{c}{\textit{Reli.} (\%)} \\
\cmidrule(lr){4-6}\cmidrule(lr){7-9}
$\epsilon$ & $\theta$ & $T$ & SLA & LA & W-CVaR & SLA & LA & W-CVaR \\
\midrule
\multicolumn{9}{c}{All Runs (AR)} \\
\midrule
$0.025$ & $0.01$ & $2$ & $94.43$ & $94.46$ & $93.75$ & $96.27$ & $96.17$ & $96.27$ \\
$0.025$ & $0.01$ & $4$ & $288.75$ & $280.60$ & $292.50$ & $95.88$ & $96.64$ & $95.58$ \\
$0.025$ & $0.05$ & $2$ & $93.22$ & $92.49$ & $91.60$ & $97.18$ & $97.21$ & $97.27$ \\
$0.025$ & $0.05$ & $4$ & $290.15$ & $284.63$ & $290.77$ & $96.93$ & $96.38$ & $95.23$ \\
$0.05$ & $0.01$ & $2$ & $95.94$ & $95.43$ & $94.37$ & $96.22$ & $96.35$ & $96.25$ \\
$0.05$ & $0.01$ & $4$ & $294.86$ & $287.93$ & $297.87$ & $95.94$ & $96.09$ & $96.76$ \\
$0.05$ & $0.05$ & $2$ & $94.97$ & $95.64$ & $92.27$ & $96.66$ & $96.46$ & $96.47$ \\
$0.05$ & $0.05$ & $4$ & $293.46$ & $287.03$ & $288.15$ & $96.37$ & $96.70$ & $95.96$ \\
\midrule
\multicolumn{9}{c}{Comparable Runs (CR)} \\
\midrule
$0.025$ & $0.01$ & $2$ & $94.97$ & $94.97$ & $94.39$ & $95.89$ & $95.90$ & $95.87$ \\
$0.025$ & $0.01$ & $4$ & $294.76$ & $294.76$ & $294.76$ & $95.50$ & $95.50$ & $95.50$ \\
$0.025$ & $0.05$ & $2$ & $91.43$ & $90.98$ & $90.78$ & $97.35$ & $97.33$ & $97.35$ \\
$0.025$ & $0.05$ & $4$ & $298.46$ & $297.69$ & $283.13$ & $91.15$ & $91.86$ & $91.57$ \\
$0.05$ & $0.01$ & $2$ & $95.26$ & $95.26$ & $94.70$ & $96.13$ & $96.13$ & $96.13$ \\
$0.05$ & $0.01$ & $4$ & $302.37$ & $302.43$ & $297.54$ & $96.22$ & $96.21$ & $96.50$ \\
$0.05$ & $0.05$ & $2$ & $96.60$ & $96.13$ & $93.45$ & $96.43$ & $96.42$ & $96.44$ \\
$0.05$ & $0.05$ & $4$ & $290.47$ & $282.33$ & $281.32$ & $96.94$ & $96.90$ & $96.94$ \\
\bottomrule
\end{tabular}
\label{tab:bilevel_reliability}
\end{table*}

\vfill
\end{document}